\documentclass{imanum}
\jno{drnxxx}
\received{\today}
\revised{}

\allowdisplaybreaks

\topmargin \dimexpr 0.65in -0.70in
\usepackage[usenames]{color}
\usepackage{graphicx}
\usepackage{enumerate,enumitem}
\usepackage{empheq}
\usepackage{booktabs}
\usepackage{multirow}
\usepackage{cite}
\usepackage{bm}

\usepackage{hyperref}
\hypersetup{%
colorlinks = true,
citecolor = black, 
linkcolor  = black
}

\def\x{{\bf x}}

\def\u{{\bf u}}
\def\v{{\bf v}}
\def\w{{\bf w}}
\def\e{{\bf e}}
\def\p{{\bf p}}
\def\q{{\bf q}}
\def\f{{\bf f}}

\def\M{{\bf M}}
\def\A{{\bf A}}

\def\C{{\bf C}}
\def\KK{{\bf K}}

\def\d{{\mathrm d}}
\def\D{{\mathbb D}}

\def\S{{\mathbb S}}
\def\I{{\mathbb I}}
\def\K{{\mathcal K}}

\def\CI{{\mathcal I}}

\def\R {{\mathbb R}}

\def\le{\leqslant}
\def\ge{\geqslant}
\def\Omega{\varOmega}
\def\Delta{\varDelta}

\def\E{{\mathcal E}}
\def\EE{\boldsymbol{\mathcal{E}}}

\def\Mrho{{\bf M}_{\boldsymbol{\rho}}}
\def\Brho{{\bf B}_{\boldsymbol{\rho}}}
\def\Amu{{\bf A}_{\boldsymbol{\mu}}}

\def\b{\color{blue}}

\graphicspath{{images/}}
\begin{document}

\title{Optimal convergence of the arbitrary Lagrangian--Eulerian interface tracking method for two-phase Navier--Stokes flow {without surface tension}}
\shorttitle{ALE FEM for two-phase NS flow}

\author{%
{\sc Buyang Li}\thanks{Email: buyang.li@polyu.edu.hk}\\[2pt]
Department of Applied Mathematics, 
The Hong Kong Polytechnic University, Hong Kong.
{\sc and}\\[6pt]
{\sc Shu Ma* and Weifeng Qiu}\thanks{Email:
shuma2@cityu.edu.hk and 
weifeqiu@cityu.edu.hk}\\[2pt]
Department of Mathematics, City University of Hong Kong, Kowloon, Hong Kong.
}
\shortauthorlist{B. Li, S. Ma \& W. Qiu}

\maketitle

\begin{abstract}
{Optimal-order convergence in the $H^1$ norm is proved for an arbitrary {Lagrangian--Eulerian} interface tracking finite element method for the sharp interface model of two-phase {Navier--Stokes} flow {without surface tension,} using high-order curved evolving mesh. In this method, the interfacial mesh points move with the fluid's velocity to track the sharp interface between two phases of the fluid, and the interior mesh points move according to a harmonic extension of the interface velocity. The error of the semidiscrete arbitrary {Lagrangian--Eulerian} interface tracking finite element method is shown to be $O(h^k)$ in the $L^\infty(0, T; H^1(\Omega))$ norm for the Taylor-Hood finite elements of degree $k \ge 2$. This high-order convergence is achieved by utilizing the piecewise smoothness of the solution on each subdomain occupied by one phase of the fluid, relying on a low global regularity on the entire moving domain. 
{Numerical experiments illustrate and complement the theoretical results.}
} 
{Two-phase {Navier--Stokes} flow, arbitrary {Lagrangian--Eulerian}, finite element method,
convergence, error estimates}
\end{abstract}

\setlength\abovedisplayskip{4pt}
\setlength\belowdisplayskip{4pt}

\section{Introduction}\label{Se:1}

Immiscible fluids mixture separated by a moving interface appears widely in nature and engineering applications. Current models for two-phase or multiphase flow can generally be categorized into two main types: diffuse interface models and sharp interface models. Diffuse interface models treat the interface as an extremely thin transition zone that separates the two fluids. In these models, {physical} quantities like density and viscosity change rapidly but smoothly within this transition zone. On the other hand, sharp interface models consider the interface as a surface without any thickness. In such models, {physical} quantities exhibit discontinuities across this surface.

\begin{figure}[htp]
\centering \includegraphics[width=1.7in]{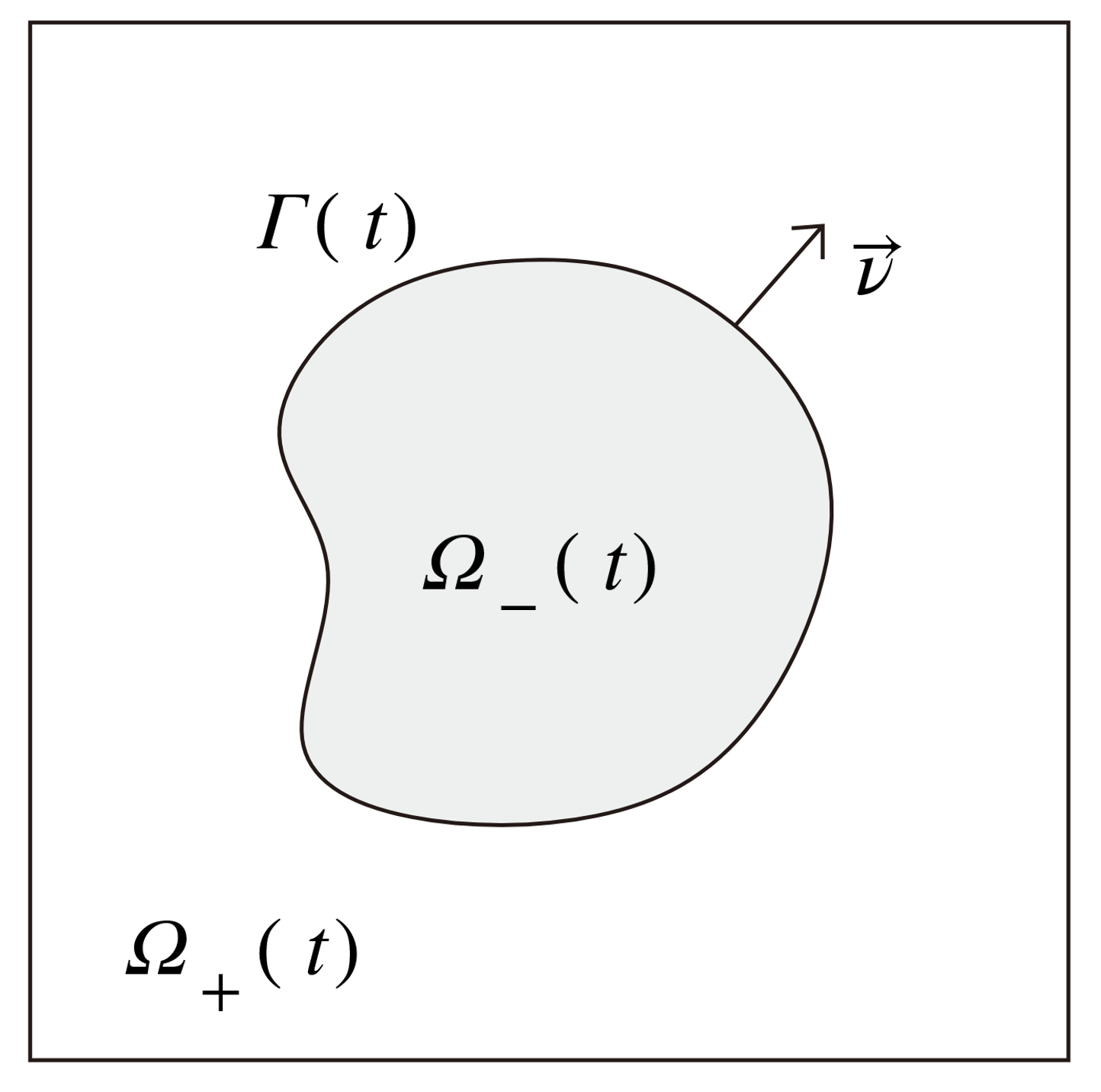}
\vspace{-8pt}
\caption{Domain $\Omega$ in the case $d = 2$.}
\label{fig_domain}
\end{figure}

We consider the sharp interface model of two-phase {Navier--Stokes (NS) flow \cite{gross2011numerical}} {without surface tension} in a bounded domain  $\Omega(t)=\Omega_+(t)\cup\Omega_-(t)\cup\Gamma(t)$ in $\R^d$, with $d\in\{2,3\}$, where $\Gamma(t)=\overline\Omega_+(t)\cap\overline\Omega_-(t)$ is a sharp interface which separates the two fluids occupying two subdomains $\Omega_+(t)$ and $\Omega_-(t)$, respectively, see Figure~\ref{fig_domain}. 
The fluid velocity $u$ and pressure $p$ are governed by NS equations in each subdomain $\Omega_\pm(t)$ and satisfy certain jump conditions on the interface $\Gamma(t)$, i.e., 
\begin{align}\label{PDE}
\left\{
\begin{aligned}
\rho_{\pm}(\partial_t u + u\cdot\nabla u) -\nabla\cdot \sigma &={f} &&\mbox{in}\,\,\,\mbox{$\bigcup\limits_{t\in[0,T]}$} \Omega_\pm(t)\times\{t\},\\
\nabla\cdot u &= 0 &&\mbox{in}\,\,\,\mbox{$\bigcup\limits_{t\in[0,T]}$} \Omega_\pm(t)\times\{t\},\\
[u]_{-}^{+}=0\quad\mbox{and}\quad
[\sigma\nu ]_{-}^{+}&= 0  &&\mbox{on}\,\,\,\mbox{$\bigcup\limits_{t\in[0,T]}$} \Gamma(t)\times\{t\},\\
u &= 0 &&\mbox{in}\,\,\,\partial\Omega\times(0,T] , \\[10pt] 
u &= u^0 &&\mbox{on}\,\,\,\Omega\times{ \{0\} ,}
\end{aligned} 
\right.
\end{align}
where {$ \rho_{\pm} $ are the densities of the fluids in $\Omega_{\pm}(t)$},  $\nu$ is the unit normal vector on $\Gamma(t)$ pointing to $\Omega_{+}(t)$; $\sigma = (2\mu_{\pm} \D(u) -p\I)$ is the stress tensor, in which $\D(u) = \frac12(\nabla u + (\nabla u)^\top)$ and $\I$ denote the deformation matrix and identity matrix, respectively; $\mu_{\pm}$ are the viscosities of the two fluids in $\Omega_{\pm}(t)$, respectively; $[u]_{-}^{+} = u_{+} - u_{-}$ and $[\sigma\nu ]_{-}^{+} = \sigma_{+}\nu - \sigma_{-}\nu$ denote the jumps of the quantities $u$ and $\sigma\nu$ across the interface, respectively. {The gravitational force is given by $f = (0,-g)^\top$ in two dimensions and $f = (0,0, -g)^\top$ in three dimensions with $g$ a constant (the gravitational acceleration).}
The interface $\Gamma(t)$ moves along with the fluid and therefore also has velocity $u$. 
If $\Gamma^0$ is the initial interface, through a flow map $X:\Gamma^0\times[0,T]\rightarrow \R^d$ we can describe the evolution of interface $\Gamma(t)=\{X(y,t):y\in\Gamma^0\}$ 
by the following equation:
\begin{align}\label{eq-interface} 
\left\{\begin{aligned}
\partial_tX(\xi,t)&= u(X(\xi,t),t) && \mbox{for}\,\,\,(\xi,t)\in \Gamma^0\times{ (0,T]},\\
X(\xi,0)&=\xi && \mbox{for}\,\,\,\xi \in \Gamma^0 . 
\end{aligned}\right. 
\end{align}
{We do not take surface tension into account on the interface $\Gamma(t)$, which differs from models that include surface tension, where the appearance of curvature leads to a coupling between the geometric properties of the interface and the fluid flow in the bulk.}
The fluid dynamics in \eqref{PDE}--\eqref{eq-interface} involve the movement of the interface $\Gamma(t)$ and two subdomains $\Omega_{\pm}(t)$ over time, while the overall domain $\Omega = \Omega(t)$ remains unchanged in shape.

Numerical methods for \eqref{PDE} can generally be classified into three categories, i.e., the Eulerian approach, the Lagrangian approach, and the arbitrary {Lagrangian--Eulerian} approach (ALE).
The arbitrary {Lagrangian--Eulerian} approach serves as a bridge between the Lagrangian and Eulerian approaches, by enabling the frame to move with an ``arbitrary" bulk velocity that fits the interface motion. At the discrete level, the mesh points can be displaced independently of the flow.  
Utilizing the ALE technique allows for flexible movement of the inner domain mesh, while the mesh on the interface can move alongside materials to precisely track the interfaces of a multi-material system.

In an early investigation of ALE methods, {the} stability analysis of {the} ALE finite element method (FEM) was firstly conducted by Formaggia and {Nobile} \cite{formaggia1999stability} for a linear parabolic equation. 
Subsequently, Gastaldi \cite{gastaldi2001priori}  developed a priori error estimates in the $L^2$ norm for ALE-FEM when solving parabolic equations. 
Many studies on the numerical analysis of ALE-FEMs for parabolic moving boundary/interface problems can be found in works such as\cite{nobile2001numerical, boffi2004stability, formaggia2004stability, badia2006analysis, lan2020novel, kovacs2018higher}.
Optimal-order error bounds of $O(h^{k+1})$ in the $L^\infty(0, T; L^2)$ norm of ALE semidiscrete FEM with curved triangles/simplices for parabolic moving boundary problems were established by Elliott and Ranner \cite{elliott2021unified}.
We also refer to \cite{Edelmann-2021} for an ALE method with harmonically evolving mesh. Optimal-order $H^1$ convergence of the ALE-FEM for PDEs coupling boundary evolution arising from shape optimization problems was proved in \cite{gong2023convergent}. These results were established for high-order curved evolving mesh. Optimal convergence of $O(h^{k+1})$ in the $L^\infty(0, T; L^2)$ norm, with flat evolving simplices in the interior and curved simplices exclusively on the boundary, was prove in \cite{li2023optimal} for the ALE semidiscrete FEM utilizing a standard iso-parametric element of degree $k$ in \cite{Lenoir-1986}.
The ALE methods for PDEs in bulk domains \cite{gong2023convergent} are also closely related to the evolving FEMs for PDEs on evolving surfaces, for which optimal-order convergence in the $L^2$ and $H^1$ norms were shown in \cite{dziuk2007finite, elliott2015error, KoLiESFEM2017}.
The above-mentioned research efforts have focused on diffusion equations.

The numerical analysis of ALE methods for the Stokes and {Navier--Stokes} equations has also noteworthy results but remained suboptimal.
Legendre and Takahashi \cite{legendre2008convergence} established an $L^2$ error estimate of $O(\tau + h^{\frac12})$ for the method of characteristics combined with a $P_{1b}-P_1$ finite element approximation to the  ALE formulation of 2D {Navier--Stokes} equations.
In \cite{san2009convergence}, by using the Taylor-Hood $P_2$-$P_1$ element, the error estimates of $O(h^2|\log h|)$ and $O(\tau + h^2+ h^2/\tau)$ were proved for the semi-discrete and fully discrete ALE-FEM for the Stokes equations in a time-dependent domain, respectively.
The error of ALE finite element solutions to the Stokes equations on a time-varying domain, with BDF-$r$ in time (for $1\le r\le 5$) and the Taylor-Hood $P_k$--$P_{k-1}$ elements in space (with degree $k\ge 2$), was proved to be $O(\tau^r + h^k)$ in the $L^2$ norm in \cite{liu2013simple}. 
Optimal-order error bounds of $O(h^{k+1})$ and $O(\tau^2+h^{k+1})$ of the semidiscrete and fully discrete ALE-FEMs in the $L^2$ norm for the Stokes equations were proved recently in \cite{Rao-Wang-Xie-arXiv}.

Many efforts have also been made in developing efficient ALE-FEMs for two-phase NS flow problems, {primarily for models with surface tension}. 
For example, Anjos et al. \cite{anjos20143d} introduced a 3D ALE-FEM scheme with adaptive mesh to achieve good control of the mesh quality by using repeated interpolations of the velocity.
Barrett et al. \cite{barrett2015on} { also considered surface tension and} designed a stable ALE-FEM in which the mesh points used to describe the interface are not mesh points of the underlying bulk finite element mesh. An implicit tangential velocity of the interface nodes is introduced in \cite{barrett2015on} to acquire good mesh quality of the interface.
{Barrett, Garcke and N\"urnberg introduced a novel weak
formulation \cite{barrett2020parametric}, which is referred to as the BGN formulation from now on.} The BGN formulation allows for tangential degrees of freedom to improve the mesh quality.
This formulation was employed for two-phase Stokes flow \cite{barrett2013eliminating} and two-phase NS flow \cite{barrett2015stable}, {both with surface tension and} with unfitted finite element approximations, which leads to an unconditional stability bound.  
Moreover, the applications of BGN formulation to the two-phase NS flow {with surface tension} with moving fitted finite element methods were also investigated in \cite{agnese2020fitted} in both Eulerian and ALE approaches, but no stability bound was available. 
Fu \cite{fu2020arbitrary} proposed an ALE-HDG FEM for two-phase NS flow {with surface tension} with an implicit tangential velocity to acquire good mesh quality of the interface, by using HDG to produce an exactly divergence-free velocity approximation. 
Recently,  Duan et al. \cite{duan2022energy} proposed an energy-diminishing ALE scheme  {with surface tension} by using Stokes extension of the velocity from the interface to the bulk domain, in order to produce an exactly divergence-free mesh velocity.
Energy stability of the ALE-FEMs {with surface tension} was discussed in \cite{barrett2015on, barrett2013eliminating, barrett2015stable, fu2020arbitrary, duan2022energy}.
We also refer to \cite{garcke2023structure} for the structure-preserving ALE-FEMs for the two-phase NS flow {with surface tension} in both the conservative and non-conservative forms. The methods introduced in \cite{garcke2023structure} were shown to satisfy unconditional stability and exact volume preservation on the fully discrete level.
The numerical investigation of convergence behavior (without error analysis) {with surface tension}, can be found in \cite{duan2022energy, fu2020arbitrary} .
{As far as we know, there is no error analysis of ALE-FEMs for the sharp interface model of two-phase NS flow problems with surface tension.
In this work, we focus on a simplified model of two-phase NS flow without surface tension.}
Two major difficulties are discussed below.

First, the error analysis of ALE-FEMs for parabolic and Stokes/NS equations in the above-mentioned articles all require the triangulation of the domain to fit the moving boundary or interface exactly, i.e., the triangulated domain $\Omega_h(t)$ for {computation  that interpolates} the exact domain $\Omega(t)$. This is the case when the boundary or interface of domain $\Omega(t)$ is given, determined either by a given flow map or by a given velocity field. However, this is the not the case in two-phase NS flow where the velocity of the interface is unknown (determined by the solution of the PDE problem). Thus the location of the exact interface is unknown {and needs to} be approximated based on the numerical solution of the two-phase NS flow. As a result, the numerically computed interface at $t>0$ generally has a distance from (therefore is not an interpolation of) the exact interface, and the triangulated {subdomains} $\Omega_{h,\pm}(t)$ are not interpolations of the exact subdomains $\Omega_{\pm}(t)$ at $t>0$. Consequently, the error analysis of ALE-FEMs in the pre-existing articles for parabolic and Stokes/NS equations with a given moving boundary/interface cannot be applied to the two-phase NS flow {without surface tension}.

Second, the error analysis of ALE-FEMs for parabolic and Stokes/NS equations in the above-mentioned articles all require the solution to be sufficiently smooth in the whole domain, while the solution of the two-phase NS equation in \eqref{PDE} is generally not smooth in the whole domain due to the existence {of the free interface} which separates the two phases of the fluid, e.g., $\nabla u$ is generally discontinuous {across the interface} $\Gamma(t)$.
Since the numerical interface and {the} exact interface differ from each other by {a} certain distance at $t>0$, the triangles of the mesh (which fit the numerical interface) generally do not fit the exact interface $\Gamma(t)$. Thus the exact solution $u$ may be nonsmooth inside a triangle of the mesh. This situation {renders it difficult to prove} high-order convergence using standard interpolation onto the mesh.

The aim of this article is to establish an optimal-order error bound in the $L^\infty(0, T; H^1(\Omega))$ norm for the semidiscrete ALE-FEM for the two-phase NS flow problem {without surface tension} in \eqref{PDE} {using high-order curved evolving mesh}, where the mesh points on the numerical interface {move} with the fluid's velocity {and is harmonically extended} to the bulk subdomains. 
We address the two major difficulties mentioned above by using the matrix-vector formulations of the ALE-FEM and the error equations. Such techniques were originally introduced in \cite{KoLiESFEM2017} for analyzing solution-driven surface evolutions. It was later developed for analyzing the convergence of evolving FEMs for surface evolution under geometric flows \cite{Bai-Li-FoCM,Hu-Li-2022,KLL-MCF,KLL-Willmore} and domain evolution under shape gradient flows \cite{gong2023convergent}. We adapt such matrix-vector formulation techniques here to address the mismatch between {the} numerical interface and {the} exact interface, using {an} intermediate auxiliary triangulated domain, as well as the missing of global regularity of the solution by using lift operations which map the numerical subdomains to the exact subdomains. This approach only requires using piecewise smoothness of the solution $u$ on the subdomains $\Omega_\pm(t)$ instead of its global smoothness of $u$ on the whole domain $\Omega(t)$. 
For a smoothly evolving sharp interface $\Gamma(t)$ without topological changes,  and a solution $(u,p)\in C([0,T];H^1(\Omega)^d)\times L^2(0,T;L^2(\Omega))$ which is piecewise smooth in the subdomains $\Omega_\pm(t)$ for $t\in[0,T]$, we prove optimal-order error bound of $O(h^k)$ in the $C([0,T];H^1(\Omega)^d)\times L^2(0,T;L^2(\Omega))$ norm for the semi-discrete ALE-FEMs scheme with Taylor-Hood $P_k$--$P_{k-1}$ elements of degree $k \ge 2$; see Theorem \ref{Main-THM}.

The rest of this article is organized as follows. Section~\ref{section:results} introduces the ALE reformulation of \eqref{PDE} and the ALE weak formulation on the evolving domain.
In Section \ref{section:consistency} we introduce the ALE weak formulation and ALE evolving FEM on the interpolated evolving domain {and provide a consistency analysis}. In Section \ref{section:consistency} we formulate the consistency error estimates based on interpolated evolving mesh and matrix-vector formulation, which serve not only for implementation purposes but also play a key role in the subsequent convergence analysis. In Section \ref{section:proof} we prove optimal-order convergence of the ALE evolving FEM considered in this article. Numerical experiments are presented in Section \ref{section:numerical-test} to illustrate the convergence of the ALE evolving FEM and its performance on a benchmark example.

\section{The ALE evolving finite element method}\label{section:results}

\subsection{The ALE reformulation of PDE} \label{section:2.1}
Let $w(\cdot,t): \Omega(t) \rightarrow \R^d$ be the harmonic extension of fluid's velocity $u$ from the interface $\Gamma(t)$ to {the} bulk subdomains $\Omega_\pm(t)$, i.e.,
\begin{align}\label{laplace-W}
- \Delta w =0 \,\,\,\mbox{in}\,\,\, \Omega_\pm(t) \quad \mbox{with} \quad
w=u  \,\,\,\mbox{on}\,\,\,\Gamma(t)\quad \mbox{and} \quad
w = 0 \,\,\,\mbox{on}\,\,\,\partial\Omega(t).
\end{align}
Harmonic extension is a common method for generating mesh velocity $w$ in the bulk subdomains; see \cite{boffi2004stability, Edelmann-2021, fu2020arbitrary}. Since $w=u$ on $\Gamma(t)$, it follows that $w$ and $u$ generate the same evolution of interface $\Gamma(t)$ and subdomains $\Omega_\pm(t)$.

{Denote the initial domain by $ \Omega^0 = \Omega(0) $ and the initial subdomains by $ \Omega_\pm^0 = \Omega_\pm(0) $.}
Let  $\phi:\Omega^0\times[0, T]\rightarrow \R^d$ be the flow map generated {by the velocity field} $w$. Then $\phi$ maps $\Omega_{\pm}^0$ to $\Omega_\pm(t)$ and $\Gamma^0$ to $\Gamma(t)$, respectively, with  
$$
\Omega(t)  = \{ \phi(\xi,t) : \xi \in \Omega^0\} \quad \mbox{with} \quad\phi(\xi,0)=\xi\quad  \mbox{for} \quad \xi \in \Omega^0 , 
$$
and satisfies the following ordinary differential equation (ODE): 
\begin{align}\label{ode-phi}
\partial_t \phi(\xi,t) = w(\phi(\xi,t),t) \quad \mbox{with} \quad\phi(\xi,0)=\xi\quad  \mbox{for} \quad \xi \in \Omega^0.
\end{align}
%
We consider an ALE moving frame associated with the flow map $\phi$, with frame velocity $w$. 
Then the original problem in \eqref{PDE} can be written as 
\begin{equation}\label{PDE-ALE}
\left \{
\begin{aligned} 
\rho_{\pm}\partial^\bullet_t u +  \rho_{\pm}((u-w)\cdot\nabla u) -\nabla\cdot \sigma &=0 &&\mbox{in}\,\,\, \Omega_\pm(t),\\
\nabla\cdot u &= 0 &&\mbox{in}\,\,\, \Omega_\pm(t),
\end{aligned}
\right . 
\end{equation}
where $\partial_t^\bullet u = \partial_t u + w \cdot \nabla u$ is the material derivative with respect to the frame velocity $w$. The subdomains $\Omega_\pm(t)$ evolve with velocity $w$.  


We assume that problem \eqref{PDE}, or equivalently \eqref{PDE-ALE}, has a solution 
{$$(u, p)\in C([0,T];H^1_0(\Omega)^d)\times C([0,T];L^2_0(\Omega)),$$}
which is piecewise smooth in subdomains $\Omega_\pm(t)$ for $t \in (0, T]$, where
\begin{align*}
H_0^1(\Omega)^d  = \{ v \in H^1(\Omega)^d : v = 0 \,\, \,\mbox{on} \, \, \,\partial \Omega\} \quad \mbox{and} \quad 
L_0^2(\Omega)  = \{ q \in L^2(\Omega) \,\,\, \mbox{and} \, \int_\Omega q \, \d x = 0 \}.
\end{align*} 
By piecewise smooth we mean that $\sum_{j=0}^2\|(\partial_t^\bullet)^j u\|_{H^{k+1}(\Omega_\pm(t))}+\sum_{j=0}^1\|(\partial_t^\bullet)^j  p \|_{H^{k}(\Omega_\pm(t))}$ is bounded uniformly with respect to $t\in[0,T]$.

We assume that the flow map $\phi(\cdot, t) \in C^{k}(\bar \Omega_{\pm}^0)$ is piecewise smooth in the subdomains $\Omega_{\pm}^0$. By pulling back to the initial domain $\Omega^0$ after composition with $\phi(\cdot, t)$ and utilizing the regularity estimate of the elliptic equation \cite[Theorem 5, Section 6.3]{evans2022partial} on the fixed domain $\Omega^0$, the frame velocity $w$ defined in \eqref{laplace-W} satisfies
$$\|w(t)\|_{H^{k}(\Omega_\pm(t))}  \le C\|u(t)\|_{H^k(\Omega_\pm(t))},$$
which is also piecewise smooth. Then, by using notations  
\begin{align*}
(u, v) = \int_{\Omega_+(t)} uv\, \d x + \int_{\Omega_-(t)}uv \, \d x,
\quad
\rho = 
\left \{
\begin{aligned}
\rho_{+} \quad \mbox{in} \,\, \Omega_{+}(t),\\
\rho_{-} \quad \mbox{in} \,\, \Omega_{-}(t),
\end{aligned}
\right .  \quad \mbox{and} \quad
\mu = 
\left \{
\begin{aligned}
\mu_{+} \quad \mbox{in} \,\, \Omega_{+}(t),\\
\mu_{-} \quad \mbox{in} \,\, \Omega_{-}(t) , 
\end{aligned}
\right .
\end{align*} 
testing \eqref{PDE-ALE} by $v$ and using the interface and boundary conditions in \eqref{PDE}, we obtain the following weak reformulation. 

{\bf ALE weak formulation}: Find $(u, p) \in C([0, T]; H_0^1(\Omega)^d) \times C([0, T]; L_0^2(\Omega)) $  with $\partial^\bullet_t u \in C([0, T]; L^2(\Omega)^d) $,  $w\in C([0,T]; H^1(\Omega_\pm)^d)$, and $\Omega_\pm(t)=\phi(\Omega_\pm^0,t)$,
such that 
\begin{align}\label{PDE-ALE-weak}
\left\{ \,\,
\begin{aligned}
&(\rho \partial^\bullet_t u, v)  + (\rho (u- w)\cdot\nabla u, v) 
+ 2(\mu \D(u), \D(v)) - (p, \nabla \cdot v)  
= {(f, v)}   && \quad \forall \,\, v \in H_0^1(\Omega)^d\\[4pt]
& (\nabla\cdot u, q) = 0 && \quad\forall \,\, q \in L_0^2(\Omega) , \\[4pt]
& (\nabla w, \nabla \chi)  = 0  \quad\forall \,\,  \chi \in H_0^1(\Omega_\pm)^d
\quad\mbox{and}\quad w= u  \quad\mbox{on}\,\,\,\Gamma(t) ,\\
& \partial_t \phi   = w \circ\phi \quad\mbox{on}\,\,\, \Omega^0 ,
\end{aligned}
\right .
\end{align}
with initial conditions $\phi(0)={\rm id}|_{\Omega^0}$ and $u(0)=u^0$.

\subsection{The ALE mapping and ALE evolving FEM}
Assume that the given smooth initial domain $\Omega^0\subset\R^d$, $d\in\{2,3\}$, is divided into a set $\K_h^0$ of shape-regular and quasi-uniform simplices with mesh size $h$, using curved finite elements of degree $k$, and fitting the interface $\Gamma^0 = \Gamma(0)$.
Every possibly curved closed simplex $K \in \K_h^0$ is the image of a unique polynomial of degree $k$ defined on the reference simplex $\hat K$, denoted by $F_{K}:\hat K\rightarrow K$, called the parametrization of $K$; see \cite[\textsection 8.6]{elliott2021unified}. Every boundary simplex $K\in\K_h^0$ (with one face or edge attached to $\Gamma^0$) contains a possibly curved face or edge to interpolate $\Gamma^0$ with accuracy of $O(h^{k+1})$. This can be obtained by using Lenoir's isoparametric finite elements \cite{Lenoir-1986}.
Thus the initial domain $\Omega^0$ is approximated by the triangulated domain 
$$
\Omega_h^0 = {\Big(\bigcup_{K\in \K_h^0}K \Big) \Big\backslash\partial \Omega_h^0} = \Omega_{h, +}^0 \cup\Omega_{h, -}^0\cup\Gamma_h^0.
$$

Let $\x^0=(\xi_1,\cdots,\xi_M)\in \R^{dM}$ be the vector that collects all the nodes $\xi_j\in \R^d$, $j=1,\ldots,M$, (which define finite elements of degree $k$) in the triangulation $\K_h^0$.      
We evolve the vector $\x^0$ in time and denote its position at time $t$ by 
$$\x(t)=(x_1(t),\cdots,x_M(t)),$$
which determines the triangulation $\K_h[\x(t)]$ and domain $\Omega_{h}[\x(t)]$, as well as subdomains $\Omega_{h,\pm}[\x(t)]$ and interface $\Gamma_h[\x(t)] $ such that 
$$
\Omega_{h}[\x(t)] = \Big( \bigcup_{K\in \K_h[\x(t)]}K \Big) \Big\backslash\partial \Omega_{h}[\x(t)]  = \Omega_{h,+}[\x(t)] \cup \Omega_{h,-}[\x(t)] \cup\Gamma_h[\x(t)] .
$$  
Similarly as the simplices on the initial domain, if $K \in \K_h[\x(t)]$ is a simplex on the evolving domain $\Omega_{h}[\x(t)] $, then we denote by $F_{K}:\hat K\rightarrow K$ the parametrization of $K$. 
For the simplicity of notation, we denote  
$$
\Omega_{h}(t) = \Omega_{h}[\x(t)],\quad \Omega_{h,\pm}(t) = \Omega_{h,\pm}[\x(t)] 
  \quad\mbox{and}\quad
\Gamma_{h}(t) = \Gamma_{h}[\x(t)]  . 
$$
Correspondingly, the scalar-valued finite element space on the evolving domain $\Omega_h(t)$ is defined as 
\begin{align*}
S_h^k[\x(t)] :=& \{v_h\in H^1(\Omega_h(t)): v_h\circ F_K \in P^k(\hat K)\,\,\mbox{for all}\,\, K\in \K_h[\x(t)]\} \\
{S_{h,\rm dc}^{k-1}[\x(t)] :=}& 
{\{v_h\in H^1(\Omega_{h,\pm}(t)): v_h\circ F_K \in P^{k-1}(\hat K)\,\,\mbox{for all}\,\, K\in \K_h[\x(t)]\} . }
\end{align*}

The ALE mapping $\phi_h(\cdot, t)$ is defined as the unique finite element function $\phi_h(\cdot,t) \in  S_h^k[\x^0]^d$ such that
$$\phi_h(\xi_j,t)=x_j(t) \quad \mbox{for} \quad j=1,\dots,M,$$
which maps the initial domain $\Omega_{h, \pm}^0$ onto the evolving domain $\Omega_{h, \pm}(t)$, and the initial interface $\Gamma_h^0$ onto the evolving interface $\Gamma_h(t)$, i.e., 
$$
\Omega_{h, \pm}(t) = \phi_h(\Omega_{h, \pm}^0, t) \quad \mbox{and} \quad 
\Gamma_h(t) = \phi_h(\Gamma_h^0, t) \quad \mbox{for} \quad t \in [0, T].
$$
The ALE mapping $\phi_h(\cdot, t)$ is the flow map generated by a mesh velocity $w_h(\cdot, t)$ and satisfies the following relation:
\begin{align}
\frac{\d}{\d t}\phi_h(\xi,t) = w_h(t) \circ \phi_h(\xi,t) \quad \mbox{for} \quad \xi\in \Omega_h^0 . 
\end{align}
Based on the ALE mapping $\phi_h(\cdot, t)$ and the mesh velocity $w_h(\cdot, t)$, 
if $v(\cdot, t)$ is a function defined on $\Omega_h(t)$, its discrete material derivative is defined as
$$
\partial_{t,h}^{\bullet} v(x,t)
= \frac{\d}{\d t} v(\phi_h(\xi,t),t) \quad\mbox{for}\,\,\, x=\phi_h(\xi,t) \in \Omega_h(t) \quad \mbox{for} \,\,\,\xi\in \Omega_h^0.
$$





Let $X_h[\x(t)] \times  M_h[\x(t)] \subset H_0^1(\Omega_h(t))^d \times L_0^2(\Omega_h(t)) $ denote 
the Taylor--Hood finite element space of degree $k\ge 2$ subject to the triangulation $\K_h[\x(t)]$ of evolving domain $\Omega_h(t)$, i.e.,
\begin{align*}
X_h[\x(t)] &= S_h^k[\x(t)]^d \cap H_0^1(\Omega_h(t))^d 
\quad\mbox{and}\quad 
M_h[\x(t)] = {S_{h,\rm dc}^{k-1}[\x(t)] \cap  L_0^2(\Omega_h(t)) }. 
\end{align*} 
We also denote by $\mathring X_{h}[\x(t)] $ the subspace of $X_h[\x(t)]$ with zero interface {and boundary conditions}, i.e.,
$$
\mathring X_{h}[\x(t)]  = \{ v_h\in X_h[\x(t)] : v_h = 0\,\,\,\mbox{on}\,\,\, \Gamma_h(t)=\Gamma_h[\x(t)]  \,\,\, {\mbox{and on}\,\,\,  \partial \Omega_{h}(t) = \partial  \Omega_{h}[\x(t)] }\} . 
$$
It is known that the Taylor--Hood finite element space satisfies the discrete inf-sup condition {(see \cite[\textsection 8.5.3]{BBF-2013} on the domain decomposition techniques)}, i.e., 
there exists positive constant $\kappa$ independent of $h$ such that 
\begin{align}\label{inf-sup} 
\kappa\|q_h\|_{L^2(\Omega_h(t))} \le  \sup_{0 \neq v_h \in X_h[\x(t)]} \frac{(q_h, \nabla \cdot v_h)_{\Omega_h}}{\|v_h\|_{H^1(\Omega_h(t))}} 
\quad\forall\, q_h\in M_h[\x(t)], 
\end{align}
where 
\begin{align*}
(u_h, v_h)_{\Omega_h} := \int_{\Omega_{h,+}(t)} u_h v_h\, \d x + \int_{\Omega_{h, -}(t)} u_h v_h\, \d x . 
\end{align*}

{Since the mesh velocity} $w_h$ should be an approximation of $w$ defined in Section \ref{section:2.1}, 
we define $w_h(\cdot, t)$ to be a discrete harmonic finite element function in the numerically computed {subdomains} $\Omega_{h, \pm}(t)$ with boundary condition matching the velocity of the fluid on {the boundary $\partial\Omega_h(t)$} and the numerically computed interface $\Gamma_h(t)$. 
This is equivalent to finding $w_{h}\in u_{h} + \mathring X_{h}[\x(t)]$ such that 
\begin{align}
\int_{\Omega_{h}(t)} \nabla w_{h} \cdot \nabla \chi_{h} \, \d x  = 
0
\quad \mbox{for} \quad \chi_{h} \in \mathring X_{h}[\x(t)] . 
\end{align}

We denote by $\rho_h$ and $\mu_h$ the piecewise constant functions in $\Omega_h(t)$, defined by 
\begin{align*}
\rho_h = 
\left \{
\begin{aligned}
\rho_{+} \quad \mbox{in} \,\, \Omega_{h, +}(t)\\
\rho_{-} \quad \mbox{in} \,\, \Omega_{h, -}(t)
\end{aligned}
\right .  \quad \mbox{and} \quad
\mu_h = 
\left \{
\begin{aligned}
&\mu_{+} \quad \mbox{in} \,\, \Omega_{h, +}(t)\\
&\mu_{-} \quad \mbox{in} \,\, \Omega_{h, -}(t).
\end{aligned}
\right .
\end{align*}
The ALE evolving FEM for \eqref{PDE-ALE-weak} is defined as follows. 
%

{\bf ALE evolving FEM}: Find a flow map $\phi_h(\cdot,t) \in S_h^k[\x^0]^d $ which determines nodal vector $\x(t) = \phi_h(\x^0, t)$ and corresponding subdomains $\Omega_{h,\pm}(t)=\Omega_{h,\pm}[\x(t)]$, as well as
$(u_h(\cdot,t) , p_h(\cdot,t) ) \in X_h[\x(t)] \times M_h[\x(t)]$ and $w_{h}(\cdot,t) \in u_h(\cdot,t) + \mathring X_h[\x(t)]$ such that 
\begin{align} \label{fem-ALE-weak}
\left\{ \,\,
\begin{aligned}
&(\rho_h \partial^\bullet_{t, h} u_h, v_h)_{\Omega_h} 
+ 2(\mu_h\D(u_h), \D(v_h))_{\Omega_h} - (p_h, \nabla \cdot v_h)_{\Omega_h}  \\[2pt]
& \hspace{100pt} + (\rho_h (u_h-w_h)\cdot\nabla u_h, v_h)_{\Omega_h}  = {(f_h, v_h)_{\Omega_h} }  \qquad &&  \forall\, v_h \in X_h[\x(t)],\\[4pt]
& (\nabla\cdot u_h, q_h)_{\Omega_h}  = 0 && \forall \, q_h \in M_h[\x(t)],\\[4pt]
& (\nabla w_h, \nabla \chi_h)_{\Omega_{h}}  = 0   && \forall \,  \chi_h \in \mathring X_{h}[\x(t)] , \\
& \partial_t \phi_h   = w_h \circ\phi_h , 
\end{aligned} 
\right . 
\end{align} 
with initial conditions $\phi_h(0) = {\rm id}|_{\Omega_h^0}$ and $u_h(0)=u_h^0$, where $u_h^0$ is the Lagrange interpolation of $u^0$ {and $ f_h $ is the Lagrange interpolation of $ f $ on $\Omega_h(t) $.}

\begin{remark}
The finite element basis functions of  $S_h^k[\x(t)]$ {are} denoted by 
$\varphi_j[\x(t)]: \Omega_h(t) \to \R$ for $j=1,\dots, M$, 
which have the property that on every triangle their pullback to the reference triangle
is polynomial of degree $k$, satisfying the identities
\begin{equation*}
\varphi_j[\x(t)](x_i(t)) = \delta_{ij} \quad\mbox{for all}\quad i,j = 1,  \dotsc, M .
\end{equation*}
The pullback of $\varphi_j[\x(t)]$ from $\Omega_h(t)$ to $\Omega_h^0$ is $\varphi_j[\x(t)]\circ \phi_h(\cdot,t)=\varphi_j[\x^0]$, which is simply the finite element basis functions on $\Omega_h^0$. Therefore, the basis functions $\varphi_j[\x(t)]$ satisfies the transport property \cite{dziuk2007finite}:
\begin{align}\label{transport-phi}
\partial_{t,h}^{\bullet} \varphi_j[\x(t)](x) = 0 \quad\mbox{for}\,\,\, x\in\Omega_h(t),\,\,\, j=1,\dots,{M.}
\end{align}
By using the basis functions $\varphi_j[\x(t)]$, the finite element space $S_h^k[\x(t)]$ and the discrete material derivative of the finite element solutions $v_h(\cdot, t) \in S_h^k[\x(t)]$ can be written as  
\begin{align*}
S_h^k[\x(t)] &=\Big\{ \sum_{j=1}^M c_j \varphi_j[\x(t)]: c_j\in \mathbb{R} \Big\},\\
\partial_{t,h}^{\bullet} v_h(x,t)
& =  \sum_{j=1}^M \dot v_j(t) \varphi_j[\x(t)](x) \quad\mbox{with} \quad
v_h(x,t)
=  \sum_{j=1}^M v_j(t) \varphi_j[\x(t)](x) \quad\mbox{for}
\,\,\, x\in \Omega_h(t) , 
\end{align*}
where the dot denotes the time derivative with respect to $t$. 
\end{remark}
{
\begin{remark}
Due to the use of higher-order isoparametric elements, the integrals in (2.8) cannot be evaluated exactly in practice and are instead approximated using numerical quadrature. The quadrature rules are chosen based on the geometry and polynomial order of the elements to ensure high accuracy.  In our analysis, quadrature errors are not considered, as they are typically negligible compared to the discretization error.
\end{remark}
}

\section{Interpolated evolving mesh and consistency analysis}\label{section:consistency}
In the error analysis, we need to compare functions on two different domains, i.e., the exact domain $\Omega(t)$, and the triangulated domain $\Omega_h(t)$ determined by the numerical solution. 
The consistency error is the defect upon inserting the Lagrange interpolation of the exact solution $u$ and the interpolated interface into the discretized equations. To this end, we define the interpolated ALE mapping and interpolated evolving mesh below. 

\subsection{Interpolated ALE mapping and interpolated evolving mesh}\label{subsec:lift}
We evolve the vector $\x^0=(\xi_1,\cdots,\xi_M)$ in time and denote $\x^*(t)=(x_1^*(t),\dots,x_M^*(t))$ the nodal vector which collects the points $x_j^*(t) = \phi(\xi_j, t)$ on the exact domain $\Omega(t)$, which determines triangulation $\K_h[\x^*(t)]$ and interpolated domain $\Omega_h^*(t)$. 
Since the initial triangulation is shape-regular and quasi-uniform, the interpolated triangulations $\K_h[\x^*(t)]$ will keep shape-regular and quasi-uniform for $t\in[0,T]$ when the flow map $\phi(\cdot,t):\Omega_\pm(0)\rightarrow\Omega_\pm(t)$ and its inverse are piecewise smooth on $\Omega_\pm(0)$ and $\Omega_\pm(t)$ for $t\in[0,T]$, respectively. 

Subject to the triangulation $\K_h[\x^*(t)]$ of domain $\Omega_h^*(t)$, the space of scalar-valued finite element functions on $\Omega_h^*(t)$ is defined as 
$$
S_h^k[\x^*(t)] := \{v_h\in H^1(\Omega_h^*(t)): v_h\circ F_K \in P^k(\hat K)\,\,\mbox{for all}\,\, K\in \K_h[\x^*(t)]\} . 
$$

We define the interpolated ALE mapping  $\phi_h^*(\cdot,t):\Omega_h^0 \to \Omega_h^*(t)$ to be the Lagrange interpolation of flow map $\phi(\cdot, t): \Omega(0) \to \Omega(t)$, which has the following properties:
$$
\phi_h^*(\xi_j,t)=x_j^*(t) = \phi(\xi_j, t) \quad\mbox{for}\quad  j=1,\dots,M,
$$
and 
$$
\Omega_{h, \pm}^*(t) = \phi_h^*(\Omega_{h, \pm}^0, t) \quad \mbox{and} \quad 
\Gamma_h^*(t) = \phi_h^*(\Gamma_h^0, t) . 
$$
If we denote by $w_h^*(\cdot, t)= I_h^*(t) w(\cdot, t)$ the Lagrange interpolations {of the exact velocity} $w(\cdot, t)$, then the following relation holds:
\begin{align}
\frac{\d}{\d t}\phi_h^*(\xi,t) = w_h^*(\phi_h^*(\xi,t), t) \quad \mbox{for} \quad \xi\in \Omega_h^0 . 
\end{align}
This equality holds because both sides are finite element functions on $\Omega_h^0$ and they are equal on the nodes. 

Let $u_h^*$ be the Lagrange interpolation {of the exact solution} $u$, and denote by $w_{h,\pm}^*$ and $u_{h,\pm}^*$ the restriction of $w_{h}^*$ and $u_{h}^*$ to {subdomains} $\Omega_{h,\pm}^*(t)$. Then the interpolated mesh velocity satisfies $w_h^*=u_h^*$ on the interpolated interface $\Gamma_h^*(t)$. Thus $w_{h,\pm}^*-u_{h,\pm}^*\in \mathring X_{h, \pm}[\x^*(t)] $. 

For a function $v(\cdot,t)$ defined on $\Omega_h^*(t)$, we can define the discrete material derivative $\partial_{t,*}^{\bullet}$ by using the interpolated ALE mapping $\phi_h^*(\cdot,t)$ and the interpolated mesh velocity $w_h^*(\cdot, t)$ as 
\begin{align}\label{dis-mat-deri-h*}
\partial_{t,*}^{\bullet} v(x,t) 
= \frac{\d}{\d t} v(\phi_h^*(\xi,t),t)
\quad\mbox{at point}\quad x=\phi_h^*(\xi,t) \in \Omega_h^*(t).
\end{align}

\subsection{Lifting operator}
The transition between the interpolated domain $\Omega_h^*(t)$ and the exact domain $ \Omega_h(t)$ is done by a lift operator. In Lenoir's isoparametric finite element approximation to $\Omega^0$, there exists a map $\Phi:\Omega_h^0\rightarrow\Omega^0$ such that  (cf. \cite{Lenoir-1986}): 
\begin{align*}
\|{\Phi(\cdot)} - {\rm id} \|_{L^\infty(\Omega_h^0)} + h\|\Phi(\cdot) - {\rm id}\|_{W^{1, \infty}(\Omega_h^0)} \le Ch^{k+1} \quad \mbox{in} \quad \Omega_h^0.
\end{align*}
For the triangulation $\K_h[\x^*(t)]$ which fits the interface $\Gamma(t)$, there exists a continuous mapping $\Phi_h^*(\cdot,t): \Omega_h^*(t)\rightarrow \Omega(t)$ given by
\begin{align}
\Phi_h^*(\cdot,t)=\phi(\cdot,t)\circ\Phi\circ \phi_h^*(\cdot,t)^{-1} \quad \mbox{in} \quad \Omega_h^*(t),
\end{align}
that maps $\Gamma_h^*(t)$ to $\Gamma(t)$ and satisfies the following estimates
\begin{align}\label{lift-approx}
\|\Phi_h^*(\cdot,t) - {\rm id} \|_{L^\infty(\Omega_h^*(t))} + h\|\Phi_h^*(\cdot,t) - {\rm id}\|_{L^\infty(\Omega_h^*(t))} \le Ch^{k+1} \quad \mbox{in} \quad \Omega_h^*(t).
\end{align}
Correspondingly, for $v_h: \Omega_h^*(t) \to \R$, we denote its lift by $v_h^{\ell}$ given by 
\begin{align}\label{lift}
v_h^{\ell} = v_h \circ \Phi_h^*(\cdot,t)^{-1}\quad \mbox{in} \quad \Omega(t).
\end{align}


\begin{enumerate}[label={\rm(\arabic*)},ref=\arabic*,topsep=2pt,itemsep=0pt,partopsep=1pt,parsep=1ex,leftmargin=10pt]

\item[$\bullet$] (Lagrange interpolation) Define $\hat I_{h, \hat K}$ be the Lagrange interpolation on the reference simplex $\hat K$, then the following relation hold:
\begin{align}
(I_h^*(t) v) \circ F_K = \hat I_{h, \hat K} \big(v\circ \Phi_h^*(\cdot, t) \circ F_K \big).
\end{align}

\item[$\bullet$] 
(Commutation of material derivative and lift) (cf. \cite[Lemma 3.5]{elliott2021unified}):
\begin{align}\label{commutation-deriv-lift}
\partial^\bullet_t (v_h^{\ell}) = \big(\partial^\bullet_{t, *} v_h \big)^{\ell}
\qquad \mbox{for} \quad v_h: \Omega_h^*(t) \to \R.
\end{align}
\item[$\bullet$] (Transport property of the lifted basis functions $\varphi_j [\x^*]^{\ell}$)  (cf.\cite[Lemma 4.1]{dziuk2013l2}):
\begin{align}\label{transport-lift-basis}
\partial^\bullet_t (\varphi_j [\x^*(t)]^{\ell}) = \big(\partial^\bullet_{t, *} \varphi_j [\x^*(t)] \big)^{\ell} = 0, \quad j = 1, \dots, M.
\end{align}

\item[$\bullet$] (Commutation of material derivative of Lagrange interpolation):
Utilizing the transport property of the lifted basis functions, we have
\begin{align}
\partial^\bullet_{t, *} I_h^*(t) v = I_h^*(t) \partial^\bullet_t v .
\end{align}
Define $\CI_h^* u := (I_h^* u)^{\ell} : \Omega(t) \to \R$. It follows from \eqref{transport-lift-basis} and \eqref{commutation-deriv-lift} that 
\begin{align}\label{commutation-interp-deriv}
\partial^\bullet_t \CI_h^* u = \CI_h^* \partial^\bullet_t u .
\end{align}
It plays a crucial role in proving the optimal-order in consistency error.

\item[$\bullet$] (Approximation properties for the Lagrange interpolation)
There exists a constant $C > 0$ independent of $h \le h_0$, with a sufficiently
small $h_0 > 0$, and $t$ such that for 
$u(\cdot,t) \in H^{k+1}(\Omega_{\pm}(t))$, 
for $0\le t \le T$ (cf.\cite[Lemma 7.5]{KoLiESFEM2017}):
\begin{align}\label{lem:interpolation-err}
\|u - \CI_h^*(t) u \|_{L^2(\Omega_{\pm}(t))} + h \|\nabla u - \nabla \CI_h^*(t) u \|_{L^2(\Omega_{\pm}(t))} \le C h^{k+1} \|u\|_{H^{k+1} (\Omega_{\pm}(t))}.
\end{align}

\item[$\bullet$](Equivalence of norms)
The $L^p$ and $W^{1, p}$ norms on the discrete and continuous domains are equivalent for $1 \le p \le \infty$, uniformly in the mesh size $h \le h_0$ (with sufficiently small $h_0 >0$) and in $t \in [0,T]$.
In particular, for $v_h: \Omega_h^*(t) \to \R$ with lift $v_h^{\ell} : \Omega(t) \to \R$,  there is a constant $C$ such that for $h \le h_0$ and $0 \le t \le T$, 
\begin{align}\label{norm-equiv-L2}
&C^{-1}\|v_h\|_{L^2(\Omega_h^*(t))} \le \|v_h^{\ell}\|_{L^2(\Omega(t))}  \le C\|v_h\|_{L^2(\Omega_h^*(t))},  \\
\label{norm-equiv-H1}
&C^{-1}\|v_h\|_{H^1(\Omega_h^*(t))} \le \|v_h^{\ell}\|_{H^1(\Omega(t))}  \le C\|v_h\|_{H^1(\Omega_h^*(t))}.    
\end{align} 
The lift operator $\ell$ maps a function on the interpolated domain $\Omega_h^*(t)$ to a function on the exact domain $\Omega(t)$, provided that $\Omega_h^*(t)$ is sufficiently close to  $\Omega(t)$.
\end{enumerate}
In the following we drop the argument $t$ when it is not essential.

\subsection{Weak formulation on the interpolated evolving mesh}
For the simplicity of notation, we denote by $\rho_h^*$ and $\mu_h^*$ the piecewise constant functions in $\Omega_h^*(t)$ defined by
\begin{align}
\rho_h^* = 
\left \{
\begin{aligned}
\rho_{+} \quad \mbox{in} \,\, \Omega_{h, +}^*(t),\\
\rho_{-} \quad \mbox{in} \,\, \Omega_{h, -}^*(t),
\end{aligned}
\right .  \quad \mbox{and} \quad
\mu_h^* = 
\left \{
\begin{aligned}
\mu_{+} \quad \mbox{in} \,\, \Omega_{h, +}^*(t),\\
\mu_{-} \quad \mbox{in} \,\, \Omega_{h, -}^*(t) {\b,}
\end{aligned}
\right .
\end{align}
and define the following inner product: 
\begin{align}
(u_h, v_h)_{\Omega_h^*} = \int_{\Omega_{h,+}^*(t)} u_h v_h\, \d x + \int_{\Omega_{h, -}^*(t)} u_h v_h\, \d x . 
\end{align}
After replacing the solution by its finite element interpolation in weak formulation \eqref{PDE-ALE-weak}, and replacing domains $\Omega(t)$ and $\Omega_{\pm}(t)$ by the interpolated domains $\Omega_{h}^*(t)$ and $\Omega_{h,\pm}^*(t)$, respectively, we obtain the following weak formulation satisfied by the interpolated functions $u_h^*$, $p_h^*$, $w_h^*$ and $\phi_h^*$:
\begin{align}\label{interpolate-ALE-weak}
\left\{ \,\,
\begin{aligned}
&(\rho_h^* \partial^\bullet_{t, *} u_h^*, v_h)_{\Omega_h^*} 
+ 2(\mu_h^* \D(u_h^*), \D(v_h))_{\Omega_h^*} - (p_h^*, \nabla \cdot v_h)_{\Omega_h^*}  \\
& \hspace{70pt} + (\rho_h^* (u_h^*-w_h^*)\cdot\nabla u_h^*, v_h)_{\Omega_h^*}  -{(f_h^*, v_h)_{\Omega_h^*}} = \E_{u1}(t, v_h) \,  &&  \forall \,\, v_h \in X_h[\x^*(t)],\\[4pt]
& (\nabla\cdot u_h^*, q_h)_{\Omega_h^*}  =  \E_{u2}(t, q_h) && \forall \,\, q_h \in M_h[\x^*(t)],\\[4pt]
& (\nabla w_h^* , \nabla \chi_h)_{\Omega_h^*}  =   \E_{w}(t, \chi_h)   && \forall \,\, \chi_h \in \mathring X_{h}[\x^*(t)] , \\
& \partial_t \phi_h^*  = w_h^* \circ\phi_h^* , 
\end{aligned}
\right .
\end{align}
with initial conditions $\phi_h^*(0) ={\rm id}|_{\Omega_h^0}$ and $u_h^*(0)=u_h^0$,  {and $ f_h^* $ is the Lagrange interpolation of $ f $ on $\Omega_h^*(t) $.}
where $\E_{u1}(t, v_h)$, $\E_{u2}(t, q_h)$ and $\E_{w}(t, \chi_h)$ are consistency errors arising from finite element approximations to the solution and domain. For example, 
\begin{align}
\E_{u1}(t, v_h) :=\, & 
(\rho_h^* \partial^\bullet_{t, *} u_h^*, v_h)_{\Omega_h^*} 
- (\rho \partial^\bullet_t u, v_h^{\ell}) \notag\\[2pt]
& +  2(\mu_h^* \D(u_h^*), \D(v_h))_{\Omega_h^*}
- 2(\mu \D(u), \D(v_h^{\ell}))   \notag\\[2pt]
&+  (\rho_h^* (u_h^*-w_h^*)\cdot\nabla u_h^*, v_h)_{\Omega_h^*}
- (\rho (u - w)\cdot\nabla u, v_h^{\ell}) \notag\\[2pt]
& - (p_h^*, \nabla \cdot v_h)_{\Omega_h^*} 
+ (p, \nabla \cdot v_h^{\ell})  - {(f_h^*, v_h)_{\Omega_h^*} +  (f, v_h^l)}, \notag 
\end{align}
which contains errors arising from approximating domain $\Omega_\pm(t)$ by $\Omega_{h,\pm}^*(t)$. 
This can be controlled by using the following estimates for the error caused by domain perturbation (see \cite{KoLiESFEM2017}): 
\begin{align*}
\Big | \int_{\Omega_{h,\pm}^*(t)}g \chi_h\d x - \int_{\Omega_\pm(t)}g^\ell\chi_h^\ell \d x\Big |&\le C h^{k}\|\chi_h\|_{L^2(\Omega_{h,\pm}^*(t))}\|g\|_{L^2(\Omega_{h,\pm}^*(t))}, \\
\Big | \int_{\Omega_{h,\pm}^*(t)} g\nabla \chi_h \d x - \int_{\Omega_\pm(t)}g^\ell \nabla \chi_h^\ell \d x \Big | & \le C h^{k} \|\chi_h\|_{H^1(\Omega_{h,\pm}^*(t))} \|g\|_{L^2(\Omega_{h,\pm}^*(t))} , \\
\Big|\int_{\Omega_{h,\pm}^*(t)}\nabla\chi_h\cdot\nabla \psi_h\d x - \int_{\Omega_\pm(t)}\nabla\chi_h^\ell\cdot \nabla\psi_h^\ell \d x\Big|& \le C h^{k}\|\chi_h\|_{H^1(\Omega_{h,\pm}^*(t))}\|\psi_h\|_{H^1(\Omega_{h,\pm}^*(t))} ,
\end{align*}
where $\chi_h,\psi_h\in S_h^k[\x^*(t)]$ and $g\in L^2(\Omega_h^*(t))$. 
By using this result and the estimates {of the interpolation error} shown in \eqref{lem:interpolation-err}, the following estimates of the consistency errors can be derived when the domain $\Omega(t)$ and interface $\Gamma(t)$ are smooth, and the solution is piecewise smooth in each subdomain $\Omega_\pm(t)$: 
\begin{align}\label{consistency-1}
&\big|\E_{u1}(t, v_h)\big| \le Ch^k \|v_h\|_{H^1(\Omega_h^*(t))}
\quad \mbox{and} \quad 
\big|\E_{u2}(t, q_h)\big| \le Ch^k \|q_h\|_{L^2(\Omega_h^*(t))}, \\
&\big|\E_{w}(t, \chi_h)\big| \le Ch^k \|\chi_h\|_{H^1(\Omega_{h}^*(t))} , 
\end{align}
for $v_h\in X_{h}[\x^*(t)] $, $q_h\in M_{h}[\x^*(t)] $ and $\chi_h\in \mathring X_{h}[\x^*(t)] $. 
In the case $\partial_{t,*}^\bullet v_h=0$ and $\partial_{t,*}^\bullet q_h=0$ on $\Omega_h^*(t)$, the derivatives of these consistency errors also satisfy similar estimates, i.e., 
\begin{align} \label{consistency-2}
\Big|\frac{\d}{\d t} \E_{u1}(t, v_h)\Big| \le Ch^k \|v_h\|_{H^1(\Omega_h^*(t))}
\quad \mbox{and} \quad 
\Big|\frac{\d}{\d t} \E_{u2}(t, q_h)\Big| \le Ch^k \|q_h\|_{L^2(\Omega_h^*(t))} .
\end{align}

\begin{remark}\upshape
From the expression of $\E_{u1}(t, v_h) $ we can see that the first consistency estimate in \eqref{consistency-2} requires not only $\|u\|_{H^{k+1}(\Omega_\pm(t))}$ and $\|\partial_t^\bullet u\|_{H^{k+1}(\Omega_\pm(t))}$ but also $\|\partial_t^\bullet\partial_t^\bullet u\|_{H^{k+1}(\Omega_\pm(t))}$ to be bounded. This is the strongest regularity required for the consistency estimates, and is assumed in Proposition \ref{lemma:consistency} for consistency estimates and Theorem \ref{Main-THM} for error estimates. 
\end{remark}

\subsection{Matrix-vector form  of the weak formulation}
For both computation and error analysis it is convenient to use the matrix-vector form of the finite element weak formulation in \eqref{fem-ALE-weak}. Through the utilization of finite element basis functions $\varphi_j[\x(t)]$ on the domain $\Omega_h(t)$, the numerical solution $u_h(t)$ can be expressed as a column vector
{$\u = (u_1(t), \dots, u_{M}(t))$}, where
\begin{align*}
& \displaystyle  u_h(t) = {\sum_{j=1}^{M}} u_j(t)\varphi_j[\x(t)] \qquad \mbox{for} \qquad  {u_j(t)\in \R^d}. 
\end{align*}

Similarly, the solutions $p_h(t)$ and $w_h(t)$ on $\Omega_h(t)$ can also be represented as column vector $\p $ and $\w$, respectively.
Then $\x(t)=(x_1(t),\dots,x_M(t)) =\phi_h(\x^0,t)$ is the vector of positions of all finite element nodes at $t$, satisfying 
\begin{align}
\frac{\d\x(t)}{\d t} = \w(t) \quad \mbox{with} \quad
\x(0) =\x^0 , 
\end{align}
where the initial value $\x^0$ is given by the triangulation of the domain $\Omega$ at time $t=0$. 

We introduce domain-dependent mass matrix $\M(\x)$, the stiffness matrix $\A(\x)$, and the matrix  $
\Mrho (\x)$ and $\Amu(\x)$ on the domain $\Omega_h(t)$ as follows: 
\begin{align*}
& \v^\top \M(\x)\u  = \sum_{K\in \mathcal{K}_h[\x]} \int_{K} u_h\cdot v_h  \, \d x, \,\,
&& \v^\top \Mrho (\x)\u
= \sum_{K\in \mathcal{K}_h[\x]} \rho_K \int_{K} u_h\cdot v_h  \, \d x,\\
& \v^\top \A(\x) \u
= \sum_{K\in \mathcal{K}_h[\x]}  \int_{K} \nabla u_h \cdot \nabla v_h \, \d x, \,\,
&&\v^\top \Amu(\x) \u
= \sum_{K\in \mathcal{K}_h[\x]} 2\mu_K \int_{K} \D (u_h) \cdot \D (v_h) \, \d x.
\end{align*} 

With identity matrix $I_d \in \R^{d\times d}$, we define $\KK(\x^n)^d$ as the Kronecker product of $I_d$ and $\KK(\x^n):=\M(\x^n)+\A(\x^n)$, i.e., 
$$
{\bf K}(\x^n)^d := I_d \otimes ({\bf M}(\x^n) + {\bf A}(\x^n)). 
$$
To simplify the notation, we will use $\KK(\x^n)$ to represent $\KK(\x^n)^d$ when the dimension of the matrix is clear and therefore no confusion arises. With the matrices defined above, the $L^2$ and $H^1$ norm of finite element functions can be expressed as quadratic forms: 
\begin{align}\label{norm-M}
& \|\u\|^2_{\M(\x)} := \u^\top \M(\x) \u = \|u_h\|^2_{L^2(\Omega_h(t))},\\
\label{norm-A}
& \|\u\|^2_{\A(\x)} :=  \u^\top \A(\x) \u = \|\nabla u_h\|^2_{L^2(\Omega_h(t))},\\
\label{norm-K}
& \|\u\|^2_{\KK(\x)} :=  \u^\top \KK(\x) \u  = \|u_h\|^2_{H^1(\Omega_h(t))}.
\end{align}
In the same way, we can define the matrices $\Brho(\x, \u)$ and $\C(\x)$ such that
\begin{align*}
&\v^\top \Brho(\x, \u){\bm \chi}
=\sum_{K\in \mathcal{K}_h[\x]}\rho_K \int_{K} (u_h\cdot\nabla \chi_h)\cdot v_h \, \d x , \\
& \q^\top \C(\x)\v
=\sum_{K\in \mathcal{K}_h [\x]}\int_{K}( \nabla\cdot v_h) q_h\, \d x = \v^\top \C(\x)^\top \q.
\end{align*} 

In situations where there is no ambiguity in the context, we omit the time dependency and rewrite
the weak formulation \eqref{fem-ALE-weak} into the following matrix-vector form:
\begin{subequations}\label{fem-vec}
\begin{empheq}[left={\empheqlbrace\,\,}]{align}
&\Mrho(\x) \dot \u 
+  \Brho(\x, \u-\w)\u 
+ \Amu(\x)\u - \C(\x)^\top \p  - {\M(\x) \f} =0 , \label{semi-u-M} \\[4pt] 
&\C(\x)\u = 0 , \label{semi-p-M} \\[4pt]
& {\bm\chi}^\top \A(\x) \w = {\bf 0} \,\,\,\mbox{for nodal vectors ${\bm\chi}$ associated to}\,\, \chi_h\in \mathring X_h[\x(t)] , \\[4pt]
&\dot \x = \w, \label{semi-x-M}
\end{empheq}
\end{subequations}
{where the column vector $ \f $ collects the nodal values of $ f_h $. }
%
%
%
In the same way, we can define the corresponding matrices and vectors over the interpolated
domain $\Omega_h^*(t)$. The matrices and vectors denoted with a superscript $*$ collect the nodal values of the interpolation of the exact solutions. For example, the vector $\v^*$ collects the nodal values of $v_h^*$. With this convention, we can rewrite the weak formulation \eqref{interpolate-ALE-weak} in the following form:
\begin{subequations}\label{cons-vec}
\begin{empheq}[left={\empheqlbrace\,\,}]{align}
&\Mrho(\x^*) \dot \u^*  + \Brho(\x^*, \u^* -\w^*)\u^* + \Amu(\x^*)\u^* - \C(\x^*)^{\top} \p^* - {\M(\x^*) \f^*}  = \EE_{\u1} , \label{u-M} \\[4pt] 
&\C(\x^*)\u^* =  \EE_{\u2}, \label{p-M} \\[4pt]
& {\bm\chi}^\top \A(\x^*) \w^* =  {\bm\chi}^\top\EE_{\w}
\,\,\,\mbox{for nodal vectors ${\bm\chi}$ associated to}\,\, \chi_h\in \mathring X_h[\x^*(t)],\\[4pt]
&\dot \x^* = \w^* , \label{x-M} 
\end{empheq}
\end{subequations}
where $\EE_{\u1}$, $\EE_{\u2}$, $\EE_{\w}$ denote the vectorized forms of the defect terms $\E_{u1}(t, \cdot)$,  $\E_{u2}(t, \cdot)$, and $\E_{w}(t, \cdot)$,
which satisfy the following relations:
\begin{align}
\v^\top \EE_{\u1}  = \E_{u1}(t, v_h^*), \qquad
\q^\top \EE_{\u2}  =  \E_{u2}(t, q_h^*), \qquad
{\bm\chi}^\top \EE_{\w}  =  \E_{w}(t, \chi_h^*)   
\end{align}
for $v_h^*\in X_h[\x^*(t)]$, $q_h^* \in M_h[\x^*(t)]$, $\chi_h^* \in \mathring X_{h}[\x^*(t)]$ and the associated nodal vectors $\v$, $\q$, $\boldsymbol{\chi}$. 

Convergence of the ALE evolving FEM is established in this article by comparing the matrix-vector formulations in \eqref{fem-vec} and \eqref{cons-vec}. By subtracting \eqref{cons-vec} from \eqref{fem-vec}, we obtain the following equations for the errors $\e_{\w} = \w - \w^*$, $\e_{\u} = \u- \u^*$, $\e_{\p} = \p - \p^*$, and $\e_{\x} = \x - \x^*$: 
\begin{subequations}\label{err-defect}
\begin{empheq}[left={\empheqlbrace\,\,}]{align}
&\Mrho(\x) \dot \e_\u 
+ \Amu(\x^*)\e_\u - \C(\x)^{\top}\e_\p 
= -(\Mrho(\x)-\Mrho(\x^*)) \dot \u^* 
\label{Error_Eq_u} \\ 
&\hspace{154pt} - (\Amu(\x)-\Amu(\x^*))\e_{\u} \notag\\[4pt]
&\hspace{154pt} - (\Amu(\x)-\Amu(\x^*))\u^* \notag \\[4pt]
&\hspace{154pt} - \big(\Brho(\x, \u -\w)\u  - \Brho(\x^*, \u^* -\w^*)\u^*\big) \notag\\[4pt]
&\hspace{154pt} + (\C(\x)^{\top}-\C(\x^*)^{\top})\p^* 
+ { (\M(\x) \f - \M(\x^*) \f^*)}
- \, \EE_{\u1} , \notag\\[4pt]
&\C(\x)\e_{\u}  = -(\C(\x)-\C(\x^*))\u^* - \EE_{\u2} , 
\label{Error_Eq_p} \\[4pt]
&\boldsymbol{\chi}^\top\A(\x) \e_{\w} = {- \boldsymbol{\chi}^\top[\A(\x) - \A(\x^*)] \w^*} - \boldsymbol{\chi}^\top \EE_{\w} ,
\label{Error_Eq_w} \\[4pt]
&\dot \e_{\x} = \e_\w ,
\label{Error_Eq_x}  
\end{empheq}
\end{subequations}
for all ${\bm\chi}$ associated to some $\chi_h\in \mathring X_h[\x^*(t)]$.


The estimates of the consistency errors in \eqref{consistency-1}--\eqref{consistency-2} can also be written into the matrix-vector form, {as shown in the following proposition. This proposition generalizes the consistency estimates found in \cite[Lemma 8.1]{KoLiESFEM2017}, as it directly follows from the estimates for the error caused by domain perturbation and the approximation properties of Lagrangian interpolation.} 
\begin{proposition}[Consistency estimates]\label{lemma:consistency}
We assume that the flow map $\phi:\Omega_\pm(0)\times[0,T]\rightarrow \R^d$ and its inverse $\phi(\cdot,t)^{-1}:\Omega_\pm(t)\rightarrow \Omega_\pm(0)$ are both sufficiently smooth so that the subdomains $\Omega_\pm(t)$ are sufficiently smooth and the interpolated triangulations $\K_h[\x^*(t)]$ are shape-regular and quasi-uniform for $t\in[0,T]$, and assume that the solution of \eqref{PDE} is piecewise smooth on these subdomains in the sense that 
$$
\max_{t\in[0,T]} \sum_{j=0}^2 \|(\partial_t^\bullet)^j u\|_{H^{k+1}(\Omega_\pm(t))} \le C_0
\quad\mbox{for some constant $C_0>0$} .
$$
Then there exist positive constants $h_0$ such that for $h\le h_0$, the following consistency estimates hold: 
\begin{align}
&\big|\v^\top\EE_{\u1}\big| \le Ch^k \|\v\|_{\KK(\x^*)}
\quad \mbox{and} \quad 
\big|\q^\top\EE_{\u2}\big| \le Ch^k \|\q\|_{\M(\x^*)}, \\
&\big|\boldsymbol{\chi}^\top\EE_{\w} \big| \le Ch^k {\|\boldsymbol{\chi}\|_{\KK(\x^*)}}
\quad \mbox{for ${\bm\chi}$ associated to some}\,\,\, \chi_h\in \mathring X_{h}[\x^*(t)] .
\end{align}
Moreover, the derivatives of these defect terms satisfy similar estimates:
\begin{align} 
\Big|\v^\top \frac{\d}{\d t}\EE_{\u1}\Big| \le Ch^k \|\v\|_{\KK(\x^*)}
\quad \mbox{and} \quad
\Big|\q^\top \frac{\d}{\d t} \EE_{\u2}\Big| \le Ch^k \|\q\|_{\M(\x^*)}.
\end{align}
\end{proposition}


\section{Convergence of the ALE evolving FEM}\label{section:proof}
We are now in the position to formulate the main result of this paper, which yields optimal-order convergence in the $H^1$ norm for the Taylor--Hood finite element semidiscretization of problem \eqref{fem-vec}. 
\begin{theorem}
[Convergence of the ALE evolving FEM] \label{Main-THM} 
We assume that the solution of \eqref{PDE} satisfies the following global regularity
$$(u,p)\in C([0,T];H^1(\Omega)^d)\times L^2(0,T;L^2(\Omega)), $$ 
and is piecewise smooth in the sense that 
\begin{align}
& \max_{t\in[0,T]} \big( \|u\|_{H^{k+1}(\Omega_\pm(t))} 
+ \|\partial_t^\bullet u\|_{H^{k+1}(\Omega_\pm(t))} 
+ \|\partial_t^\bullet\partial_t^\bullet u\|_{H^{k+1}(\Omega_\pm(t))}\big) \le C_0 \\
&{\max_{t\in[0,T]}( \|p\|_{H^{k}(\Omega_\pm(t))}+ \|\partial_t^\bullet p\|_{H^{k}(\Omega_\pm(t))}) \le C_0 }
\quad\mbox{for some constant $C_0>0$}, 
\end{align}
and assume that flow map $\phi:\Omega_\pm(0)\times[0,T]\rightarrow \R^d$ and its inverse $\phi(\cdot,t)^{-1}:\Omega_\pm(t)\rightarrow \Omega_\pm(0)$ are both sufficiently smooth so that the subdomains $\Omega_\pm(t)$ are sufficiently smooth and the interpolated triangulations $\K_h[\x^*(t)]$ are shape-regular and quasi-uniform for $t\in[0,T]$. 
%
%
Then there exists a positive constant $h_0$ such that for $h\le h_0$, the finite element method in \eqref{fem-ALE-weak} with Taylor--Hood finite elements of degree $k\ge 2$ has following error bound: 
\begin{align} \label{main-error-est}
\max_{t\in[0,T]} \|\phi_h(\cdot, t) - \phi_h^*(\cdot, t)\|_{H^1(\Omega_h^0)} 
&+\max_{t\in[0,T]} \|u_h \circ \phi_h(\cdot, t) - u_h^* \circ \phi_h^*(\cdot, t) \|_{H^1(\Omega_h^0)}   \notag\\
&+ \max_{t\in[0,T]} \|w_h \circ \phi_h(\cdot, t) - w_h^* \circ \phi_h^*(\cdot, t) \|_{H^1(\Omega_h^0)}  \notag\\
&+\|p_h \circ \phi_h(\cdot, s) - p_h^* \circ \phi_h^*(\cdot, s) \|_{L^2(0,T;L^2(\Omega_h^0))} 
\le Ch^{k} ,
\end{align} 
where $\phi_h^*\in S_h^k[\x^0]^d$, $(u_h^*, p_h^*) \in X_h[\x^*(t)] \times M_h[\x^*(t)]$ and  $w_h^* \in u_h^* + \mathring X_h[\x^*(t)]$ are the Lagrange interpolations of the exact solutions. 
\end{theorem}
The proof of Theorem~\ref{Main-THM} will be presented in the next two subsections based on the technique of a homotopy map between two different domains $\Omega_h^*(t)$ and $\Omega_h(t)$.

\subsection{Comparison of integrals on two different domains}\label{auxiliary-stability}
To establish the convergence of the numerical schemes, it is necessary to compare the integrals over two different domains, i.e, the triangulated domain $\Omega_h^*(t)$  obtained from interpolating the exact domain $\Omega(t)$, and the triangulated domain $\Omega_h(t)$ determined by the numerical solution. In the spirit of the techniques in \cite{KoLiESFEM2017,Edelmann-2021}, we can obtain a similar sequence of results employing shape derivatives by constructing a homotopy map {between $\Omega_h^*(t)$} and $\Omega_h(t)$. 

To handle the rate of change of an integral over a moving domain, we will frequently make use
of the following lemma (cf. {\cite[Lemma 5.7]{shapes-of-things}}):
\begin{lemma}\label{leibniz}
If the domain $\Omega(t)$ moves with velocity $w\in W^{1,\infty}(\Omega(t))$, then 
\begin{align}
\frac{\d}{\d t}\int_{\Omega(t)} f \d x & = \int_{\Omega(t)} (\partial_t^\bullet f  + f\nabla\cdot  w) \d x,
\end{align}
where the material derivative $\partial_t^\bullet f = \partial_t f + \nabla f \cdot w$.
\end{lemma}
\begin{remark}
The commutation of the material derivative and gradient is crucial in the error analysis. Through direct computation, {we obtain the following identities}:
\begin{align}\label{commutation-material-derivative-gradient}
& \partial_t^\bullet (\nabla v) =  \nabla (\partial_t^\bullet  v) - \nabla v \nabla w
\quad  \mbox{and} \quad 
\partial_t^\bullet (\nabla \cdot v) =  \nabla \cdot (\partial_t^\bullet  v) - {\rm{tr}}[\nabla v \nabla w].
\end{align}
\end{remark}

%





We establish a linear continuous deformation from $\Omega_h^*(t)$ to $\Omega_h(t)$. Using the basis functions $\varphi_j[\x(0)]$ of $S^k_h(\Omega_h^0)^d$, we uniquely determine the domains $\Omega_h^*(t)$ and $\Omega_h(t)$ by vectors $\x^*$ and $\x$, respectively.
The vector 
\begin{align}
\x^\theta = (1-\theta) \x^* + \theta \x = \x^* + \theta \e_\x,  \quad \mbox{for} \quad 0\le \theta \le 1, 
\end{align}
then defines an intermediate domain $\Omega_h^{\theta}(t)$ that changes continuously from $\Omega_h^*(t)$ to $\Omega_h(t)$ when the parameter $\theta$ varies over the interval $[0,1]$. 
For the vector {$\u^* = (u_1^*(t), \dots, u_{M}^*(t))$}, we define the finite element function on $\Omega_h^\theta(t)$ as
\begin{align}\label{func-domian-theta-u}
u_h^{*, \theta} = {\sum_{j=1}^{M}} u_j^*(t) \varphi_j[\x^\theta(t)] \quad \mbox{for} \quad  {u_j^*(t)\in \R^d} \quad \mbox{and} \quad  0\le \theta \le 1.
\end{align}
Similarly, we can define functions $e_x^{\theta}$, $e_u^{\theta}, e_w^\theta$, $w_h^{*, \theta}$ and $(\partial_{t,*}^\bullet u_h^*)^\theta$, which are defined on intermediate domain $\Omega_h^\theta(t)$ and share the nodal vectors $\e_\x, \e_\u, \e_\w$, $\w^*$ and $\dot \u^*$, respectively. These notations are consistently applied to other functions without redundant introductions.

Let $\partial_{\theta}^\bullet$ denote the material derivative with respect to the velocity field $e_x^{\theta}$. 
By using the transport property, we have
\begin{align}\label{dtustar}
\partial_{\theta}^\bullet u_h^{*,\theta} = 0 \qquad \mbox{and} \qquad  \partial_{\theta}^\bullet e_u^{\theta} = 0.
\end{align}
Then the following relations follow from \eqref{dtustar} and \eqref{commutation-material-derivative-gradient}:  
\begin{align}
& \partial_{\theta}^\bullet (\nabla e_u^\theta) =  -  \nabla e_u^\theta\nabla e_x^\theta
\qquad \,\,\, \mbox{and} \quad \quad 
\partial_\theta^\bullet (\nabla \cdot e_u^\theta) =  - {\rm tr}[ \nabla e_u^\theta \nabla e_x^\theta ],\\
&\partial_{\theta}^\bullet (\nabla u_h^{*, \theta}) =  - \nabla u_h^{*, \theta} \nabla e_x^\theta \quad \, \mbox{and} \quad \quad 
\partial_{\theta}^\bullet (\nabla \cdot u_h^{*, \theta}) =  - {\rm tr}[\nabla u_h^{*, \theta} \nabla e_x^\theta].
\end{align}
We denote by $\rho_h^\theta$ and $\mu_h^\theta$ the piecewise constant functions in $\Omega_h^\theta(t)$ defined by
\begin{align}
\rho_h^\theta = 
\left \{
\begin{aligned}
\rho_{+} \quad \mbox{in} \,\, \Omega_{h, +}^\theta(t),\\
\rho_{-} \quad \mbox{in} \,\, \Omega_{h, -}^\theta(t),
\end{aligned}
\right .  \quad \mbox{and} \quad
\mu_h^\theta = 
\left \{
\begin{aligned}
\mu_{+} \quad \mbox{in} \,\, \Omega_{h, +}^\theta(t),\\
\mu_{-} \quad \mbox{in} \,\, \Omega_{h, -}^\theta(t).
\end{aligned}
\right .
\end{align}
In combination with the fundamental theorem of calculus, Lemma \ref{leibniz} and the transport property \eqref{transport-phi}, the following lemma was proved in \cite[Lemma 5.1]{Edelmann-2021}. 
\begin{lemma}\label{difference-uv}
In the above setting the following identities hold (for the corresponding finite element functions and the associated nodal vectors): 
\begin{align}
(\u^*)^\top  (\M(\x)-\M(\x^*))\v^* & = \int_0^1  \int_{\Omega_h^\theta(t)} u^{*, \theta}_h (\nabla\cdot e_x^{\theta}) v_h^{*, \theta} \d x\d \theta,\\
(\u^*)^\top  (\A(\x) - \A(\x^*))\v^* & = \int_0^1  \int_{\Omega_h^\theta(t)} \nabla u_h^{*, \theta} (D_{\Omega_h^{\theta}} e_x^{\theta})\cdot \nabla v_h^{*, \theta}\d x \d \theta,\\
%
(\u^*)^\top   (\Mrho(\x)-\Mrho(\x^*))\v^* & = \int_0^1  \int_{\Omega_h^\theta(t)} \rho_h^\theta u^{*, \theta}_h (\nabla\cdot e_x^{\theta}) v_h^{*, \theta} \d x\d \theta,\\%
(\u^*)^\top  (\Amu(\x) - \Amu(\x^*))\v^* &  = 2
\int_0^1  \int_{\Omega_h^\theta(t)} \mu_h^\theta  \Big[\D (u_h^{*, \theta} ) (\nabla \cdot e_x^{\theta})  \cdot \D (v_h^{*, \theta} )  \\
&\quad \, - 
\S(\nabla u_h^{*, \theta}  \nabla e_x^{\theta})
\cdot \D (v_h^{*, \theta} ) 
-\D (u_h^{*, \theta} ) \cdot \S( \nabla v_h^{*, \theta} \nabla e_x^{\theta})\Big] \, \d x \d \theta, \notag\\
(\u^*)^\top  (\C(\x) - \C(\x^*))^\top \p^* & 
= \int_0^1  \int_{\Omega_h^\theta(t)} 
\Big[(\nabla \cdot u_h^{*, \theta})(\nabla \cdot e_x^\theta) - \mbox{tr}(\nabla u_h^{*, \theta} \nabla e_x^{\theta})\Big] p_h^{*, \theta} \d x \d \theta,
\end{align}
where $D_{\Omega_h^{\theta}} e_x^{\theta} = \nabla \cdot e_x^{\theta} - 2\S(\nabla e_x^{\theta})$ 
with $\S(E) =  \frac12(E + E^\top) $.
\end{lemma}

A direct consequence of Lemma~\eqref{difference-uv} is the following conditional equivalence of
norms, i.e., if $\nabla e_x^{\theta}$ is small, the norms of finite element functions on the two domains $\Omega_h^*(t)$ and $\Omega_h(t)$ with same nodal vectors are equivalent. 
This lemma was proved in \cite[Lemma 5.2]{Edelmann-2021} by using the Gronwall's inequality. 
\begin{lemma}\label{equi-MA}
If $\|\nabla\cdot e_x^{\theta}\|_{L^\infty(\Omega_h^\theta(t))} \le \alpha$ for $0 \le \theta \le 1$, then
\begin{align*}
\|\v\|_{\M(\x^\theta)} \le e^{\alpha/2} \|\v\|_{\M(\x^*)}.
\end{align*}
If $\|D_{\Omega_h^\theta} e_x^{\theta} \|_{L^\infty(\Omega_h^\theta(t))} \le \eta$ for $0 \le \theta \le 1$, then 
\begin{align*}
\|\v\|_{\A(\x^\theta)} \le e^{\eta/2} \|\v\|_{\A(\x^*)}.
\end{align*}
\end{lemma}

\begin{remark}[Norm equivalence]\label{remark-norm-equivalence}
By using the Poincar\'e's inequality, the norms $\| \cdot \|_{\A(\x^*)}$ and $\| \cdot \|_{\KK(\x^*)}$ are equivalent for finite element functions in $H_0^1(\Omega_h^*(t))$, and the equivalence is independent of $h$ since there is a one-to-one $W^{1,\infty}$-uniformly bounded lift map from $\Omega_h^*(t)$ onto $\Omega_h(t)$.

In combination with the above Lemmas, the definition of the norms in \eqref{norm-M}--\eqref{norm-K} implies the following conclusion under the condition in Lemma~\ref{equi-MA}:
\begin{align}\label{norm-equivalence-theta}
\begin{aligned}
&\mbox{the norms $\| \cdot  \|_{\M(\x^\theta)}$ are $h$-uniformly equivalent for $0 \le  \theta \le  1$,}\\
&\mbox{and so are the norms $\| \cdot  \|_{\A(\x^\theta)}$ and  $\| \cdot  \|_{\KK(\x^\theta)}$.} 
\end{aligned}
\end{align}
Then by using Korn's inequality, the estimate \eqref{lift-approx}, Lemma~\ref{equi-MA} and the equivalence \eqref{norm-equivalence-theta}, the following norms are $h$-uniformly equivalent when $h$ is sufficiently small for $0 \le  \theta \le  1$ under the condition in Lemma~\ref{equi-MA}:
\begin{align}\label{norm-equivalence-theta-MK}
\|\cdot\|_{\M(\x^\theta)}  \sim \|\cdot\|_{\Mrho(\x^\theta)}
\qquad \mbox{and} \qquad
\|\cdot\|_{\KK(\x^\theta)}  \sim {\|\cdot\|_{\A(\x^\theta)}}  \sim  \|\cdot\|_{\Amu(\x^\theta)}.
\end{align}
\end{remark}

The following lemma was proved in \cite[Lemma 5.3]{Edelmann-2021}, which says that the condition for $e_x^\theta$ in Lemma \ref{equi-MA} can be reduced to $\theta = 0$. 
\begin{lemma}\label{equi-nabla}
Let $ e_x^* = e_x^\theta |_{\theta = 0}$ be the finite element error function on $\Omega_h^*(t)$ with nodal vectors $\e_{\x}$.
If $\|\nabla e_x^*\|_{L^{\infty}(\Omega_h^*(t))} \le \frac12$, then the finite element function $v_h^{*, \theta}(t)$ defined on $\Omega_h^\theta(t)$, with $0 \le \theta \le 1$, satisfies the following estimate:  
\begin{align*}
\|\nabla v_h^{*, \theta}\|_{L^p(\Omega_h^\theta(t))} \le C_p \|\nabla v_h^*\|_{L^p(\Omega_h^*(t))} \quad \mbox{for}\quad 1\le p \le \infty,
\end{align*}
where $C_p$ depends only on $p$.
\end{lemma}

We also need results that bound the time derivatives of the mass and stiffness matrices. 
The following results can be shown in a similar way as the analogous results on a surface \cite[Lemma 4.6]{KoLiESFEM2017}.
\begin{lemma}\label{equi-MA-deriv}
For the nodal vectors $\u$, $\v \in \R^{dM}$, there holds 
\begin{align}
\u^\top \Big(\frac{\d}{\d t} \M(\x^*)\Big) \v &\le C\|\u\|_{\M(\x^*)}\|\v\|_{\M(\x^*)},\\
\u^\top \Big(\frac{\d}{\d t} \A(\x^*) \Big)\v &\le C\|\u\|_{\A(\x^*)}\|\v\|_{\A(\x^*)},\\
\u^\top \Big(\frac{\d}{\d t} \Amu(\x^*) \Big)\v &\le C\|\u\|_{\A(\x^*)}\|\v\|_{\A(\x^*)},
\end{align}
where $C$ depends only on a bound of the $W^{1,\infty}$ norm of the domain velocity $w_h^*$.
\end{lemma}

\subsection{The proof of Theorem~\ref{Main-THM}}\label{section:stability}
Let $ e_x^*\in S_h^k[\x^*(t)]^d$,  $e_u^* \in X_h[\x^*(t)]$, $e_p^* \in M_h[\x^*(t)]$ and $e_w^* \in X_{h, \pm}[\x^*(t)](e_u^*)$ be the finite element error functions on $\Omega_h^*(t)$ with nodal vectors $\e_{\x}$, $\e_\u$, $\e_\p$ and $\e_\w$, respectively. 

Let $t^* \in(0,T]$ be the maximal time such that the following inequalities hold for $t \in [0, t^*]$:
\begin{subequations}\label{ex-W1infty}
\begin{align}
&\|e_x^*(t)\|_{W^{1,\infty}(\Omega_h^*(t))} \le  h^{-\frac{d}{4}}h^{\frac{k}{2}},\label{ex-W1infty-a}\\
&\|e_u^*(t)\|_{W^{1,\infty}(\Omega_h^*(t))} \le  1, \label{ex-W1infty-b}\\ 
& \|e_w^*(t)\|_{W^{1,\infty}(\Omega_h^*(t))} \le 1. \label{ex-W1infty-c}
\end{align}
\end{subequations}
{A key issue in the proof is to ensure that the $W^{1, \infty}$ norm of the position error $e_x^*(t)$ remains small. Since \eqref{ex-W1infty-a} is small only for $k \ge 2$, the estimates in this proof cannot be derived for the linear case $k = 1$.}
In fact, $t^*$ is positive because the inequalities above hold for $t = 0$ when $h$ is smaller than some sufficiently small constant. 

\begin{enumerate}[label={\rm(\arabic*)},ref=\arabic*,topsep=2pt,itemsep=0pt,partopsep=1pt,parsep=1ex,leftmargin=20pt]
\item[(1)] Since $\x^*(0) = \x^0$ and $\u^*(0) = \u^0$, it follows that  $\e_{\x}(0) = 0$ and $\e_{\u}(0) = 0$, which gives $e_{x}^*(0) = 0$ and $e_{u}^*(0) =  0$, respectively.

\item[(2)] Testing \eqref{Error_Eq_w} with $\e_{\w}-\e_{\u}$ at $t=0$ and using the fact that $\e_\x(0)=\e_\u(0)=0$ yields the following relation:
\begin{align*}
\|\e_\w(0) \|_{\A(\x^*(0))}^2 = - \e_\w(0) ^\top \EE_\w(0).
\end{align*}
\end{enumerate}
The relation above, together with Lemma~\ref{lemma:consistency}, implies that $\|\e_\w(0) \|_{\KK(\x^*)} \le Ch^k$. By the inverse inequality of finite element function, we have
\begin{align*}
&\|e_w^*(0)\|_{W^{1,\infty}(\Omega_h^*(t))} 
\le  Ch^{-\frac{d}{2}} \| e_w^*(0) \|_{H^{1}(\Omega_h^*(t))}
\le  Ch^{-\frac{d}{2}} \| \e_\w(0) \|_{\KK(\x^*)}\le Ch^{k-\frac{d}{2}}\le 1.
\end{align*}
Since $k \ge 2 > d/2$, it follows that \eqref{ex-W1infty-b}--\eqref{ex-W1infty-c} hold for $t = 0$ when $h$ is sufficiently small.

We first prove that error estimate \eqref{main-error-est} holds for $t\in[0,t^*]$
under condition \eqref{ex-W1infty}. Then we complete the proof by showing that $t^*=T$.
\medskip

\noindent (A) \textit{Estimates for $\e_\x$}:  
Multiplying \eqref{Error_Eq_x} by matrix $\A(\x^*)$, testing with $\e_{\x}$, and dropping the omnipresent argument
$ t \in [0, t^*]$, we have
\begin{align}
(\e_{\x})^\top \A(\x^*)\dot \e_{\x} 
= (\e_{\x})^\top\A(\x^*) \e_\w. 
\end{align}
In order to apply Gronwall's inequality, we relate $\frac{\d}{\d t}\| \e_{\x}\|_{\A(\x^*)}^2$ 
to $(\e_{\x})^\top \A(\x^*)\dot \e_{\x}$ as follows
\begin{align}
\frac12 \frac{\d}{\d t}\| \e_{\x}\|_{\A(\x^*)}^2 =\frac12  \frac{\d}{\d t} \big((\e_{\x})^\top \A(\x^*) \e_{\x}  \big) = (\e_{\x})^\top \A(\x^*)\dot \e_{\x} + \frac12 (\e_{\x})^\top\Big(\frac{\d}{\d t} \A(\x^*) \Big)\e_\x.
\end{align}
By using Lemma~\ref{equi-MA-deriv}, {we have}
\begin{align}\label{stability-ex-H1-}
\frac12 \frac{\d}{\d t}\| \e_{\x}\|_{\A(\x^*)}^2 = \, &
(\e_{\x})^\top \A(\x^*) \e_{\w} + \frac12 (\e_{\x})^\top \Big(\frac{\d}{\d t} \A(\x^*) \Big)\e_\x \\
\le \, & \int_{\Omega_h^*(t)}  \nabla e_w^* \cdot \nabla e_x^* \, \d x 
+ C\|\e_{\x}\|_{\A(\x^*)}^2
\notag\\
\le \, & C\|\nabla e_w^*\|_{L^2(\Omega_h^*(t))} \|\nabla e_x^*\|_{L^2(\Omega_h^*(t))} + C\|\e_{\x}\|_{\A(\x^*)}^2
\notag\\
\le \, & C\|\e_{\w}\|_{\A(\x^*)}^2 + C\|\e_{\x}\|_{\A(\x^*)}^2.
\notag
\end{align}

\noindent (B) \textit{Estimates for $\e_\w$ under condition \eqref{ex-W1infty}}:
Since $e_w^*- e_u^* = 0$ on $\Gamma_h^*(t)$, i.e., $e_w^*- e_u^*\in \mathring X_h[\x^*(t)]$, we can choose ${\bm\chi}=\e_{\w} - \e_\u$ in \eqref{Error_Eq_w}. This yields the following result: 
\begin{align} \label{Error_Energy_w}
\| \e_{\w} \|_{\A(\x^*)}^2
&=  - (\e_{\w} - \e_\u)^\top (\A(\x) - \A(\x^*)) \w  +  (\e_\u)^\top \A(\x^*) \e_{\w}  \\
&\quad \, - (\e_{\w} - \e_\u)^\top \EE_{\w} \notag\\
&= : \mathcal{H}_1 + \mathcal{H}_2 + \mathcal{H}_3. \notag
\end{align}

The first term on the right hand side of \eqref{Error_Energy_w} can be estimated by using Lemma \ref{difference-uv} and H\"older's inequality, i.e.,
\begin{align}
\mathcal{H}_1 
&\le C\| \nabla w_h^\theta  \|_{L^\infty(\Omega_h^\theta(t))} \|\nabla e_x^\theta\|_{L^2(\Omega_h^\theta(t))}\|\nabla e_w^\theta - \nabla e_u^\theta \|_{L^2(\Omega_h^\theta(t))}\notag\\
&\le C \big(1 + \|\nabla e_w^*\|_{L^\infty(\Omega_h^*(t))} \big) \|\nabla e_x^*\|_{L^2(\Omega_h^*(t))}
\big(\|\nabla  e_w^*\|_{L^2(\Omega_h^*(t))} + \| \nabla e_u^* \|_{L^2(\Omega_h^*(t))} \big) \notag\\
&\le C\|\e_\x\|_{\A (\x^*)} ( \|\e_\w\|_{\A (\x^*)}  + \|\e_\u \|_{\A (\x^*)} ) \notag\\
&\le C\|\e_\x\|_{\A (\x^*)}^2 + C\|\e_\u\|_{\A (\x^*)}^2 + \frac18 \|\e_\w\|_{\A (\x^*)}^2,\notag
\end{align} 
{where we have used 
$ \| \nabla w_h^\theta  \|_{L^\infty(\Omega_h^\theta(t))} =
\| \nabla e_w^\theta - \nabla w_h^{*, \theta} \|_{L^\infty(\Omega_h^\theta(t))} \le C \big(1 + \|\nabla e_w^*\|_{L^\infty(\Omega_h^*(t))}\big)$.}

The second term on the right hand side of \eqref{Error_Energy_w} can be bounded by
\begin{align}
\mathcal{H}_2 
&\le C \|\e_\u\|_{\A(\x^*)} \|\e_\w \|_{\A (\x^*)} 
\le C \|\e_\u\|_{\A(\x^*)}^2 +  \frac18\|\e_\w \|_{\A (\x^*)}^2. \notag
\end{align}
The last term on the right hand side of \eqref{Error_Energy_w} can be estimated similarly, i.e.,
\begin{align}
\mathcal{H}_3 
\le Ch^k \|\e_\w - \e_\u\|_{\KK (\x^*)} 
& \le Ch^k \|\e_\w - \e_\u\|_{\A (\x^*)} \\
& \le Ch^{2k}  + \frac18 \|\e_\u\|_{\A (\x^*)}^2 + \frac18\|\e_\w\|_{\A (\x^*)}^2 {\b,} \notag
\end{align}
where we have converted $ \|\e_\w - \e_\u\|_{\KK (\x^*)} $ (i.e., the full $H^1$ norm) to $ \|\e_\w - \e_\u\|_{\A (\x^*)} $ (i.e., the $H^1$ semi-norm) because the corresponding function $e_w^*-e_u^*$ satisfies the zero boundary condition in $\Omega_{h,\pm}^*(t)$. 
Altogether we obtain the bound as
\begin{align}\label{stability-ew-H1}
\|\e_{\w}\|_{\A(\x^*)}^2
\le Ch^{2k} + C \|\e_\x\|_{\A(\x^*)}^2 + C\|\e_\u\|_{\A(\x^*)}^2. 
\end{align}
Then, substituting \eqref{stability-ew-H1} into \eqref{stability-ex-H1-}, we obtain
\begin{align}\label{stability-ex-H1-2}
\frac{\d}{\d t}\| \e_{\x}\|_{\A(\x^*)}^2 
\le Ch^{2k} + C\|\e_{\u}\|_{\A(\x^*)}^2 + C\|\e_{\x}\|_{\A(\x^*)}^2.
\end{align}
By integrating \eqref{stability-ex-H1-2} on both sides with respect to time from $0$ to $t$, we obtain
\begin{align}\label{stability-ex-H1}
\| \e_{\x}\|_{\A(\x^*)}^2  \le  Ch^{2k} +  \int_0^t \Big(\|\e_{\x}(s)\|_{\A(\x^*)}^2 +
\|\e_\u(s)\|_{\A(\x^*)}^2 \Big) \,\d s . 
\end{align}

\noindent (C) \textit{Estimates for $\e_\u$ under condition \eqref{ex-W1infty}}: 
Testing \eqref{Error_Eq_u} by $\dot \e_{\u}$, we obtain
\begin{align} \label{Error_Energy_u}
&\hspace{-13pt} (\dot \e_\u)^\top \Mrho(\x)\dot \e_\u
+ \frac12 \frac{\d}{\d t} \Big( (\e_\u)^\top \Amu(\x^*)\e_\u  \Big) 
- \frac12 (\e_\u)^\top \Big(\frac{\d}{\d t}\Amu(\x^*)\Big)\e_\u  \\
=\, &  
- (\dot \e_\u)^\top (\Mrho(\x)-\Mrho(\x^*)) \dot \u^* \notag \\
\, &  
-  (\dot \e_\u)^\top (\Amu(\x)-\Amu(\x^*))\u^* \notag \\
\, &  
- (\dot \e_\u)^\top (\Amu(\x)-\Amu(\x^*))\e_{\u} \notag \\
& - (\dot \e_\u)^\top \Big(\Brho(\x, \u -\w)\u - \Brho(\x^*, \u^* -\w^*) \u^*   \Big) \notag\\
&+  (\dot \e_\u)^\top (\C(\x)^{\top}-\C(\x^*)^{\top})\p^*
+  (\dot\e_\u)^\top \C(\x)^\top  \e_\p \notag\\
& + {(\dot\e_\u)^\top (\M(\x) \f - \M(\x^*) \f^*)}
- (\dot \e_\u)^\top\EE_{\u1}\notag\\
=: & {\sum_{j=1}^7} S_j - (\dot \e_\u)^\top \EE_{\u1}. \notag 
\end{align}

The last term on the left-hand side of \eqref{Error_Energy_u} can be estimated by using Lemma~\ref{equi-MA-deriv}, i.e., 
\begin{align}\label{stability-eu-e2}
\Big| \frac12 (\e_\u)^\top \Big(\frac{\d}{\d t} \Amu(\x^*)\Big)\e_\u \Big| 
\le  C \|\e_\u\|_{\A(\x^*)}^2.
\end{align}
In the following, we present estimates for each term on the right-hand side of \eqref{Error_Energy_u}.

The first term on the right hand side of \eqref{Error_Energy_u} can be estimated by using Lemma \ref{difference-uv} and H\"older's inequality, and the norm equivalence in Lemma \ref{equi-MA}, i.e., 
\begin{align} \label{est-S_1}
S_1
&= (\dot \e_\u)^\top (\Mrho(\x) - \Mrho(\x^*)) \dot \u^*  \\
&= \int_0^1 \int_{\Omega_h^\theta(t)} \rho_h^\theta  (\partial_{t,*}^{\bullet} e_u^*)^\theta (\nabla\cdot e_x^{\theta}) (\partial_{t,*}^{\bullet} u_h^*)^\theta  \, \d x \d\theta \notag\\
&\le 
C\|\nabla\cdot   e_x^{\theta}\|_{L^2(\Omega_h^\theta(t))} \|(\partial_{t,*}^{\bullet} u_h^*)^\theta \|_{L^\infty(\Omega_h^\theta(t))} \|(\partial_{t,*}^{\bullet} e_u^*)^\theta\|_{L^2(\Omega_h^\theta(t))}\notag\\
&\le 
C\|\nabla e_x^*\|_{L^2(\Omega_h^*(t))} \| \partial_{t,*}^{\bullet} u_h^*\|_{L^\infty(\Omega_h^*(t))} \|\partial_{t, *}^{\bullet} e_u^*\|_{L^2(\Omega_h^*(t))}\notag\\
&\le C\|\e_\x\|_{\A(\x^*)} \|\dot \e_\u\|_{\M(\x)}\notag\\
&\le \varepsilon \|\dot \e_\u\|_{\Mrho(\x)}^2 + C \varepsilon^{-1} \|\e_\x\|_{\A(\x^*)}^2 {\b,} \notag
\end{align}
where the positive constant $\varepsilon$ can be arbitrarily small and we have used $\| \partial_{t, *}^{\bullet} u_h^*\|_{L^\infty(\Omega_h^*(t))} \le C$.

The second term on the right hand side of \eqref{Error_Energy_u} can be estimated by using integration by parts in time, i.e., 
\begin{align}\label{stability-eu-e3-}
S_2 =\, &- (\dot \e_\u)^\top (\Amu(\x)-\Amu(\x^*))\u^* \\
=\, & -\frac{\d}{\d t} \Big( (\e_\u)^\top (\Amu(\x)-\Amu(\x^*))\u^* \Big)
+ (\e_\u)^\top \Big(\frac{\d}{\d t}(\Amu(\x)-\Amu(\x^*))\Big)\u^* \notag \\
& + (\e_\u)^\top (\Amu(\x)-\Amu(\x^*))\dot \u^*  \notag\\
=:\, & S_{21} + S_{22} + S_{23} ,\notag
\end{align}
where the first term on the right-hand side of \eqref{stability-eu-e3-} can be bounded after integrating it over $[0, t]$ for $0 \le t \le t^*$, i.e., 
\begin{align}\label{stability-eu-e33}
\int_0^t S_{21}(s)\, \d s
=\,&- (\e_\u)^\top (\Amu(\x)-\Amu(\x^*))\u^* \\
=\,&- 2\int_0^1  \int_{\Omega_h^\theta(t)} \mu_h^\theta \Big[ 
\D (u_h^{*, \theta} )  (\nabla \cdot e_x^{\theta}) \cdot \D (e_u^{\theta} ) 
\quad\mbox{(Lemma \ref{difference-uv} is used)} \notag \\
\,&\hspace{80pt} - \S(\nabla u_h^{*, \theta} \nabla e_x^{\theta})\cdot \D ( e_u^{\theta} )
+ \D (u_h^{*, \theta} ) \cdot \S(\nabla e_u^{\theta} \nabla e_x^{\theta}) \Big]\, \d x \d \theta \notag\\
\le\,& C\|\nabla u_h^*(t)\|_{L^\infty(\Omega_h^*(t))} \|\nabla e_x^*(t)\|_{L^2(\Omega_h^*(t))}\|\nabla e_u^*(t)\|_{L^2(\Omega_h^*(t))} 
\quad\mbox{(Lemma \ref{equi-MA} is used)} \notag\\
\le \, &  C\|\e_\x(t)\|_{\A (\x^*)}  \|\e_\u(t)\|_{\A (\x^*)}  \notag\\
\le\,& \varepsilon \|\e_\u(t)\|_{\Amu (\x^*)}^2  + C \varepsilon^{-1}  \|\e_\x(t)\|_{\A (\x^*)}^2. \notag
\end{align}

In order to estimate $S_{22}$, we define $w_h^\theta = w_h^{*, \theta} + \theta e_w^{\theta}$ {to be the} finite element function on $\Omega_h^\theta(t)$ with nodal vector $ \w^\theta = \w^* + \theta \e_\w$.
Note that $\frac{\d}{\d t} \x^\theta = \dot \x =  \w^\theta$ and $\frac{\d}{\d \theta} \x^\theta = \e_\x$. By using the transport formula, Lemma~\ref{leibniz} and \eqref{commutation-material-derivative-gradient}, we have
\begin{align*}
&\frac{\d}{\d t}\int_{\Omega_h^\theta(t)} f  \d x  = \int_{\Omega_h^\theta(t)} (\partial_{t, \theta}^\bullet f  +  f \nabla\cdot  w_h^\theta) \d x, \\
&\frac{\d}{\d \theta}\int_{\Omega_h^\theta(t)} f \d x  = \int_{\Omega_h^\theta(t)} (\partial_{\theta}^\bullet  f  +  f \nabla\cdot  e_x^\theta) \d x,
\end{align*}
where 
$\partial_{\theta}^\bullet$  denotes the material derivative with respect to the velocity field $e_x^{\theta}$ with
\begin{align}
& \partial_\theta^\bullet (\nabla w_h^\theta) =  \nabla e_w^\theta -  \nabla w_h^\theta \nabla e_x^\theta 
\quad \mbox{and} \quad
\partial_{\theta}^\bullet (\nabla \cdot w_h^\theta) =  \nabla \cdot e_w^\theta  - {\rm tr}[\nabla w_h^\theta \nabla e_x^\theta],
\end{align}
and $\partial_{t, \theta}^\bullet$ denotes the material derivative with respect to the velocity field $w_h^{\theta}$ with 
\begin{align}
&\partial_{t, \theta}^{\bullet}  \D(\eta_h^\theta) 
= -\S(\nabla \eta_h^\theta \nabla e_x^\theta)  \qquad \mbox{for} \qquad  \partial_{t, \theta}^{\bullet} \eta_h^\theta= 0,\\
&\partial_{t, \theta}^{\bullet}\,  \S(\nabla \eta_h^\theta \nabla e_x^{\theta}) 
= \S(\nabla \eta_h^\theta \nabla e_w^{\theta}) - \S(\nabla \eta_h^\theta (\nabla w_h^\theta  \nabla e_x^{\theta} + \nabla e_x^{\theta}\nabla w_h^\theta)),  
\end{align}
where we note $\partial_{t, \theta}^{\bullet} e_x^\theta = e_w^\theta$ and $\eta_h^\theta$ is defined from \eqref{func-domian-theta-u}, which is the finite element function on $\Omega_h^\theta(t)$ and shares the nodal vector ${\boldsymbol{\eta}} \in \R^{dM}$ with $ \dot {\boldsymbol{\eta}} = 0$. 

The second term on the right hand side of \eqref{stability-eu-e3-} can be estimated by using the above-mentioned setting and Lemma~\ref{difference-uv}, we have 
\begin{align}\label{stability-eu-e32}
{|S_{22}|} =
\, & \Big| (\e_\u)^\top \Big(\frac{\d}{\d t}(\Amu(\x) - \Amu(\x^*))\Big)\u^* \Big|\\
=\,& 
\Big|\frac{\d}{\d t} \int_0^1  \int_{\Omega_h^\theta(t)} 2\mu_h^\theta \Big[\D (u_h^{*, \theta} )  (\nabla \cdot e_x^\theta) \cdot \D (e_u^{\theta} ) \notag\\
& \qquad \qquad\quad  -
\S(\nabla u_h^{*, \theta} \nabla e_x^{\theta})
\cdot \D ( e_u^{\theta} )
- \D (u_h^{*, \theta} ) \cdot 
\S(\nabla e_u^\theta \nabla e_x^{\theta}) \Big] \d x \d \theta\Big| \notag\\
=\,
&  \Big| \int_0^1 \int_{\Omega_h^\theta(t)} 2\mu_h^\theta \Big[ - \S (u_h^{*, \theta} w_h^\theta )(\nabla \cdot e_x^\theta) \cdot \D (e_u^{\theta} )
-  \D (u_h^{*, \theta})  (\nabla \cdot e_x^\theta) \cdot \S (e_u^{\theta} w_h^\theta)\notag \\
&  \qquad \qquad  \quad + \D (u_h^{*, \theta} )  \big(\nabla \cdot e_w^\theta - \mbox{tr}[\nabla e_x^\theta \nabla w_h^\theta] \big) \cdot \D (e_u^{\theta} )\Big] \d x \d \theta \notag\\
&+ \int_0^1  \int_{\Omega_h^\theta(t)} 2\mu_h^\theta 
\Big[
\Big(\S(\nabla u_h^{*, \theta} \nabla e_w^{\theta}) - \S( \nabla u_h^{*, \theta} (\nabla e_x^{\theta}\nabla w_h^\theta + \nabla w_h^\theta \nabla e_x^{\theta})) \Big) \cdot \D ( e_u^{\theta}) \notag \\
& \qquad \qquad  \quad +
\D ( u_h^{*, \theta}) \cdot
\Big(\S(\nabla e_u^{\theta} \nabla e_w^{\theta}) - \S(\nabla e_u^{\theta}(\nabla e_x^{\theta}\nabla w_h^\theta + \nabla w_h^\theta \nabla e_x^{\theta})) \Big) \notag \\
&  \qquad \qquad \quad -  
\S(\nabla u_h^{*, \theta}\nabla e_x^{\theta}) 
\cdot  \S(\nabla e_u^{\theta}\nabla w_h^\theta) - \S (u_h^{*, \theta} \nabla w_h^\theta ) \cdot 
\S(\nabla e_u^\theta\nabla e_x^{\theta})\Big] \d x \d \theta \notag \\
&+ \int_0^1  \int_{\Omega_h^\theta(t)} 2\mu_h^\theta \Big[\D (u_h^{*, \theta} )  (\nabla \cdot e_x^\theta) \cdot \D (e_u^{\theta} ) \notag \\
& \qquad \qquad  -
\S(\nabla u_h^{*, \theta} \nabla e_x^{\theta})
\cdot \D ( e_u^{\theta} )
- \D (u_h^{*, \theta} ) \cdot 
\S(\nabla e_u^\theta \nabla e_x^{\theta}) \Big](\nabla \cdot w_h^\theta) \, \d x \d \theta \Big|\notag\\
\le \,  & C\|\nabla u_h^*\|_{L^\infty(\Omega_h^*(t))}\big(1 + \|\nabla e_w^*\|_{L^\infty(\Omega_{h}^*(t))}\big)\|\nabla e_x^*\|_{L^2(\Omega_h^*(t))} \|\nabla  e_u^*\|_{L^2(\Omega_h^*(t))} \notag\\
& + C\|\nabla u_h^*\|_{L^\infty(\Omega_h^*(t))}\big(1 + \|\nabla e_w^*\|_{L^\infty(\Omega_{h}^*(t))}\big)
\|\nabla e_w^*\|_{L^2(\Omega_{h}^*(t))} \|\nabla  e_u^*\|_{L^2(\Omega_h^*(t))} \notag\\
\le \, & C \big(\|\nabla e_w^*\|_{L^2(\Omega_{h}^*(t))}  + \|\nabla e_x^*\|_{L^2(\Omega_h^*(t))}  \big) \|\nabla e_u^*\|_{L^2(\Omega_h^*(t))} \notag\\
\le \, &  C \|\e_\u\|_{\A (\x^*)}^2  + C \|\e_\x\|_{\A (\x^*)}^2 + C \|\e_\w\|_{\A (\x^*)}^2.  \notag
\end{align}
%
%
Similarly as the estimates of \eqref{stability-eu-e33}, the third term on the right hand side of \eqref{stability-eu-e3-} can be estimated as
\begin{align}
{|S_{23}|}
&= \big| (\e_\u)^\top (\Amu(\x)-\Amu(\x^*))\dot\u^* \big| \\
&\le  C\|\nabla \partial_{t, *}^{\bullet} u_h^*\|_{L^\infty(\Omega_h^*(t))}
\|\nabla e_x^*\|_{L^2(\Omega_h^*(t))}\|\nabla e_u^*\|_{L^2(\Omega_h^*(t))} \notag\\
&\le  C\|\e_\u\|_{\A (\x^*)}^2  + C \|\e_\x\|_{\A (\x^*)}^2.  \notag
\end{align}
The estimates of $\int_0^t S_{21}(s)\d s$, $S_{22}$ and $S_{23}$ show that 
\begin{align}\label{Estimate-eu-e33}
&\Big| \int_0^t S_{2}(s)\, \d s \Big| \\
&\le \varepsilon \|\e_\u(t)\|_{\Amu (\x^*)}^2  + C \varepsilon^{-1}  \|\e_\x(t)\|_{\A (\x^*)}^2 
+ C \int_0^t \Big( \|\e_\u(s)\|_{\A (\x^*)}^2  +  \|\e_\x(s)\|_{\A (\x^*)}^2 +  \|\e_\w(s)\|_{\A (\x^*)}^2\Big) \d s . \notag
\end{align}

The third term on the right hand side of \eqref{Error_Energy_u} can be written as 
\begin{align}  \label{stability-eu-S3}
S_3 
= &\, -  (\dot \e_\u)^\top (\Amu(\x)-\Amu(\x^*))\e_{\u} \\
= &\, -  \frac{\d}{\d t} \Big( \frac12 (\e_\u)^\top (\Amu(\x)-\Amu(\x^*))\e_{\u} \Big) 
        +  \frac12 (\e_\u)^\top \Big( \frac{\d}{\d t} (\Amu(\x)-\Amu(\x^*)) \Big) \e_{\u} \notag\\
= &\!:\, S_{31} + S_{32} ,  \notag
\end{align} 
where $S_{31}$ and $S_{32}$ can be treated similarly as $S_{21}$ and $S_{22}$, with the following estimates: 
\begin{align}  \label{stability-eu-S31}
\Big| \int_0^t S_{31}(s)\, \d s \Big| 
=\,&\Big|  \frac12 (\e_\u)^\top (\Amu(\x)-\Amu(\x^*))\e_{\u} \Big|  \\
\le\,& C\|e_x^*\|_{W^{1,\infty}(\Omega_h^*(t))} \|\e_\u(t)\|_{\Amu (\x^*)}^2 \notag \\ 
\le\,& Ch^{\frac{k}{2}-\frac{d}{4}} \|\e_\u(t)\|_{\Amu (\x^*)}^2  \notag\\
\label{stability-eu-S32}
{|S_{32}|} = 
\, & \Big| (\e_\u)^\top \Big(\frac{\d}{\d t}(\Amu(\x) - \Amu(\x^*))\Big)\e_{\u} \Big|\\
\le \, &  C\|e_u^*\|_{W^{1,\infty}(\Omega_h^*(t))}  \big(\|\nabla e_w^*\|_{L^2(\Omega_{h}^*(t))}  + \|\nabla e_x^*\|_{L^2(\Omega_h^*(t))}  \big) \|\nabla e_u^*\|_{L^2(\Omega_h^*(t))} \notag\\
\le \, &  C \|\e_\u\|_{\A (\x^*)}^2  + C \|\e_\x\|_{\A (\x^*)}^2 + C \|\e_\w\|_{\A (\x^*)}^2 , \notag
\end{align}
where we have used the results {$\|e_x^*\|_{W^{1,\infty}(\Omega_h^*(t))} \le h^{-\frac{d}{4}}h^\frac{k}{2}$} and $\|e_u^*\|_{W^{1,\infty}(\Omega_h^*(t))}\le 1$ in \eqref{ex-W1infty}. Therefore, 
\begin{align}\label{Estimate-eu-S3}
\Big| \int_0^t S_{3}(s)\, \d s \Big| 
&\le Ch^{\frac{k}{2}-\frac{d}{4}} \|\e_\u(t)\|_{\Amu (\x^*)}^2 
+ C\int_0^t \Big( \|\e_\u(s)\|_{\A (\x^*)}^2  +  \|\e_\x(s)\|_{\A (\x^*)}^2 +  \|\e_\w(s)\|_{\A (\x^*)}^2 \Big) \d s . 
\end{align}

The fourth term on the right hand side of \eqref{Error_Energy_u} can be estimated similarly, i.e.,
\begin{align}  \label{stability-eu-e5}
S_4 =\, & - (\dot \e_\u)^\top \Big(\Brho(\x, \u -\w)\u - \Brho(\x, \u^* -\w^*) \u^* \Big) \\
& -  (\dot \e_\u)^\top\Big(\Brho(\x, \u^* -\w^*)\u^* - \Brho(\x^*, \u^* -\w^*) \u^* \Big) \notag \\
= \, & -\int_{\Omega_h(t)} \rho_h 
 \big[ (e_u^\theta - e_w^\theta)\cdot\nabla u_h^{*,\theta}
 +(u_h^\theta - w_h^\theta)\cdot\nabla e_u^\theta \big]\cdot (\partial_{t, *}^{\bullet} e_u^{*})^\theta\Big|_{\theta = 1} \d x  \notag \\
& +\int_0^1  \int_{\Omega_h^\theta(t)} \rho_h^\theta 
\Big[ (u_h^{*, \theta}  - w_h^{*, \theta}) \cdot \big(\nabla u_h^{*, \theta} \nabla e_x^\theta
- \nabla u_h^{*, \theta}  (\nabla \cdot e_x^\theta)  \big) \Big] \cdot (\partial_{t, *}^{\bullet} e_u^{*})^\theta \d x \d \theta \notag \\
\le \, &  C  \big(\|e_u^*\|_{L^2(\Omega_h^*(t))}  + \|e_w^*\|_{L^2(\Omega_h^*(t))} \big)
\big(\| \nabla u_h^*\|_{L^\infty(\Omega_h^*(t))} + \| \nabla e_u^*\|_{L^\infty(\Omega_h^*(t))}\big) 
\|\partial_{t, *}^{\bullet} e_u^{*}\|_{L^2(\Omega_h^*(t))} \notag\\
& +  C \| u_h^* - w_h^* \|_{L^\infty(\Omega_h^*(t))}   \| \nabla u_h^*\|_{L^\infty(\Omega_h^*(t))}  \|\nabla e_x^*\|_{L^2(\Omega_h^*(t))} 
\|\partial_{t, *}^{\bullet} e_u^{*}\|_{L^2(\Omega_h^*(t))} \notag\\
\le \, & C \|\dot \e_\u\|_{\M(\x^*)} \big(\|\e_\u\|_{\A(\x^*)} + \|\e_\w\|_{\A(\x^*)} + \|\e_\x\|_{\A(\x^*)} \big)
\notag\\
\le \, &  \varepsilon \|\dot \e_\u\|_{\Mrho(\x)}^2 + C \varepsilon^{-1} \big( \|\e_\u\|_{\A(\x^*)}^2  + \|\e_\w\|_{\A(\x^*)}^2 + \|\e_\x\|_{\A(\x^*)}^2 \big). \notag
\end{align}

The fifth term on the right hand side of \eqref{Error_Energy_u} can be written as 
\begin{align}\label{stability-eu-e7-}
S_5 =\, &(\dot \e_\u)^\top (\C(\x)^{\top}-\C(\x^*)^{\top})\p^* \\
=\, & \frac{\d}{\d t} \Big((\e_\u)^\top (\C(\x)^{\top}-\C(\x^*)^{\top})\p^*  \Big)
- (\e_\u)^\top \Big(\frac{\d}{\d t}(\C(\x)^{\top}-\C(\x^*)^{\top})\Big)\p^* \notag \\
& - (\e_\u)^\top (\C(\x)^{\top}-\C(\x^*)^{\top})\dot \p^* \notag\\
=:\, & S_{51} + S_{52} + S_{53}, \notag
\end{align}
where the first term on the right-hand side of \eqref{stability-eu-e7-} can be bounded after integrating it over $[0, t]$ for $0 \le t \le t^*$, i.e., 
\begin{align}
\int_0^t S_{51}(s)\, \d s 
&=(\e_\u)^\top (\C(\x)^{\top}-\C(\x^*)^{\top})\p^*  \\
&=  \int_0^1 \int_{\Omega_h^\theta(t)} 
\Big( (\nabla\cdot e_x^\theta) (\nabla\cdot e_u^\theta) -{\rm tr}\big[ \nabla e_u^\theta \nabla e_x^\theta]
\Big) p_h^{*,\theta}  \, \d x \, \d \theta \notag \\
&\le 
C\|p_h^*\|_{L^\infty(\Omega_h^*(t))} \|\nabla e_x^*\|_{L^2(\Omega_h^*(t))} \|\nabla e_u^*\|_{L^2(\Omega_h^*(t))} \notag\\
&\le  \varepsilon \|\e_\u(t)\|_{\Amu (\x^*)}^2  + C \varepsilon^{-1} \|\e_\x(t) \|_{\A (\x^*)}^2.  \notag
\end{align}
The second and third terms on the right-hand side of \eqref{stability-eu-e7-} can be treated similarly as $S_{22}$ and $S_{23}$, with the following estimates:  
%
\begin{align} \label{stability-eu-e62}
|S_{52}|
&\le  C\|\e_\u\|_{\A (\x^*)}^2  + C \|\e_\w\|_{\A (\x^*)}^2 + C \|\e_\x\|_{\A (\x^*)}^2 ,   \\ 
|S_{53}|
&\le  C\|\e_\u\|_{\A (\x^*)}^2  + C \|\e_\x\|_{\A (\x^*)}^2.  
\end{align}

The sixth term on the right-hand side of \eqref{Error_Energy_u} can be estimated by differentiating \eqref{Error_Eq_p} in time. This yields the following relation: 
\begin{align} \label{Eq-div-dt-eu}
\C(\x)  \dot \e_\u = - \Big(\frac{\d}{\d t}(\C(\x)-\C(\x^*))\Big) \u^* 
-(\C(\x)-\C(\x^*))\dot\u^* 
- \Big(\frac{\d}{\d t}\C(\x)\Big)  \e_\u - \frac{\d}{\d t} \EE_{\u2}.
\end{align}
Testing this equation with $\e_\p$ and using relation $\q^\top \C(\x)\v = (\C(\x)\v)^\top \q$, we have 
\begin{align} \label{stability-dteu-eq}
S_6 = (\C(\x)  \dot \e_\u)^\top \e_\p 
&= - (\u^*)^\top\Big(\frac{\d}{\d t}(\C(\x)^\top -\C(\x^*)^\top)\Big)  \e_\p
-(\dot \u^*)^\top(\C(\x)^\top -\C(\x^*)^\top)\e_\p \\
&\quad \, - (\e_\u)^\top\Big(\frac{\d}{\d t}\C(\x)^\top \Big)\e_\p   - (\e_\p)^\top \Big( \frac{\d}{\d t} \EE_{\u2}\Big), \notag
\end{align}
which can be estimated similarly as $S_{22}$, $S_{23}$, and so on. In particular, we have 
\begin{align}
S_{6}
&\le \varepsilon \|\e_\p\|_{\M(\x^*)}^2
+ C \varepsilon^{-1} \big( \|\e_\u\|_{\A (\x^*)}^2 + \|\e_\w\|_{\A (\x^*)}^2 + \|\e_\x\|_{\A (\x^*)}^2 
+ h^{2k} \big).
\end{align}
{
The seventh term on the right hand side of \eqref{Error_Energy_u} can be written as 
\begin{align}
S_7 
&= (\dot\e_\u)^\top (\M(\x) \f - \M(\x^*) \f^*)\\
& =  (f_h , (\partial_{t,*}^{\bullet} e_u^*)^\theta|_{\theta = 1} )_{\Omega_h(t)}
- (f_h^*, \partial_{t,*}^{\bullet} e_u^*)_{\Omega_h^*(t)}    \qquad (\text{here we use} \, (\partial_{t,*}^{\bullet} e_u^*)^\theta|_{\theta = 0} = \partial_{t,*}^{\bullet} e_u^*)\notag\\
& = (f_h - f, (\partial_{t,*}^{\bullet} e_u^*)^\theta|_{\theta = 1} )_{\Omega_h(t)}
- (f_h^* - f, \partial_{t,*}^{\bullet} e_u^*)_{\Omega_h^*(t)}
+ \Big((f, (\partial_{t,*}^{\bullet} e_u^*)^\theta|_{\theta = 1} )_{\Omega_h(t)} - (f, \partial_{t,*}^{\bullet} e_u^* )_{\Omega_h^*(t)}\Big) \notag\\
& =: S_{71} + S_{72} + S_{73}. \notag
\end{align} 
By using the approximation properties of the Lagrange interpolations $ f_h $  and $ f_h^* $ for $ f $ on the domains $ \Omega_h(t) $ and $ \Omega_h^*(t) $, respectively,  the terms $ S_{71} $ and $ S_{72} $ can be estimates as
\begin{align}
|S_{71}| +  |S_{72}|\le Ch^k \|f\|_{H^{k-1}(\Omega(t))} \|\partial_{t, *}^{\bullet} 
e_u^*\|_{L^2(\Omega_h^*(t))} \le Ch^k\|\dot \e_\u\|_{\M(\x)} \le \varepsilon \|\dot \e_\u\|_{\Mrho(\x)}^2 + C \varepsilon^{-1} h^{2k}.
\end{align} 
By applying  $ \partial_\theta^{\bullet} f  =  \nabla f \cdot e_x^{\theta}$ on $ \Omega_h^\theta(t) $, the term $ S_{73} $ can be bounded similarly to the estimates for \eqref{est-S_1} as
\begin{align} 
|S_{73}|
&=\Big|\int_0^1 \int_{\Omega_h^\theta(t)} (\nabla f \cdot e_x^{\theta})  (\partial_{t,*}^{\bullet} e_u^*)^\theta  \, \d x \d\theta 
+  \int_0^1 \int_{\Omega_h^\theta(t)} f (\nabla\cdot e_x^{\theta})  (\partial_{t,*}^{\bullet} e_u^*)^\theta  \, \d x \d\theta\Big| \\
&\le C \Big(\|e_x^*\|_{L^2(\Omega_h^*(t))} \| \nabla f\|_{L^\infty(\Omega_h^*(t))} +
\|\nabla e_x^*\|_{L^2(\Omega_h^*(t))} \| f\|_{L^\infty(\Omega_h^*(t))} \Big)\|\partial_{t, *}^{\bullet} 
e_u^*\|_{L^2(\Omega_h^*(t))} \notag\\
&\le C\|\e_\x\|_{\A(\x^*)} \|\dot \e_\u\|_{\M(\x)} \notag\\
&\le \varepsilon \|\dot \e_\u\|_{\Mrho(\x)}^2 + C \varepsilon^{-1} \|\e_\x\|_{\A(\x^*)}^2. \notag
\end{align}
}

The last term on the right hand side of \eqref{Error_Energy_u} can also be estimated by using integration by parts in time, i.e., 
\begin{align}
- (\dot \e_\u)^\top \EE_{\u1} 
= -\frac{\d}{\d t} \big( (\e_\u)^\top\EE_{\u1}  \big)
+ (\e_\u)^\top \Big( \frac{\d}{\d t} \EE_{\u1} \Big)
=: {S_{81} + S_{82} }. 
\end{align}
Then, using the consistency estimates in Proposition \ref{lemma:consistency}, we have 
\begin{subequations}\label{Estimate-S7}
\begin{align}
\Big| \int_0^t S_{81}(s)\, \d s \Big| 
&\le Ch^k \|\e_\u(t)\|_{\KK(\x^*)} 
 \le \varepsilon \|\e_\u(t)\|_{\Amu(\x^*)}^2 + C \varepsilon^{-1} h^{2k}, \\
|S_{82}| & \le Ch^{k}\|\e_\u\|_{\KK(\x^*)}  \le C\| \e_\u\|_{\A(\x^*)}^2 +   C h^{2k}  . 
\end{align}
\end{subequations}

Integrating the inequality \eqref{Error_Energy_u} in time from $0$ to $t$ and using the estimates of $S_j$, $j=1,\dots,6$, as well as the estimates in \eqref{Estimate-S7}, with a sufficiently small $\varepsilon$ in these estimates, the terms $\varepsilon \|\e_\u(t)\|_{\Amu (\x^*)}^2$ and $\varepsilon \int_0^t \|\dot \e_\u(s)\|_{\Mrho(\x)}^2\d s$ on the right-hand side can be absorbed by the left-hand side.
For $k \ge 2$ and sufficiently small $h$, the term $Ch^{\frac{k}{2}-\frac{d}{4}} \|\e_\u(t)\|_{\Amu (\x^*)}^2$ on the right-hand side can also be absorbed by the left-hand side. Then we obtain the following result: 
\begin{align}
\label{stability-eu-on-ep}
&\|\e_\u(t)\|_{\A(\x^*)}^2  + \int_0^t \|\dot \e_\u(s)\|_{\M(\x^*)}^2 \, \d s \\ 
&\le \varepsilon \int_0^t  \|\e_\p(s)\|_{\M(\x^*)}^2\, \d s + C\varepsilon^{-1} h^{2k} \notag\\
&\quad\, + C\varepsilon^{-1}\int_0^t \Big(\|\e_\x(s)\|_{\A(\x^*)}^2 + \|\e_\u(s)\|_{\A(\x^*)}^2 + \|\e_\w(s)\|_{\A(\x^*)}^2 \Big) \, \d s.
\notag
\end{align}
The first term on the right-hand side of \eqref{stability-eu-on-ep} is estimated in the next part of the proof.

\noindent (D) \textit{Estimates for $\e_\p$ under condition \eqref{ex-W1infty}}: 
To estimate $\|\e_\p\|_{\M(\x^*)}^2$, we test \eqref{Error_Eq_u} by a vector $\v$ associated to a finite element function  $v_h$ on $\Omega_h(t)$ to get 
\begin{align} \label{ep-test-dx-vh}
(\v)^\top \C(\x)^\top \e_\p 
& = (\v)^\top  \Mrho(\x)\dot \e_\u  +  (\v)^\top (\Mrho(\x)- \Mrho(\x^*)) \dot \u^* \\[4pt]
&\quad\,+ (\v)^\top \Amu(\x)\e_\u
+(\v)^\top (\Amu(\x)-\Amu(\x^*))\u^* \notag \\[4pt]
&\quad\, 
+ (\v)^\top  \big(\Brho(\x, \u -\w)\u - \Brho(\x^*, \u^* -\w^*)\u^*\big) \notag\\
&\quad\,
- (\v)^\top (\C(\x)^{\top}-\C(\x^*)^{\top})\p^* - {(\v)^\top (\M(\x) \f - \M(\x^*) \f^*)} + (\v)^\top  \EE_{\u1}. \notag
\end{align}

The estimates of the terms on the right-hand side of \eqref{ep-test-dx-vh} are similar as that in part (C) of this proof and therefore omitted. Specifically, we have 
\begin{align} \label{Gj}
|(\v)^\top \C(\x)^\top \e_\p |
&\le C\big(\|\dot \e_\u\|_{\M(\x)} + \|\e_\x\|_{\A(\x^*)}+\|\e_\u\|_{\A(\x^*)} + \|\e_\w\|_{\A(\x^*)} +  h^{k} \big)\|\v\|_{\KK(\x)} .
\end{align}

By using the inf-sup condition \eqref{inf-sup} in the expression for $\|\e_\p\|_{\M(\x)}$, we {can estimate} $\|\e_\p\|_{\M(\x)}$ as follows:
\begin{align} \label{ep-L2-1}
\|\e_\p\|_{\M(\x)}
&\le  \sup_{ \v \neq 0} \frac{(\v)^\top \C(\x)^\top \e_\p }{\|\v\|_{\A(\x)}}\\
&\le C\big(\|\dot \e_\u\|_{\M(\x)} + \|\e_\x\|_{\A(\x^*)}+\|\e_\u\|_{\A(\x^*)} + \|\e_\w\|_{\A(\x^*)} + h^k \big), \notag
\end{align}
which together with the norm equivalence in Remark~\ref{remark-norm-equivalence} gives 
\begin{align} \label{stability-ep-L2-}
\|\e_\p\|_{\M(\x^*)}^2
\le \, & C\big(\|\dot \e_\u\|_{\M(\x^*)}^2 + \|\e_\x\|_{\A(\x^*)}^2 +\|\e_\u\|_{\A(\x^*)}^2 + \|\e_\w\|_{\A(\x^*)}^2 + h^{2k} \big) .
\end{align} 
By integrating the above inequality with respect to time from $0$ to $t$, we have 
\begin{align} \label{stability-ep-L2}
\int_0^t\|\e_\p(s)\|_{\M(\x^*)}^2\, \d s
\le \, & Ch^{2k}  +  C\int_0^t \|\dot \e_\u(s)\|_{\M(\x^*)}^2\, \d s \\
&+ C\int_0^t \Big( \|\e_\x(s)\|_{\A(\x^*)}^2 +\|\e_\u(s)\|_{\A(\x^*)}^2 + \|\e_\w(s)\|_{\A(\x^*)}^2 
\Big)\, \d s . \notag
\end{align} 

\noindent {(E) \textit{Combination of the estimates}:}
Summing up inequalities \eqref{stability-ex-H1}, \eqref{stability-eu-on-ep} and $\varepsilon\times$\eqref{stability-ep-L2}, substituting estimate \eqref{stability-ew-H1} into the resulted inequality and choosing a sufficiently small $\varepsilon$, we obtain the following combined estimate: 
\begin{align}\label{stability-combination-2}
&\|\e_\u(t)\|_{\A(\x^*)}^2  + \|\e_\x(t)\|_{\A(\x^*)}^2 + \int_0^t \|\dot \e_\u(s)\|_{\M(\x^*)}^2 \, \d s \\
&\le Ch^{2k} + C\int_0^t \Big(\|\e_\x(s)\|_{\A(\x^*)}^2 + \|\e_\u(s)\|_{\A(\x^*)}^2  \Big) \, \d s. \notag
\end{align}
By applying Gronwall's inequality, we obtain 
\begin{align} \label{stability-eu-H1-final}
\sup_{0 \le t \le t^*} \Big(\|\e_\u(t)\|_{\KK(\x^*)}^2 +  \|\e_\x(t)\|_{\KK(\x^*)}^2
+ \int_0^t \|\dot \e_\u(s)\|_{\M(\x^*)}^2\, \d s \Big)
\le Ch^{2k}.
\end{align}
This result can be substituted into \eqref{stability-ew-H1} and \eqref{stability-ep-L2} to yield an error estimate for $\e_\w$ and $\e_\p$, i.e., 
\begin{align} \label{stability-ew-H1-final}
\sup_{0 \le t \le t^*} \Big( \|\e_\w(t)\|_{\KK(\x^*)}^2  + 
\int_0^t \| \e_\p(s)\|_{\M(\x^*)}^2\, \d s \Big)\le Ch^{2k}.
\end{align}

It remains to show that $t^* = T$ when $h$ is sufficiently small (smaller than some constant). In fact, for $0 \le t \le t^*$ we can use the inverse inequality and \eqref{stability-eu-H1-final} to bound the left-hand side in \eqref{ex-W1infty-a} as
\begin{align}
\|\nabla e_x^*(t)\|_{L^{\infty}(\Omega_h^*(t))} \le  Ch^{-\frac{d}{2}}\|\e_{\x}(t)\|_{\KK(\x^*)}
\le  Ch^{k -\frac{d}{2}},
\end{align}
where $C$ is a constant dependent on $T$, but independent of $t^*$. Since $k \ge 2$, we have the following
inequality for $0 \le t \le t^*$ when $h$ is sufficiently small:
\begin{align}
\|\nabla e_x^*(t)\|_{L^{\infty}(\Omega_h^*(t))} \le \frac12 h^{-\frac{d}{4}}h^{\frac{k}{2}}.
\end{align}
For the same reason, for $k \ge 2$ and sufficiently small $h$, we have
\begin{align}
\|\nabla e_u(t)\|_{L^{\infty}(\Omega_h^*(t))}  \le \frac12 \quad \mbox{and} \quad \|\nabla e_w(t)\|_{L^{\infty}(\Omega_h^*(t))}  \le \frac12 \quad \mbox{for} \quad 0\le t\le t^*.
\end{align}
Hence, we can extend the bound \eqref{ex-W1infty-a}--\eqref{ex-W1infty-c} beyond  $t^*$, which contradicts the maximality of $t^*$ unless $t^* = T$. This proves that the error estimates in \eqref{stability-eu-H1-final}--\eqref{stability-ew-H1-final} actually hold for $t^*=T$. 
Since the flow map $\phi_h^*(\cdot, t)$ and its inverse {are bounded} in the $W^{1,\infty}$ norm, it follows that error estimates can be pulled back to the initial domain $\Omega_h^0$ after composition with $\phi_h^*(\cdot, t)$. This proves the desired error estimates in Theorem \ref{Main-THM}.
\hfill
%

\section{Numerical results}\label{section:numerical-test}
{ In this section we present numerical results for a benchmark example without surface tension. One of the most well-known benchmark examples for the two-phase NS flow problem was proposed in \cite[Table I]{hysing2009quantitative} for simulating a moving bubble in fluid.
This example includes surface tension and has been} widely adopted to test the convergence and robustness of numerical methods for two-phase NS flows; see \cite{aland2012benchmark, barrett2015stable, fu2020arbitrary}. {In our study, we modify this model by considering a simplified version without surface tension and with no-slip boundary conditions.}
All the computations are performed using the software package NGSolve ({\tt https://ngsolve.org}), {and the mesh is generated by Gmsh ({\tt https://gmsh.info/})}. 

The computational domain is $\Omega = (0, 1)\times(0, 2)$ and the no-slip pieces of the boundary are given by
$$
\Gamma_1 =[0,1]\times \{0,2\} \quad \mbox{and} \quad \Gamma_2 = \{0,1\}\times[0,2].
$$
The initial interface is a circle with radius 0.25, i.e.,
$$
\Gamma(0) = \big\{(x, y) : |(x, y) - (0.5, 0.5)| = 0.25 \big \}.
$$
The gravitational force is $f = (0, - 0.98)^\top$, and the physical parameters are shown in Table~\ref{tab_parameters}.
\begin{table}[htp!]
\centering
\caption{\it Physical parameters for the benchmark problem}
\setlength{\tabcolsep}{5mm}{
\begin{tabular}{ccccc}
\toprule
& $\rho_+$ & $\rho_-$ & $\mu_+$ & $\mu_-$ \\
\cmidrule(r){2-5}
BP-1 & $10^3$ & $10^2$ &$10$ &$1$ \\
\bottomrule
\end{tabular}}
\label{tab_parameters}
\end{table}

We first present numerical results to illustrate the convergence rate of the proposed method \eqref{fem-ALE-weak} for BP-1 in Table \ref{tab_parameters}. {We apply the time discretization of the proposed method \eqref{fem-ALE-weak} with a linearly semi-implicit Euler method, which only requires solving several decoupled linear systems at every time level.}
Since the exact solution of the problem is unknown, we evaluate the convergence rate in space by the formula
$$
\mbox{conv. rate in space} = 
\log\Big(\frac{\|u_{\tau, h} -u_{\tau, h/m}\|}{\|u_{\tau, h/m}-u_{\tau, h/m^2}\|}\Big)/\log(m),
$$
based on the finest three meshes, where $u_{\tau, h}$ denotes the numerical solution at $t_N=T$ computed by using a stepsize $\tau$ and mesh size $h$.

\begin{figure}[htp]
\centerline{
\includegraphics[width=3.0in]{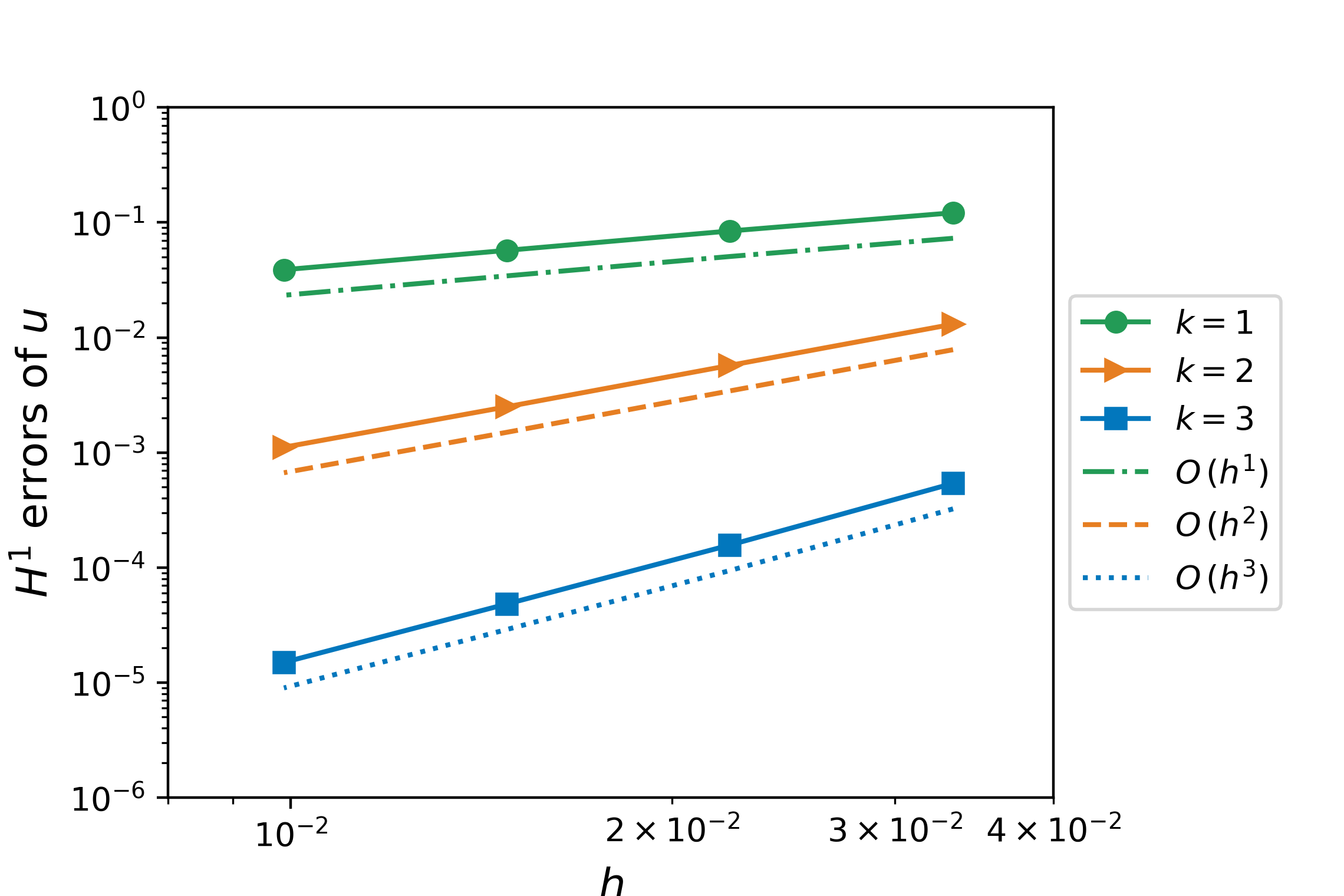}
\hspace{-4pt}
\includegraphics[width=3.0in]{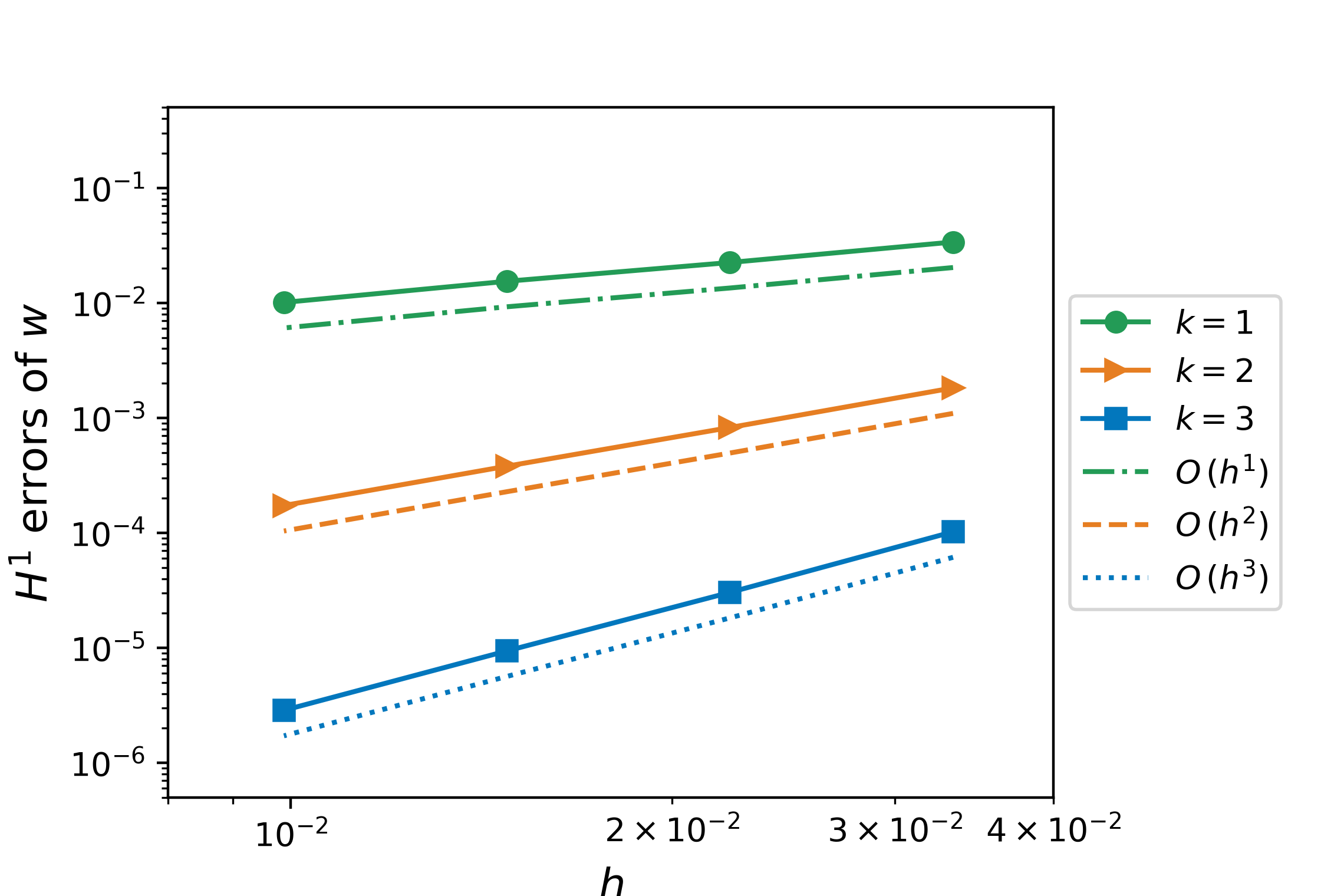}
}
\centerline{
\includegraphics[width=3.0in]{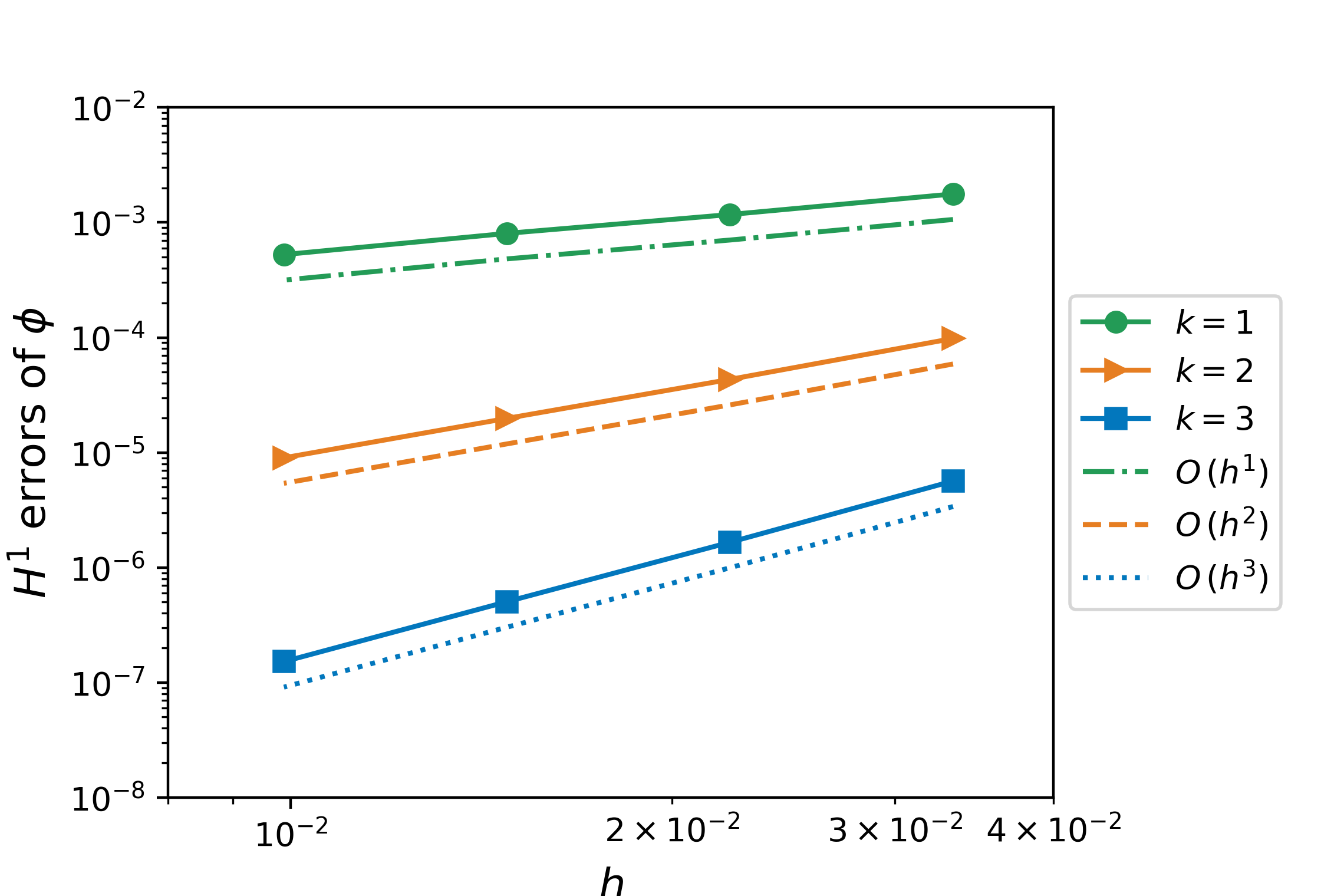}
\hspace{-4pt}
\includegraphics[width=3.0in]{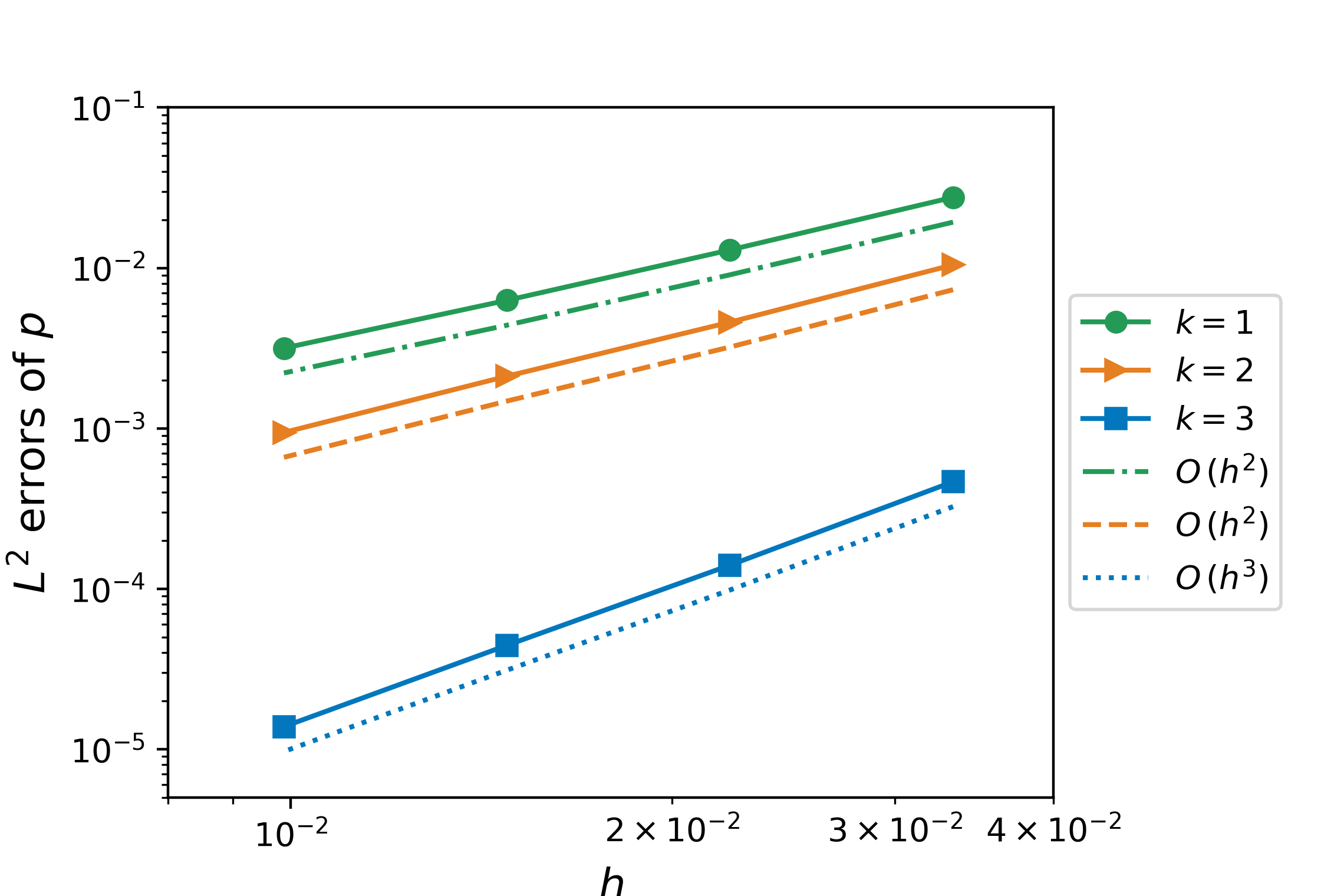}
}
\vspace{-6pt}
\caption{\it Spatial discretization errors of the numerical solutions at $T = 0.1$.}
\label{fig_errors}
\centerline{
\includegraphics[width=3.0in]{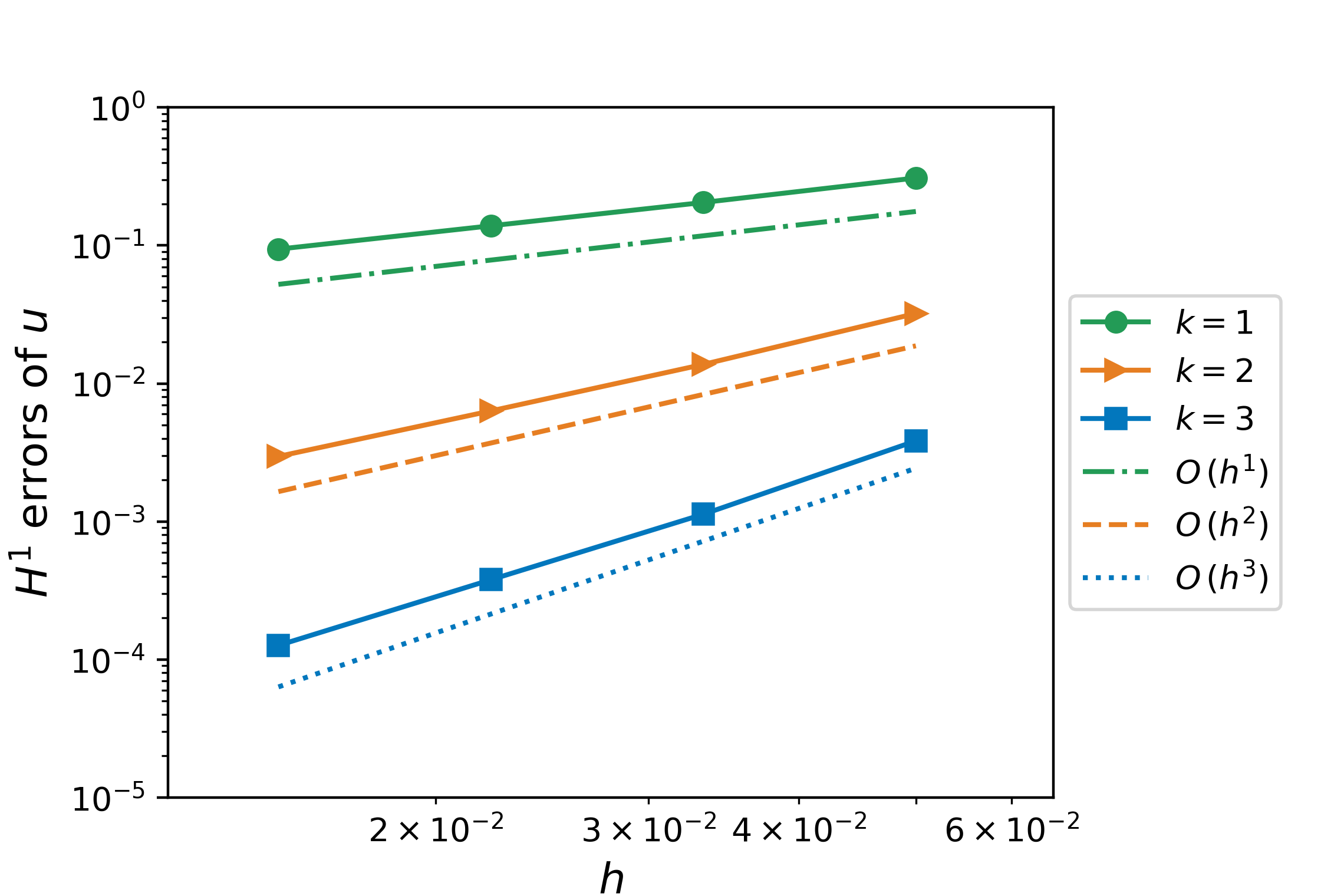}
\hspace{-4pt}
\includegraphics[width=3.0in]{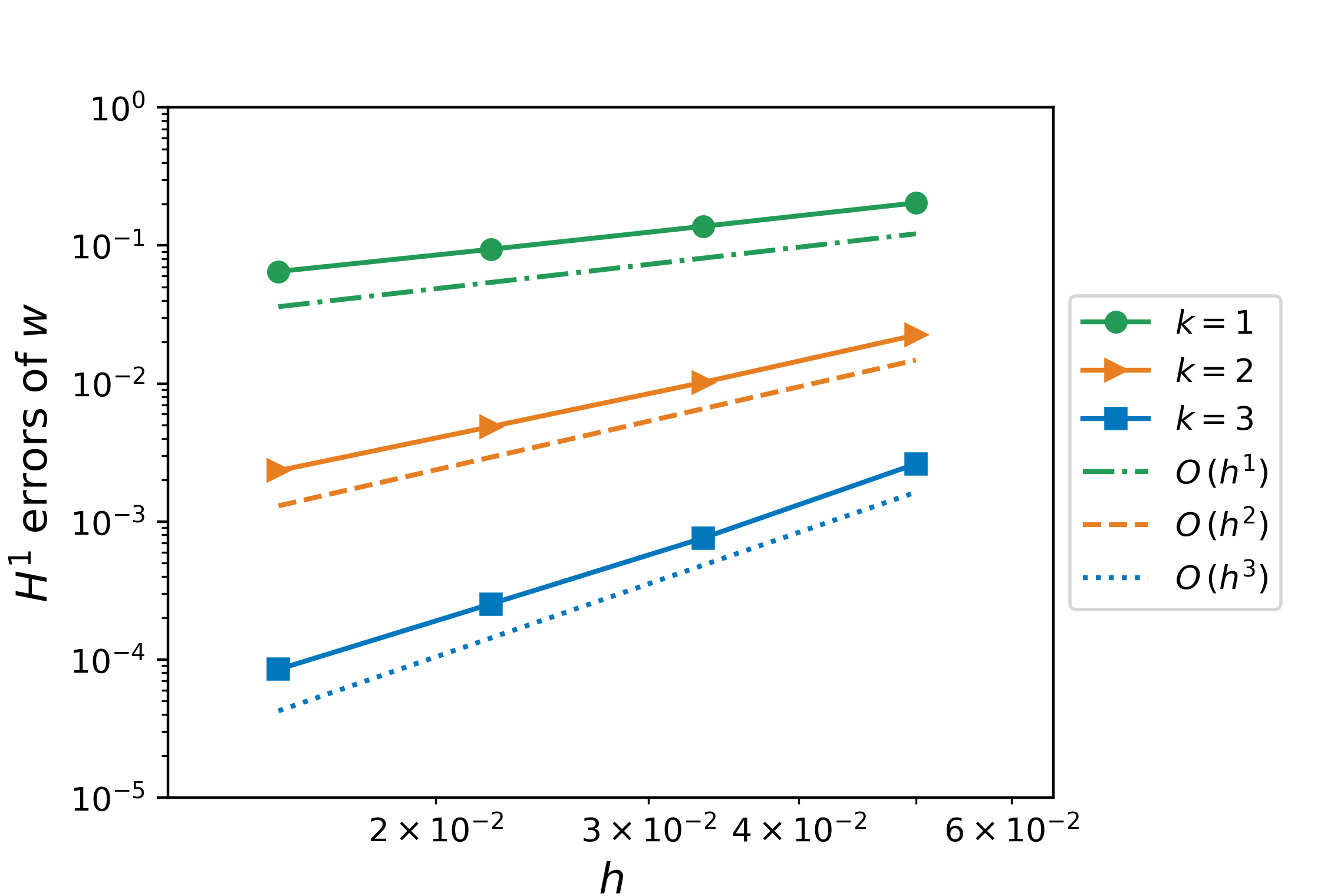}
}
\centerline{
\includegraphics[width=3.0in]{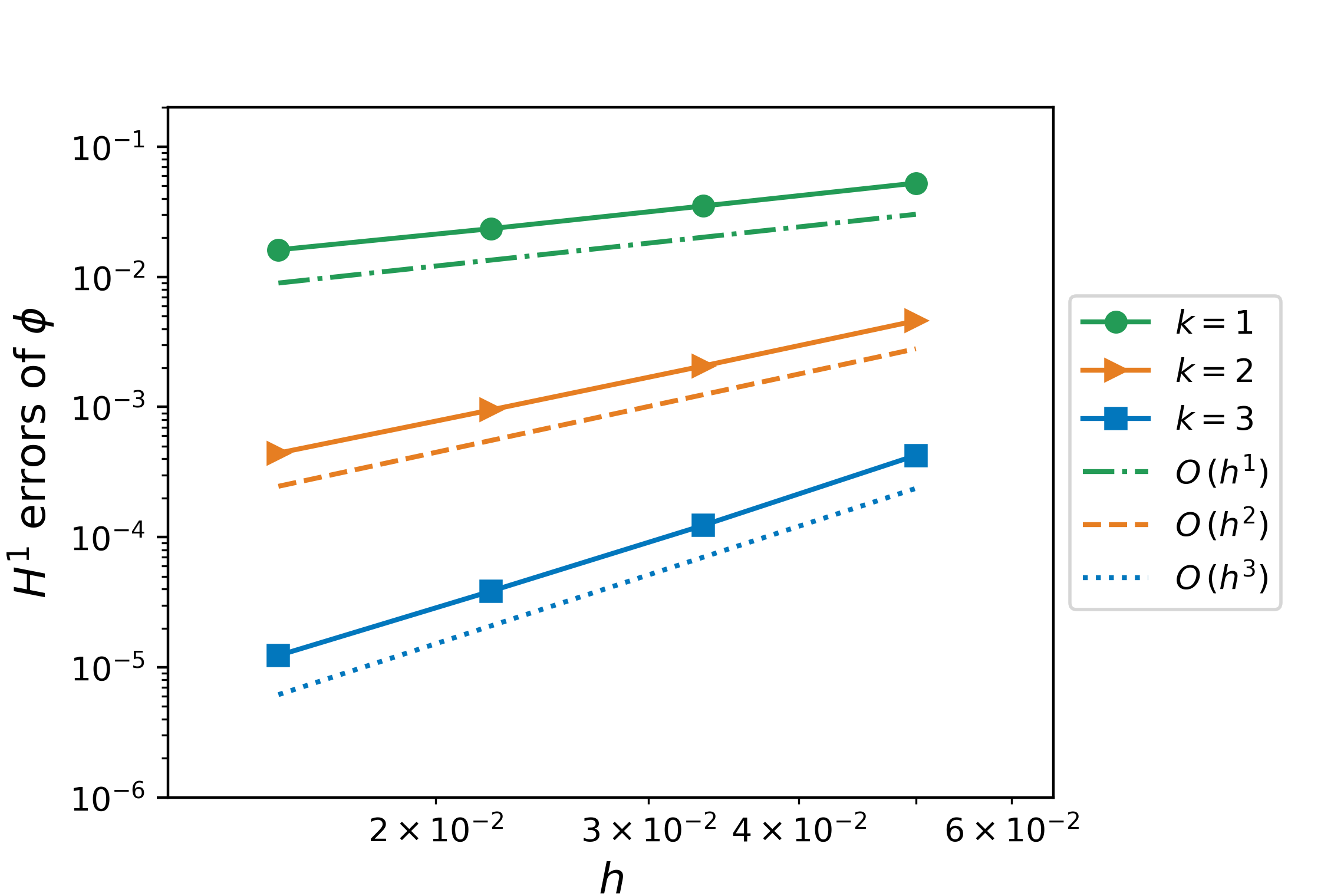}
\hspace{-4pt}
\includegraphics[width=3.0in]{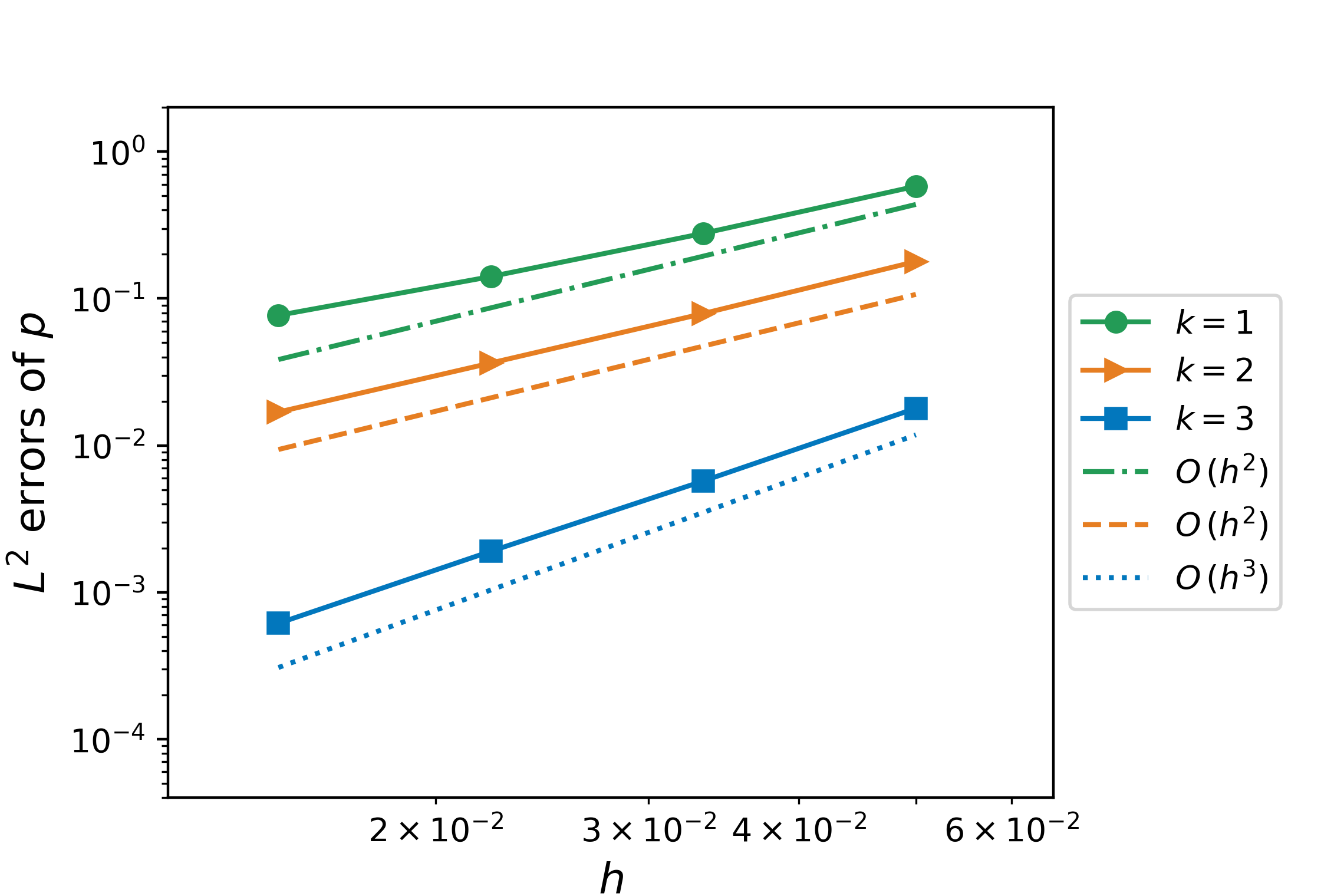}
}
\vspace{-6pt}
\caption{\it Spatial discretization errors of the numerical solutions at $T = 0.5$.}
\label{fig_errors_t5}
\end{figure}

The spatial discretization errors of the numerical solutions at $T = 0.1$ {and $ T = 0.5 $} are presented in Figure~\ref{fig_errors} {and Figure~\ref{fig_errors_t5}, respectively}, where we have used the following finite elements
\begin{align*}
&k = 1:\,\,\mbox{(Mini element $P_{1b}$--$P_1$)}\\ 
&k = 2:\,\, \mbox{(Taylor--Hood $P_2$--$P_1$)} \qquad\mbox{and}\qquad
k = 3: \,\,\mbox{(Taylor--Hood $P_3$--$P_2$)} 
\end{align*}
{with sufficiently small time step size} $\tau = 10^{-3}$ so that the temporal discretization errors are negligibly small compared to the spatial errors. 
From {Figures~\ref{fig_errors}--\ref{fig_errors_t5}}, we see that the error of space discretization is $O(h^{k})$, which is consistent with the theoretical results proved in Theorem~\ref{Main-THM}.

Next, we present simulations of the bubble's circularity, center of mass, rise velocity, and the total energy of the two-phase NS system with parameters BP-1 in Table \ref{tab_parameters}. We focus on the following benchmark quantities in \cite{hysing2009quantitative} with respect to time. 
\begin{itemize}
\item[(a)] {\it Center of mass}: 
\begin{align*}
X_d = \frac{1}{|\Omega_{-}|} \int_{\Omega_{-}} x_d \, \d x.
\end{align*}
\item[(b)] {\it Circularity} (referred as {\it sphericity} in 3-D): 
In two dimensions it is defined as the quotient between {\it perimeter of area-equivalent circle}
and {\it perimeter of the bubble}.
In three dimensions it is defined as the quotient between {\it surface area of volume-equivalent
ball} and {\it surface area of the bubble}, i.e., 
\begin{align*}
Cir = \frac{2\sqrt{|\Omega_{-}| \pi}}{|\Gamma|}  \quad(\mbox{in 2-D}) \qquad \mbox{and}
\qquad Cir = \frac{ \pi^{1/3} |\Omega_{-}|^{2/3}}{|\Gamma|} \quad(\mbox{in 3-D}).
\end{align*}
\item[(c)] {\it Rise velocity}:
\begin{align*}
V_d = \frac{1}{|\Omega_{-}|} \int_{\Omega_{-}} u_h \cdot e_d \, \d x,
\end{align*}
{where the integrand $ u_h \cdot e_d$ refers to the vertical component of the velocity $u_h$.}
\end{itemize}

\begin{figure}[htp]
\centerline{
\includegraphics[width=1.7in]{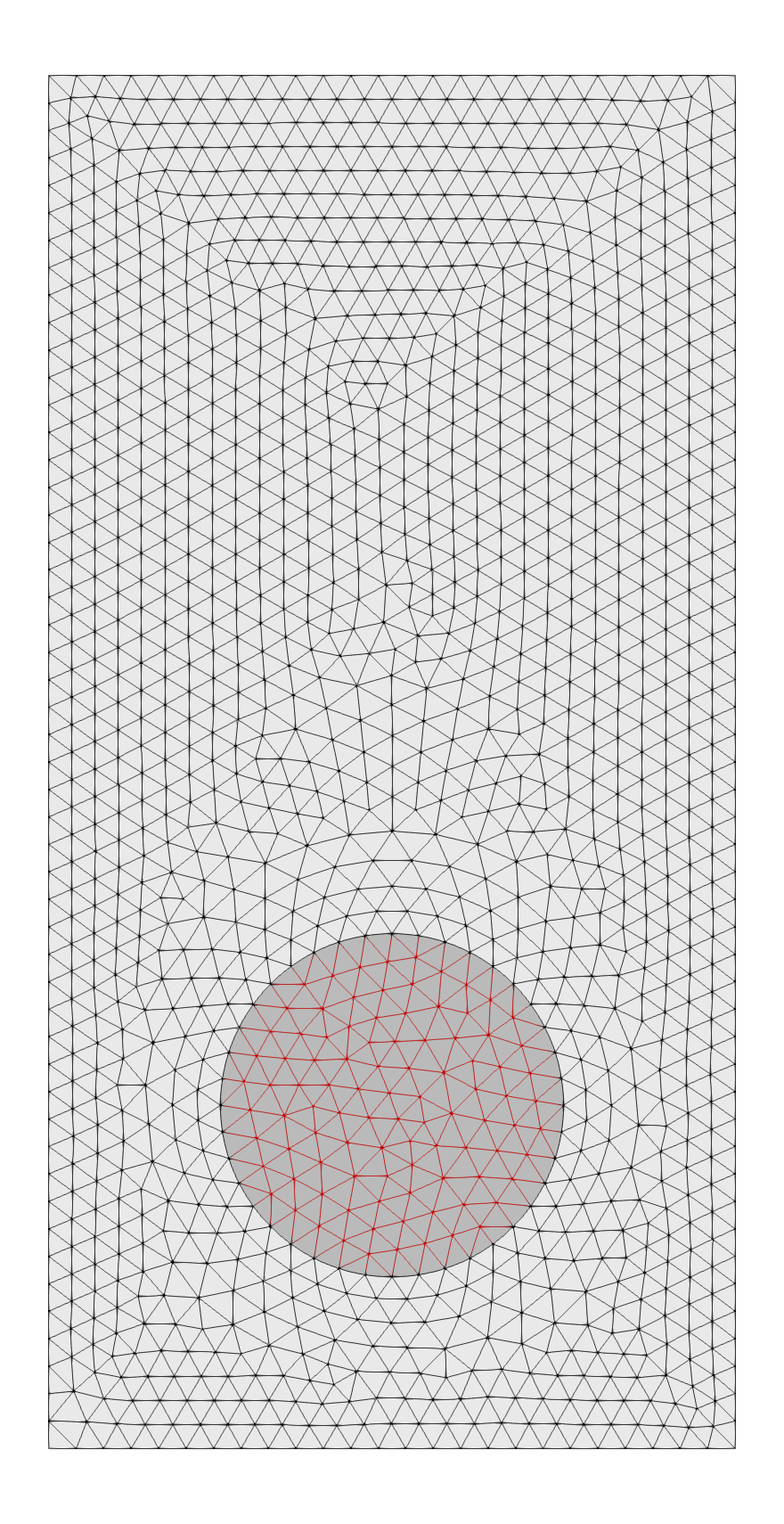}
\includegraphics[width=1.7in]{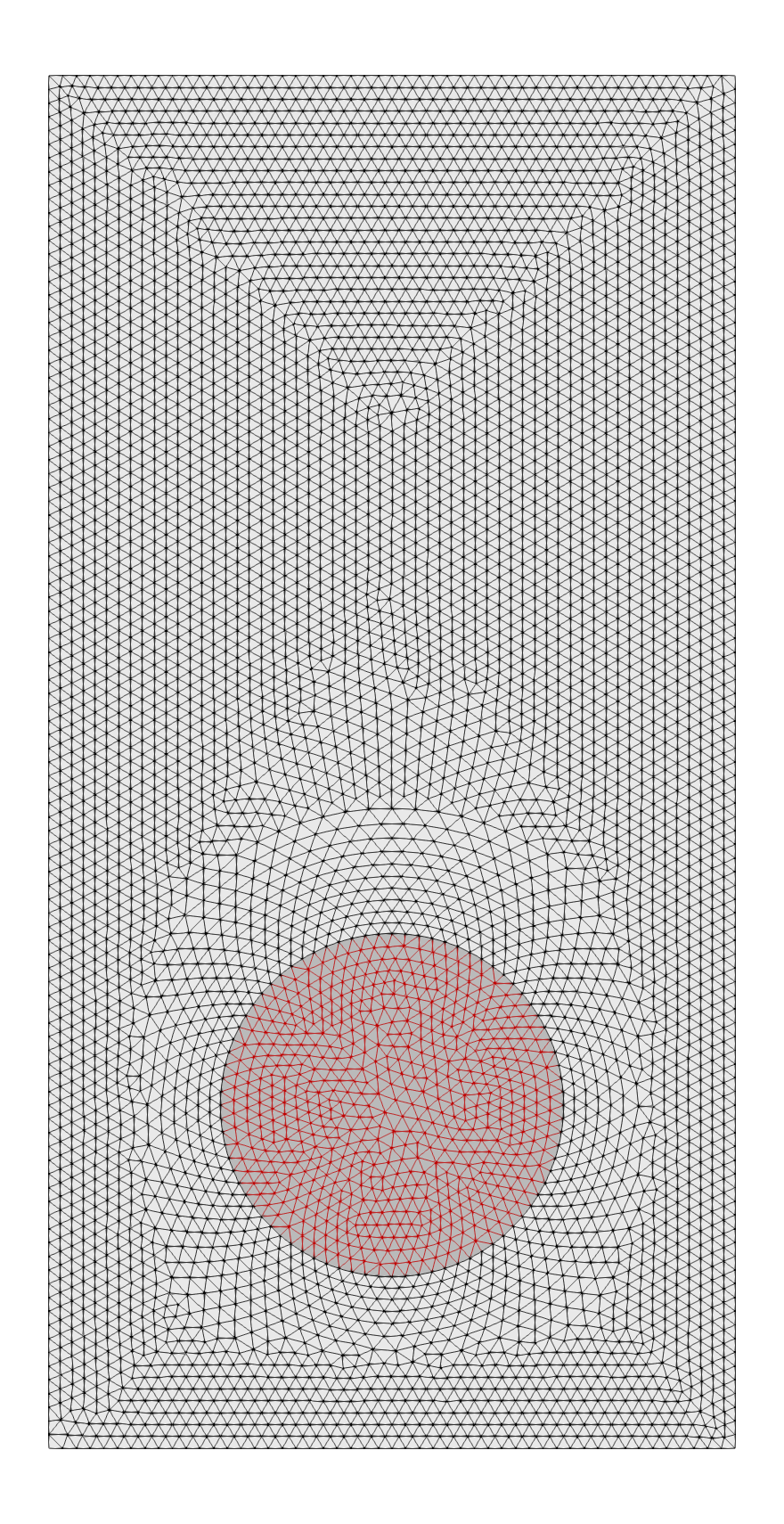}
}
\vspace{-8pt}
\caption{\it Left: mesh $M_1$. Right: mesh $M_2$. Dark color: domain $\Omega_{-}(0)$.}
\label{fig_three_mesh}
\vspace{20pt}
\centerline{
\includegraphics[width=1.5in]{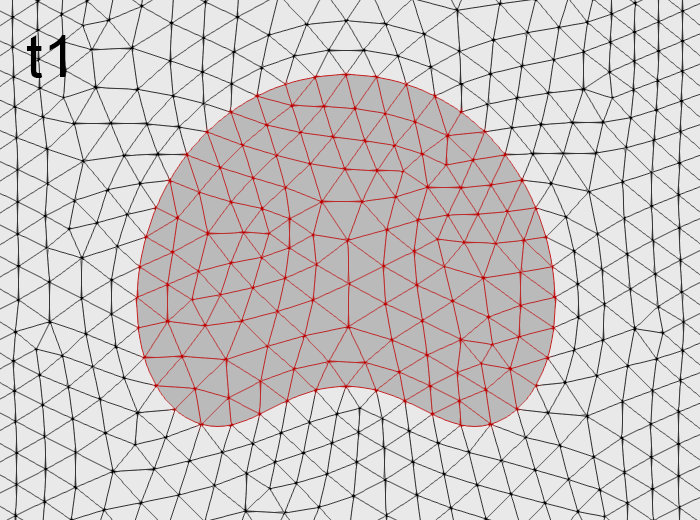}
\hspace{8pt}
\includegraphics[width=1.5in]{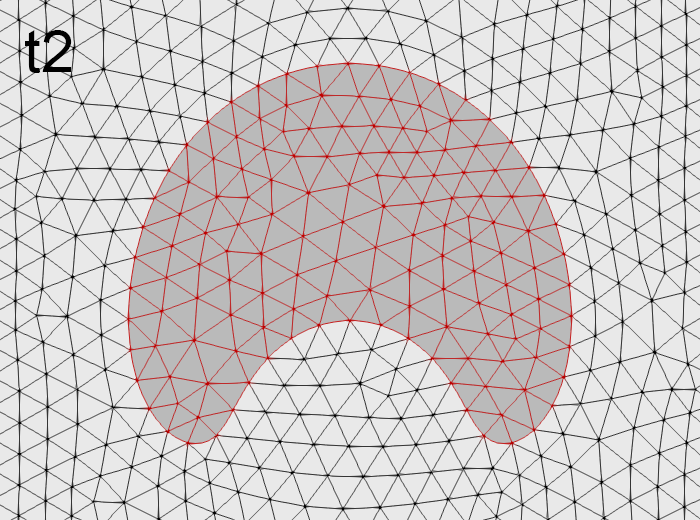}
\hspace{8pt}
\includegraphics[width=1.5in]{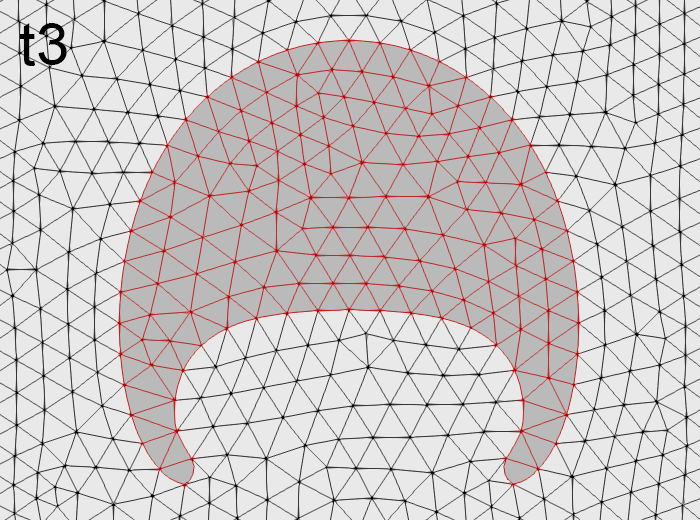}
}
\vspace{4pt}
\centerline{
\includegraphics[width=1.5in]{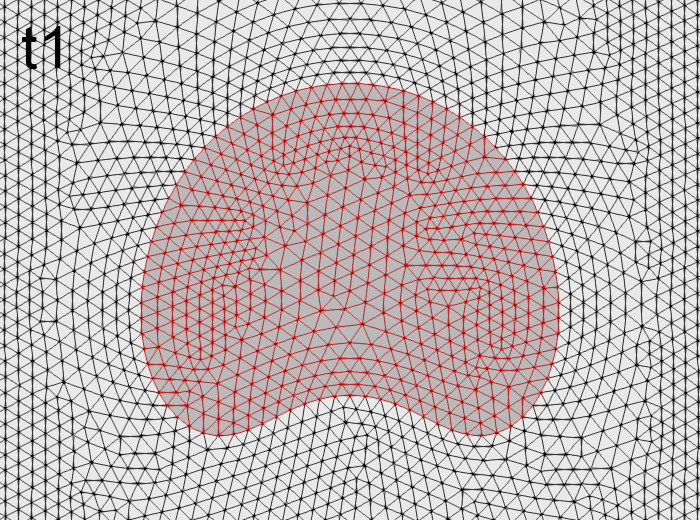}
\hspace{8pt}
\includegraphics[width=1.5in]{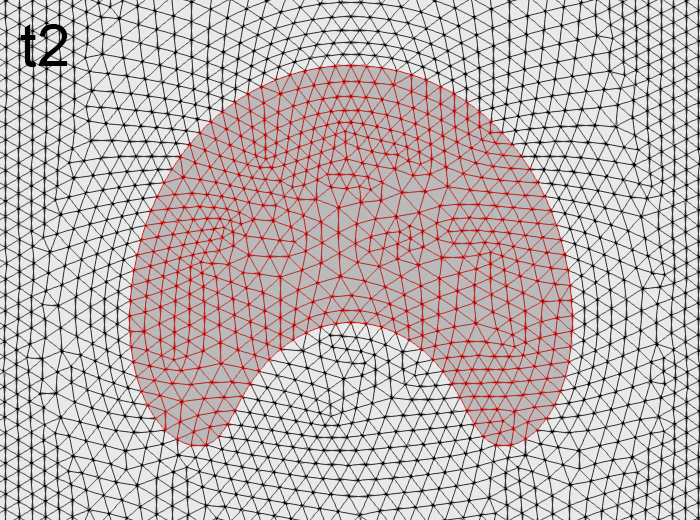}
\hspace{8pt}
\includegraphics[width=1.5in]{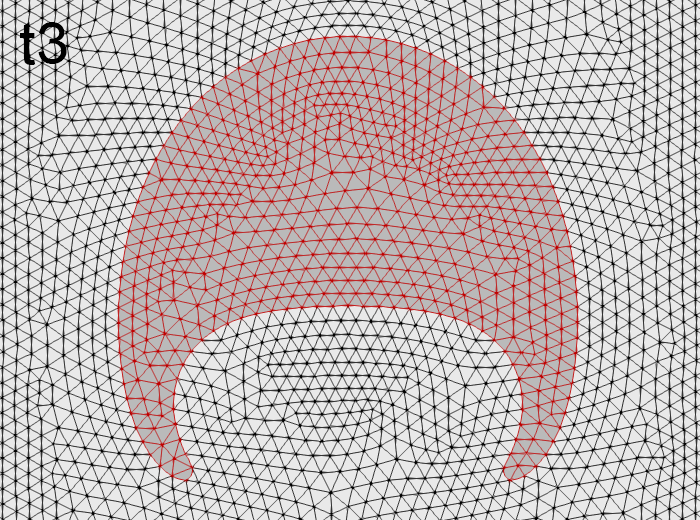}
}
\caption{\it Zoom-in of the deformed mesh around $\Omega_{-}(t)$ {with $ k=2 $}. Top: $M_1$. Bottom: $M_2$. Left: $t_1$ = 1. Middle: $t_2$ = 1.5. Right: $t_3$ = 2.5.}
\label{fig_zoom_in}
\end{figure}

{We consider the proposed method \eqref{fem-ALE-weak} with polynomial degree $k = 2$ and quadratic isoparametric triangular elements} on the two different meshes as depicted in Figure~\ref{fig_three_mesh}, where the mesh $M_1$ has mesh size $h_1 = 0.04$ and the mesh $M_2$ has mesh size $h_2 = 0.02$. Zoom-in of the deformed mesh around $\Omega_{-}(t)$ at $t = 1$, $1.5$ and $2.5$  are shown in Figure~\ref{fig_zoom_in}. 

\begin{figure}[htp]
\centerline{
\includegraphics[width=2.8in]{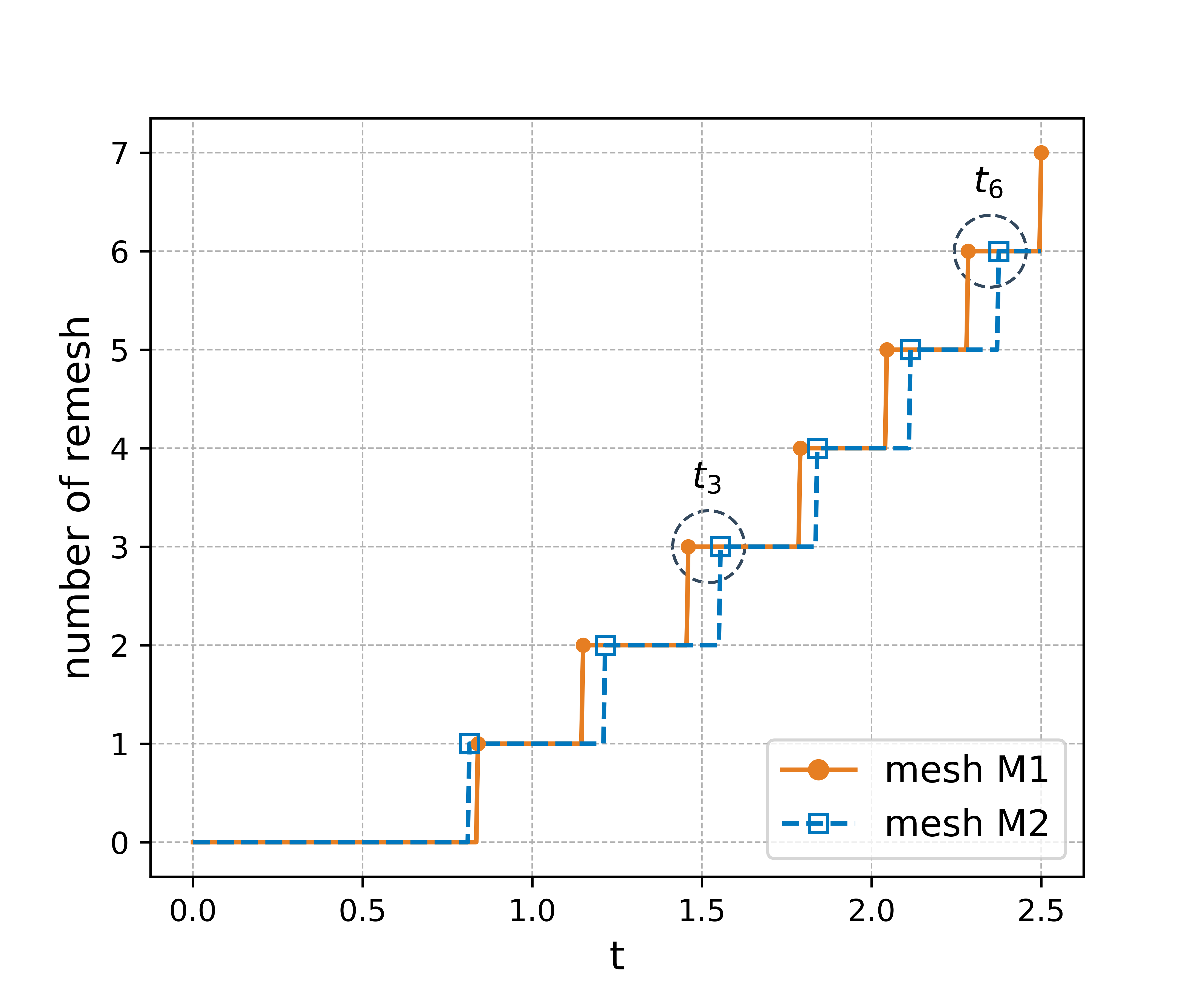}
}
\vspace{-8pt}
\caption{\it The number of remeshing steps with respect to time.}
\label{fig_remesh_number}
\vspace{10pt}
\centerline{
\includegraphics[width=1.45in]{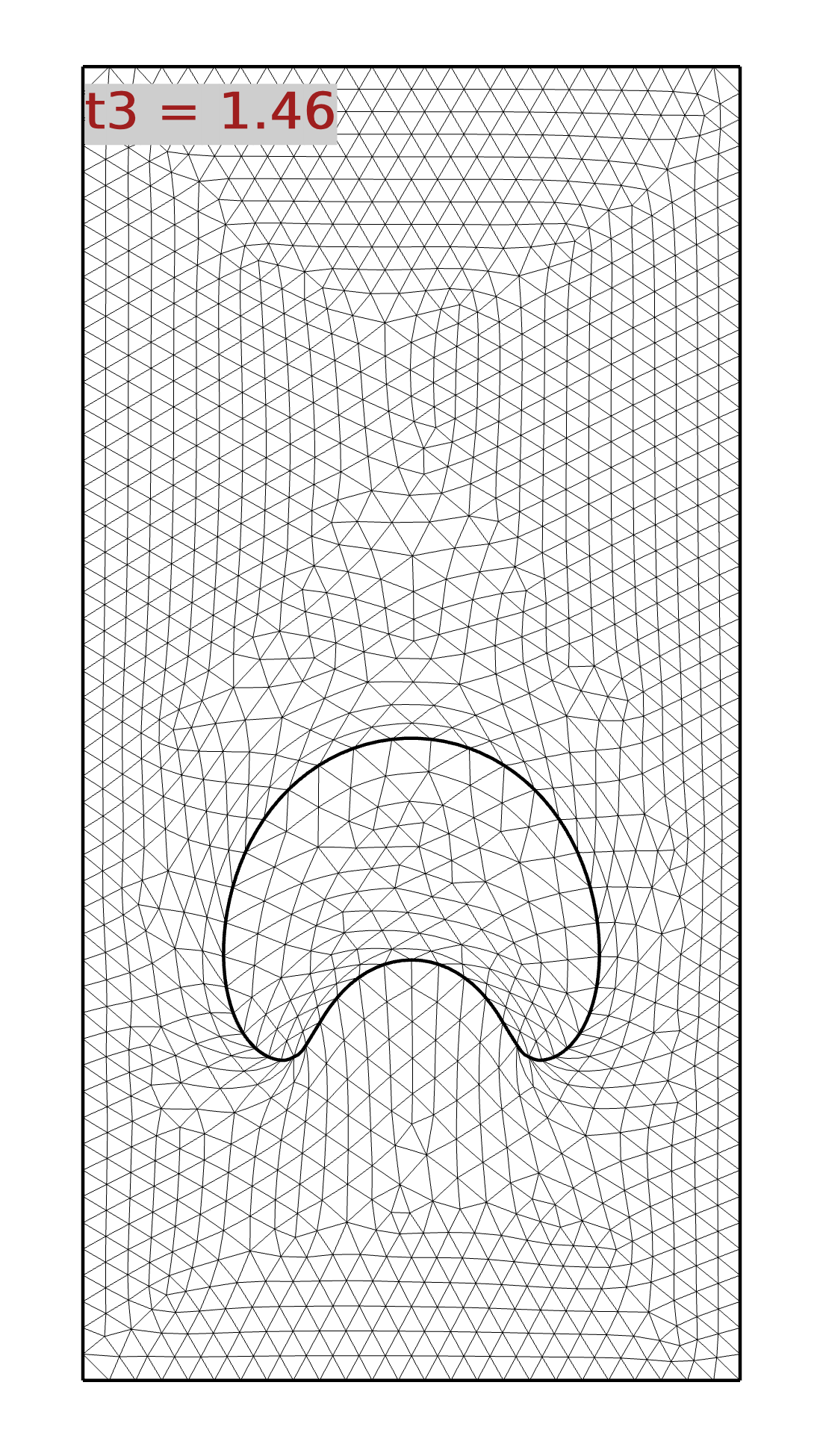}
\hspace{-12pt}
\includegraphics[width=1.45in]{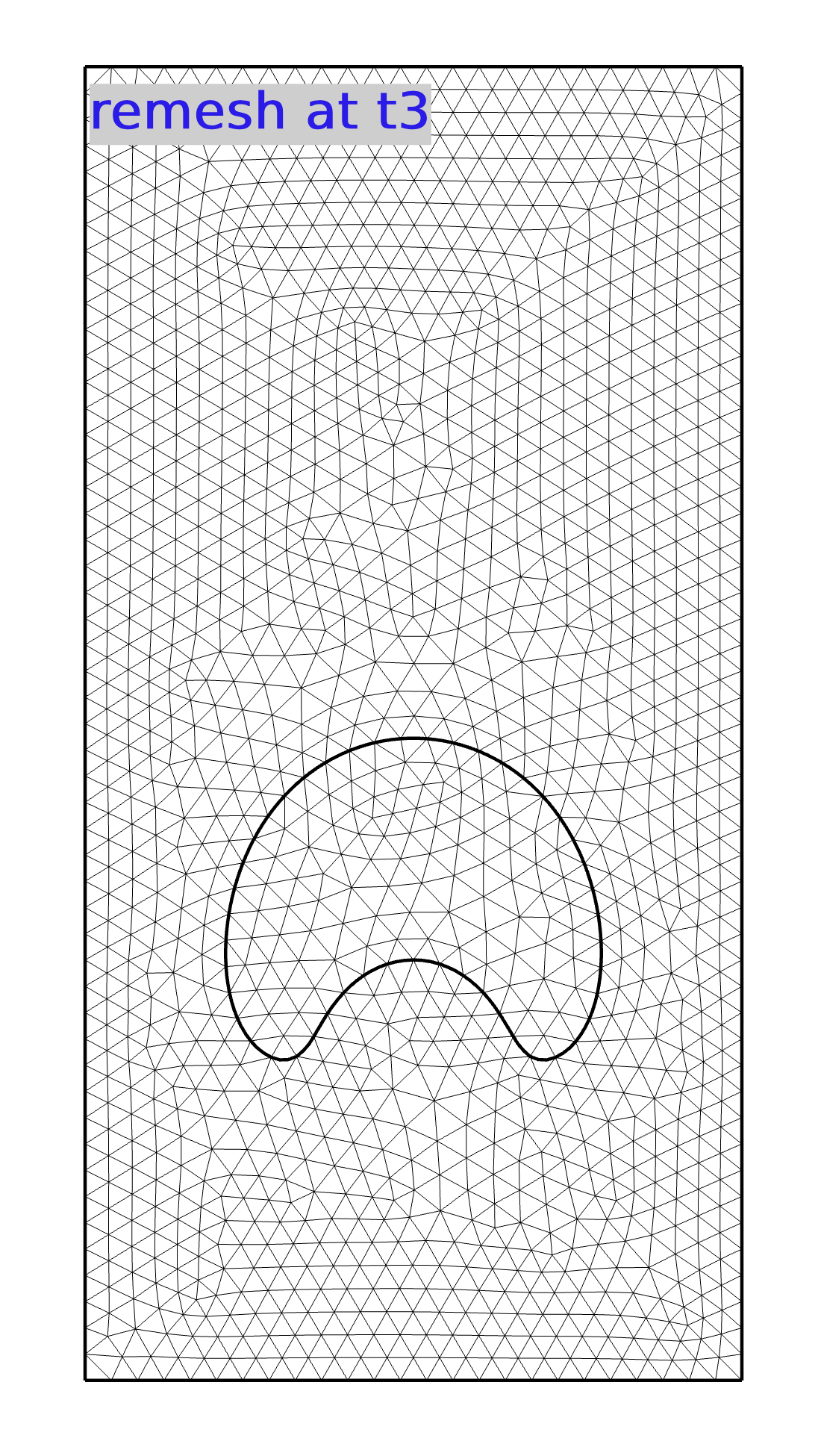}
\hspace{-12pt}
\includegraphics[width=1.45in]{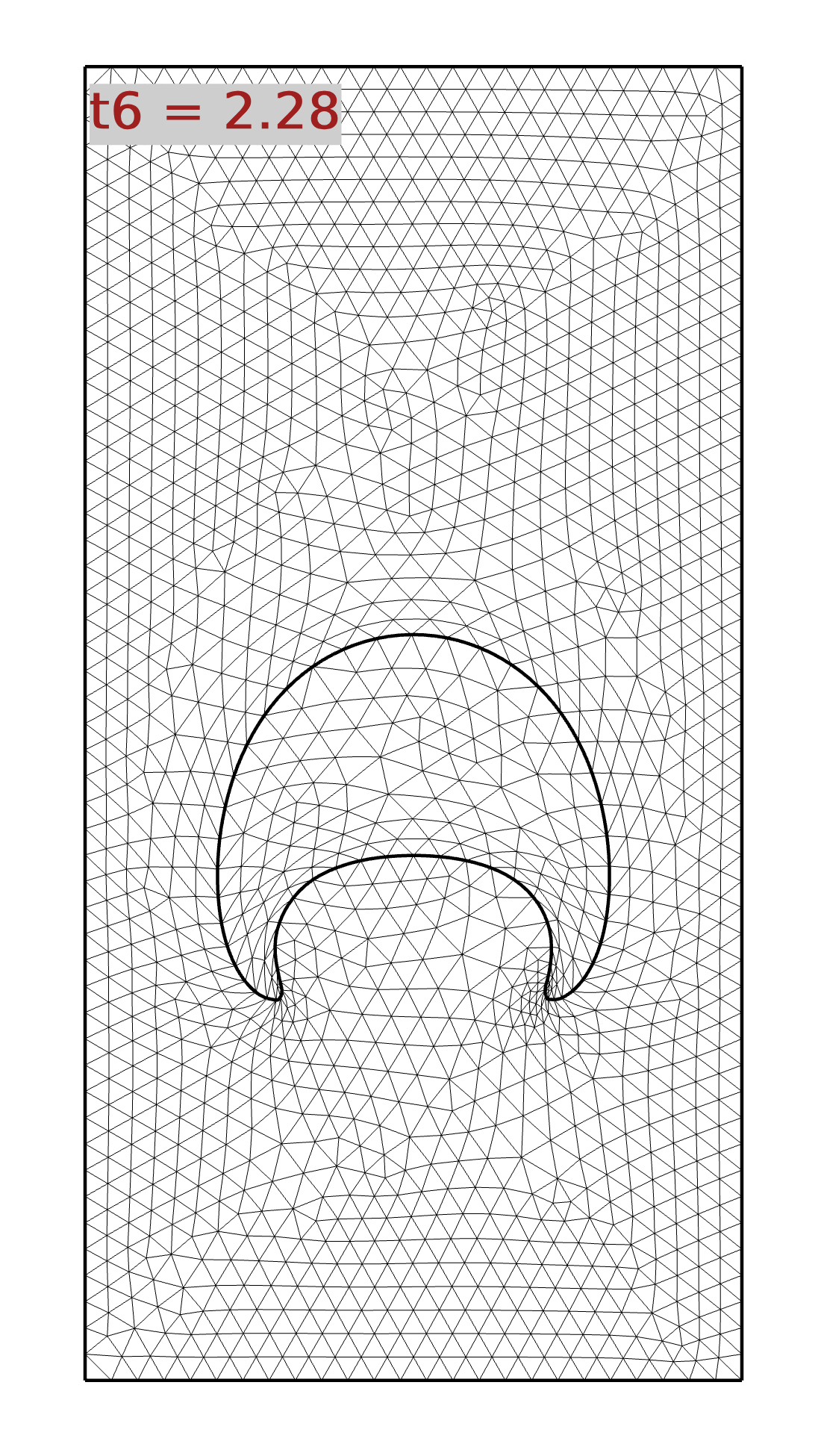}
\hspace{-12pt}
\includegraphics[width=1.45in]{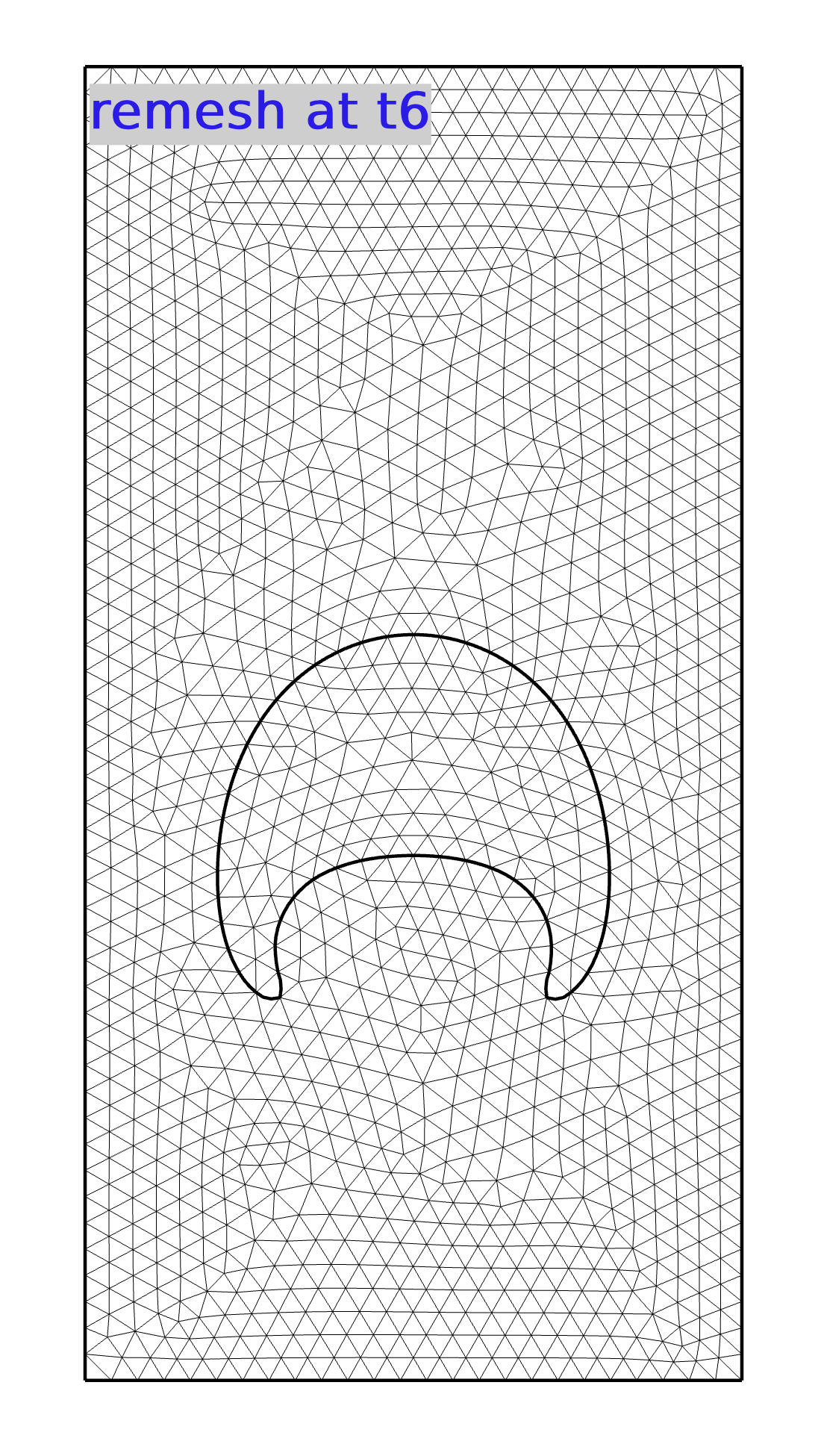}
}
\vspace{-10pt}
\caption{\it Meshes before and after a remeshing step for mesh $M_1$ with $k=2$. }
\label{fig_remesh_compare_M1}
\vspace{10pt}
\centerline{
\includegraphics[width=1.45in]{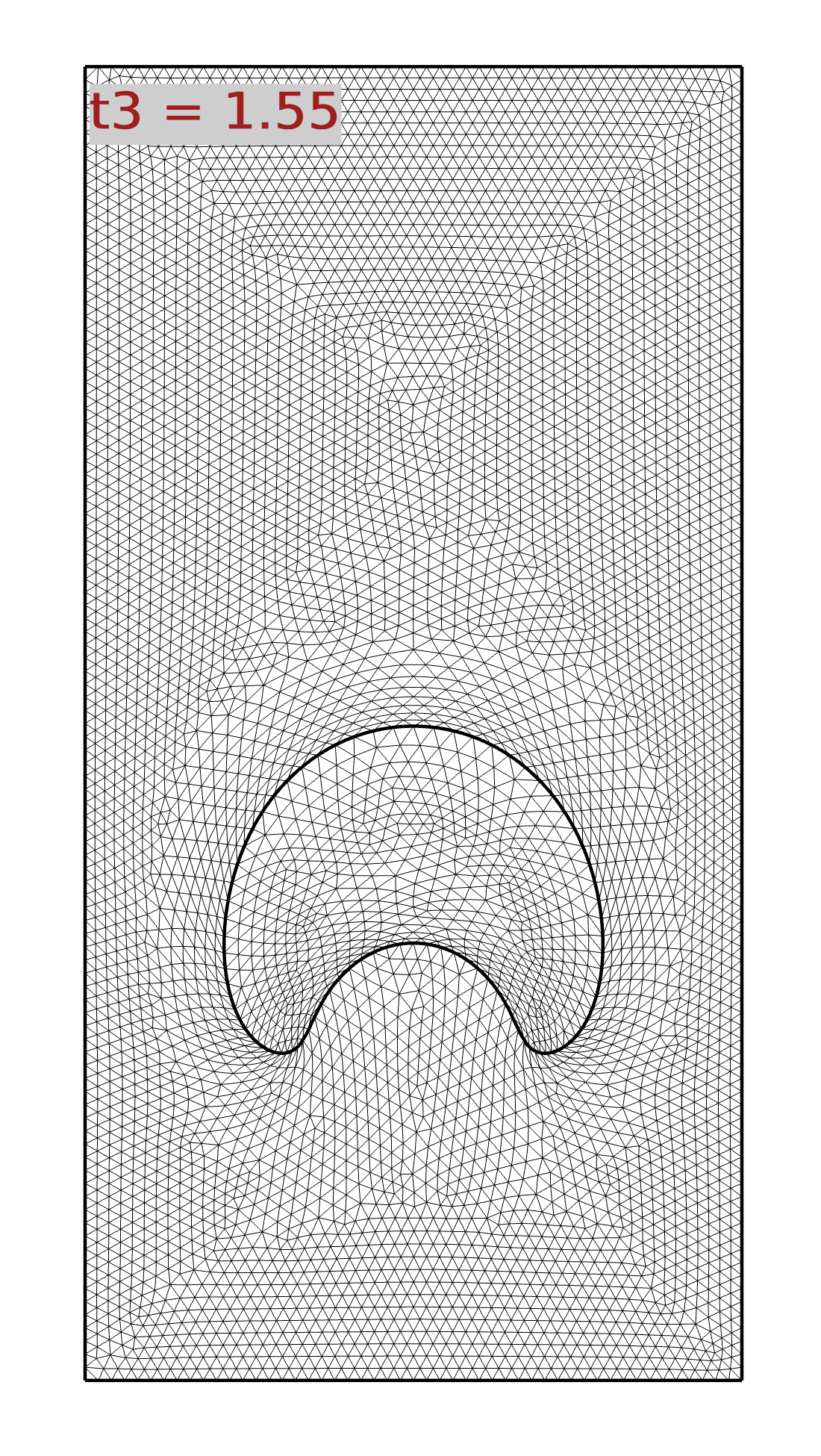}
\hspace{-12pt}
\includegraphics[width=1.45in]{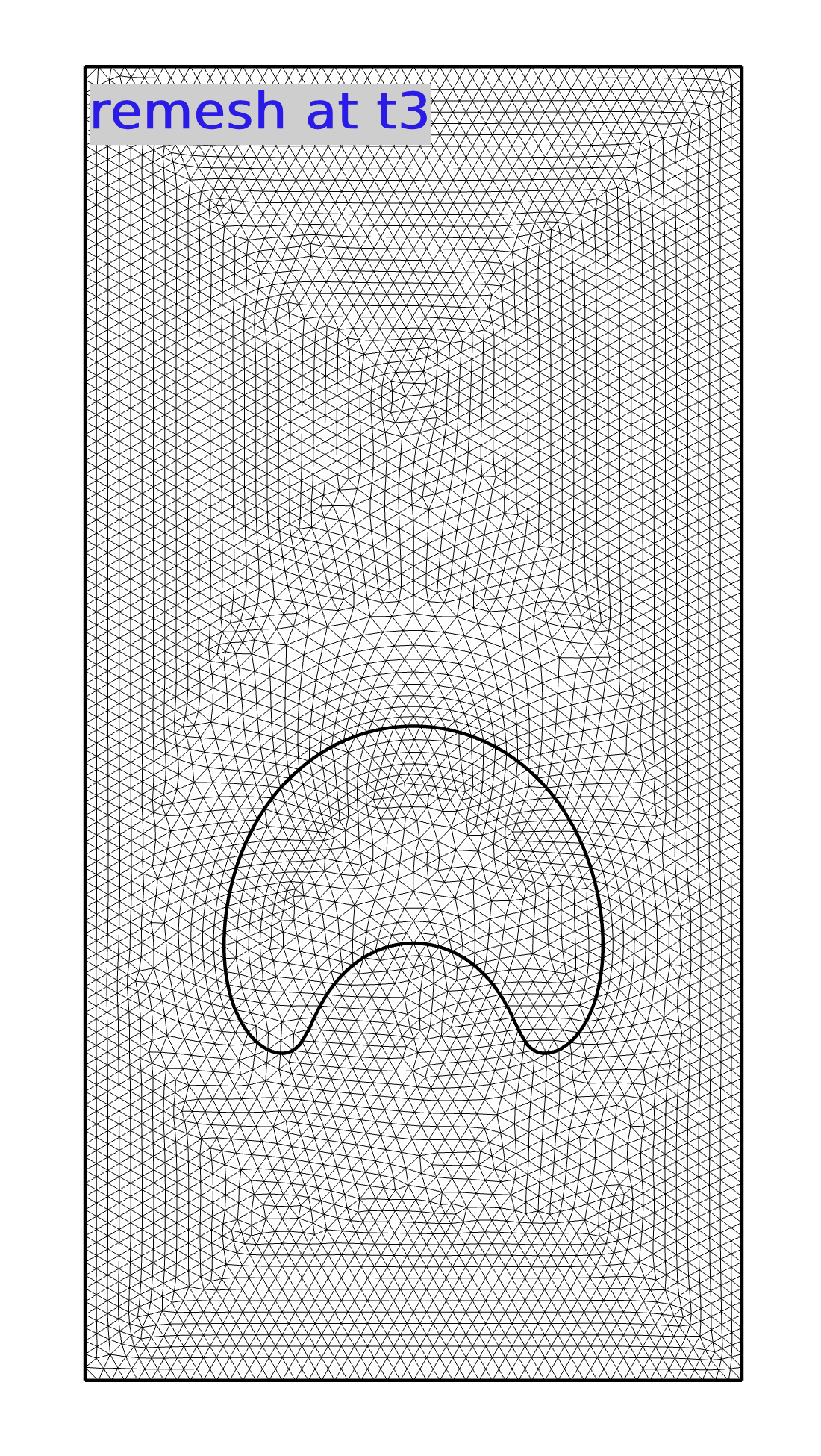}
\hspace{-12pt}
\includegraphics[width=1.45in]{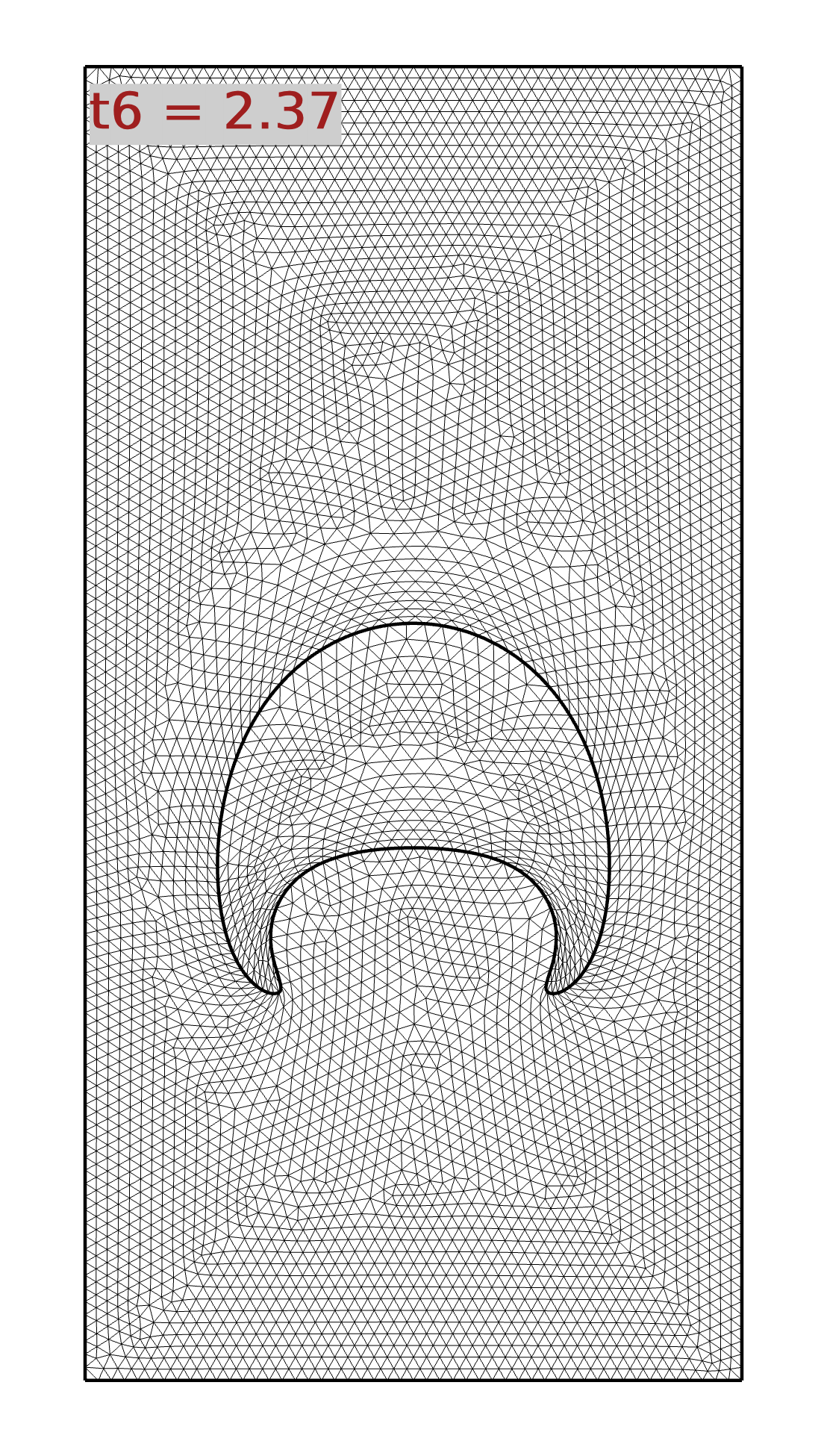}
\hspace{-12pt}
\includegraphics[width=1.45in]{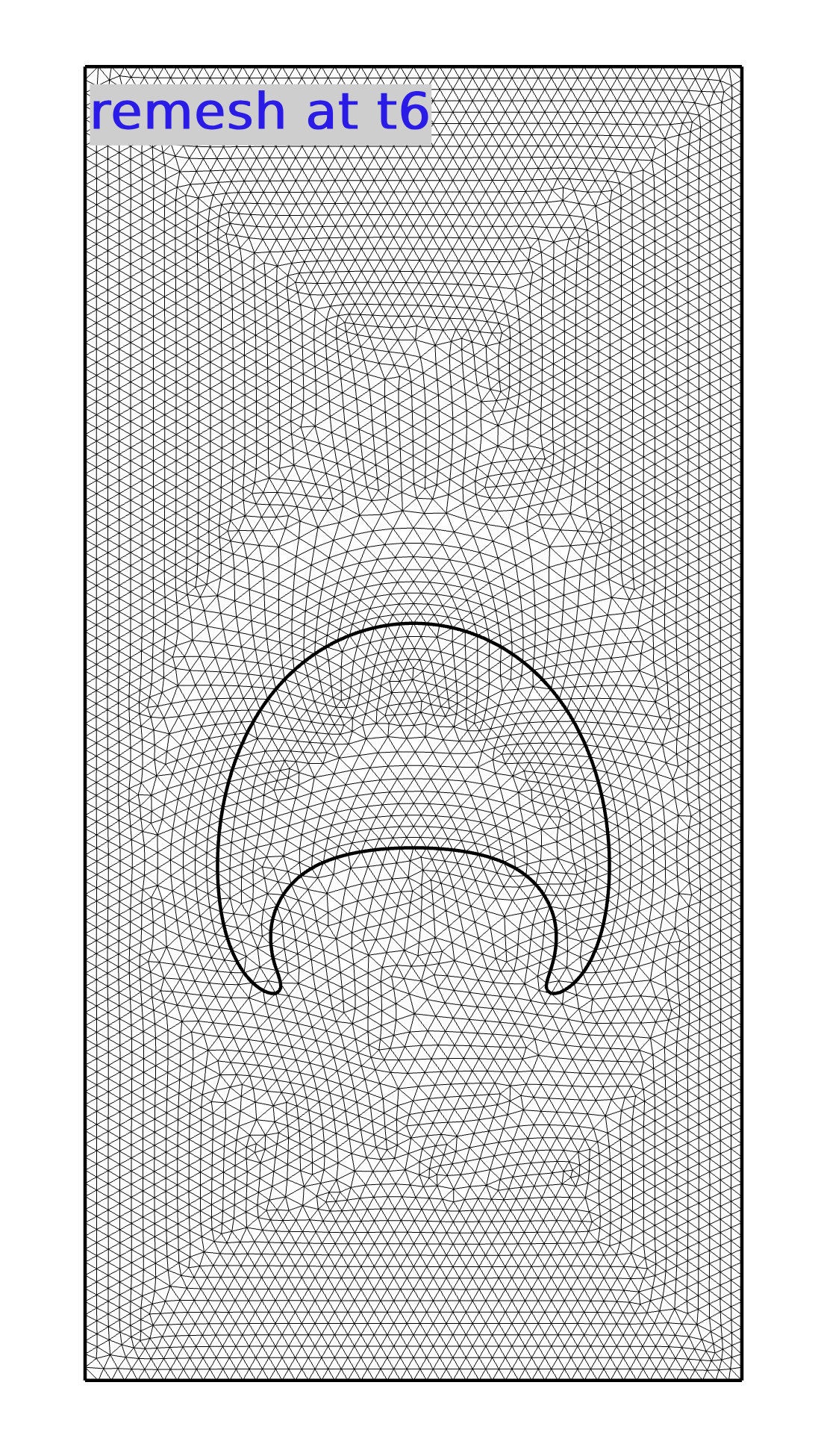}
}
\vspace{-10pt}
\caption{\it Meshes before and after a remeshing step for mesh $M_2$ with $k=2$.  }
\label{fig_remesh_compare_M2}
\end{figure}

\begin{figure}[htp!]
\vspace{-20pt}
\centerline{
\includegraphics[width=2.4in]{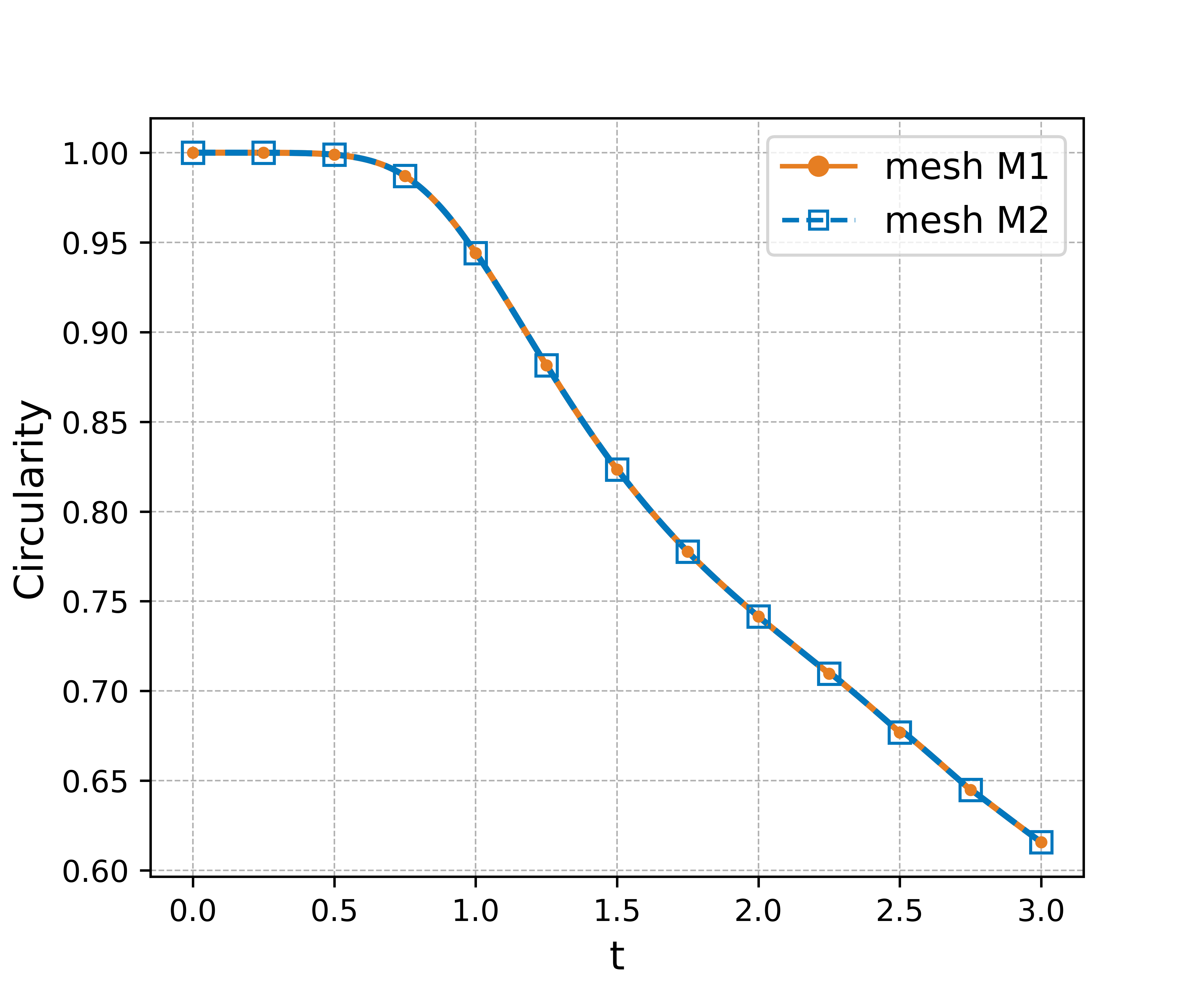}\hspace{-8pt}
\includegraphics[width=2.4in]{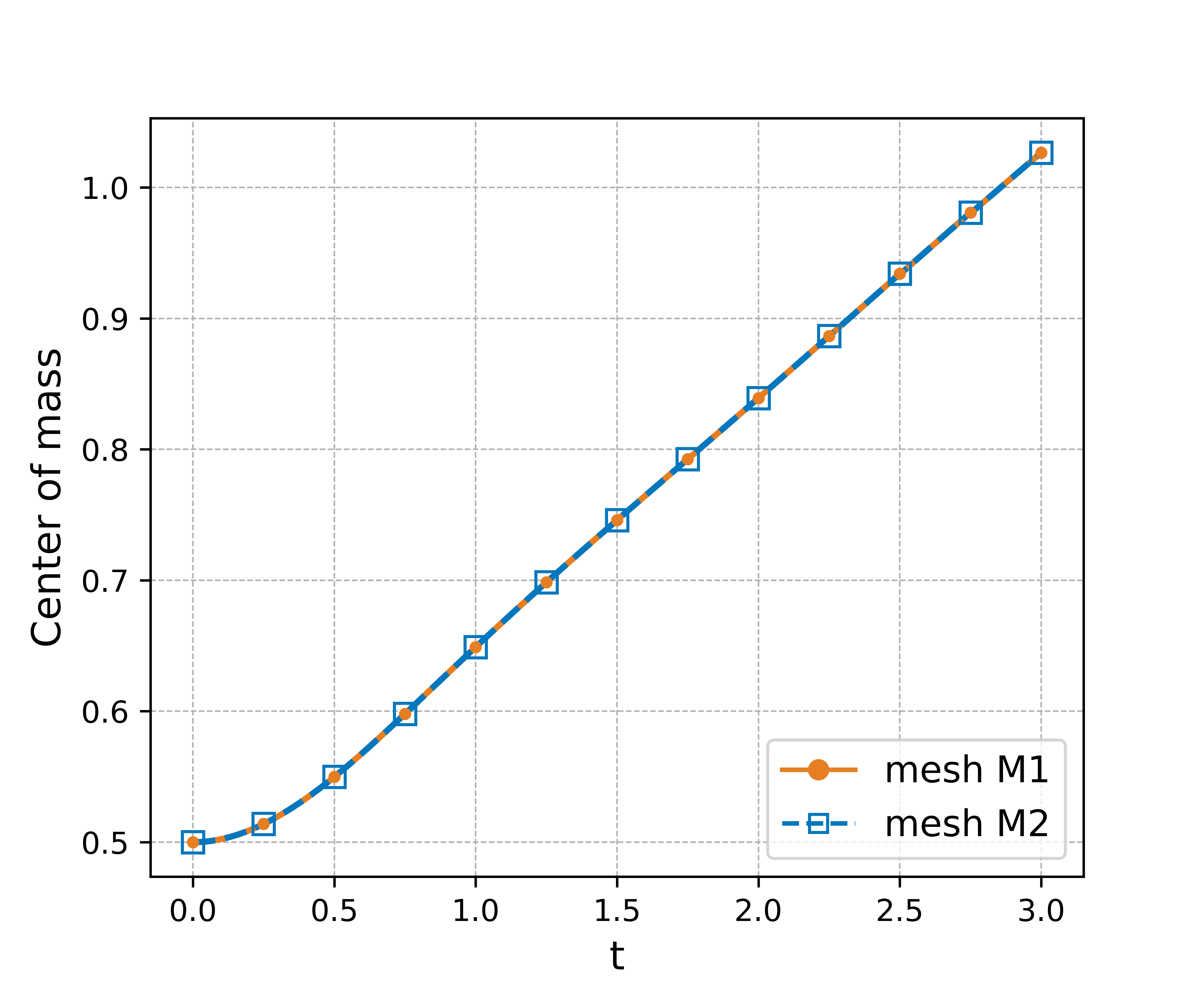}
}
\vspace{-3pt}
\centerline{
\includegraphics[width=2.4in]{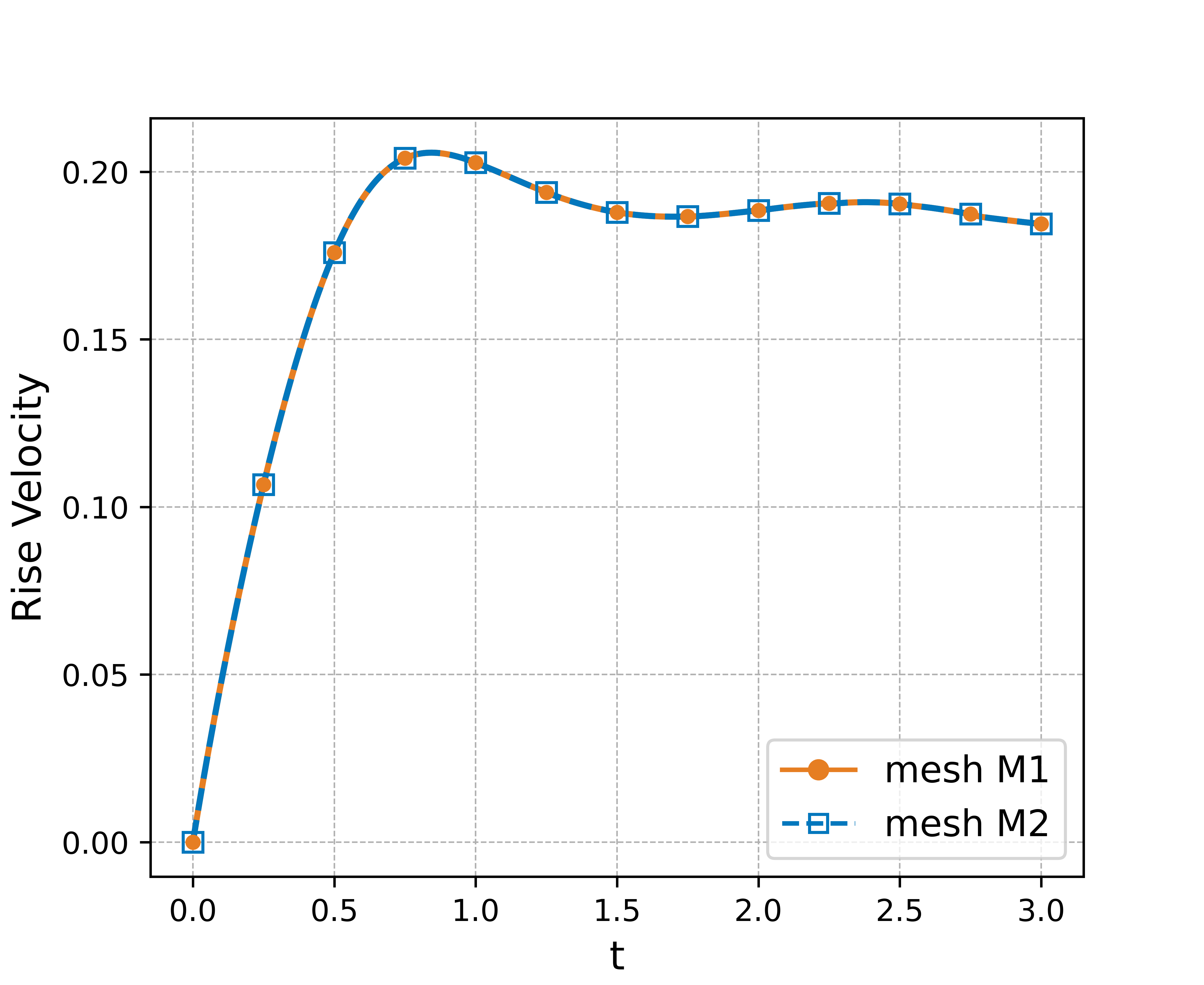}\hspace{-8pt}
\includegraphics[width=2.4in]{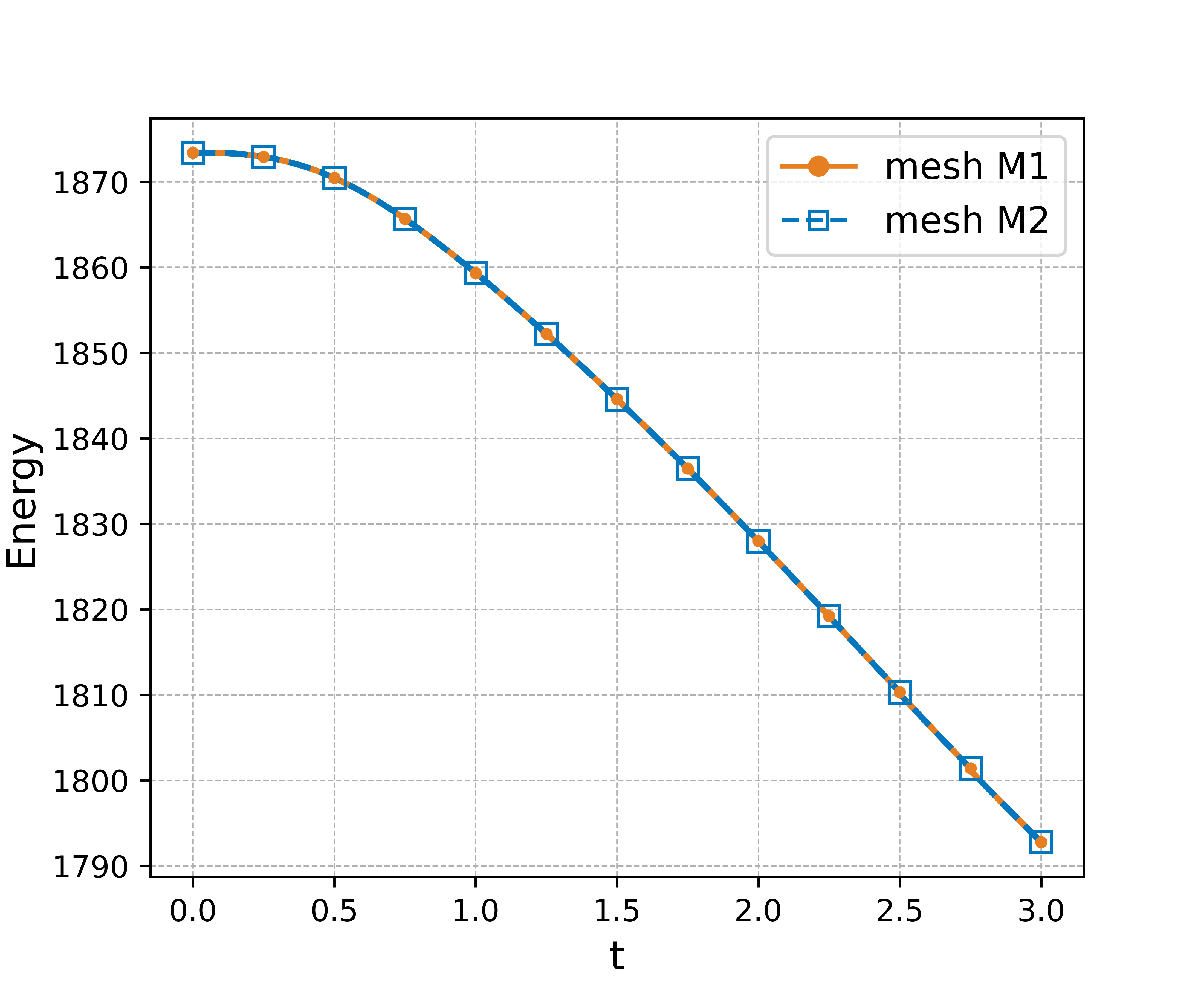}
}
\vspace{-6pt}
\caption{\it Benchmark quantities with respect to time}
\label{fig_benchmark_quantities}

\centerline{
\includegraphics[width=1.9in]{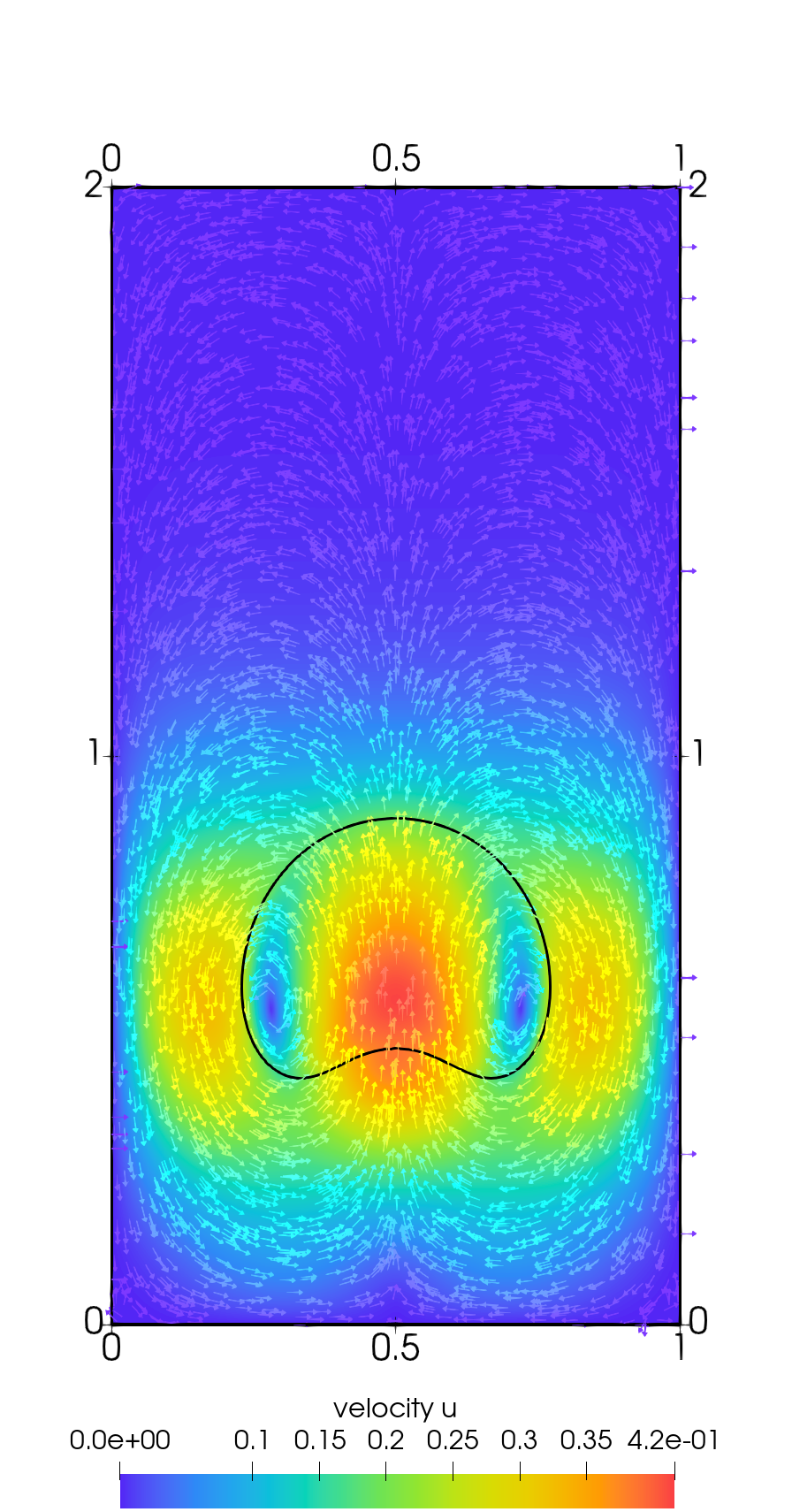}
\hspace{-20pt}
\includegraphics[width=1.9in]{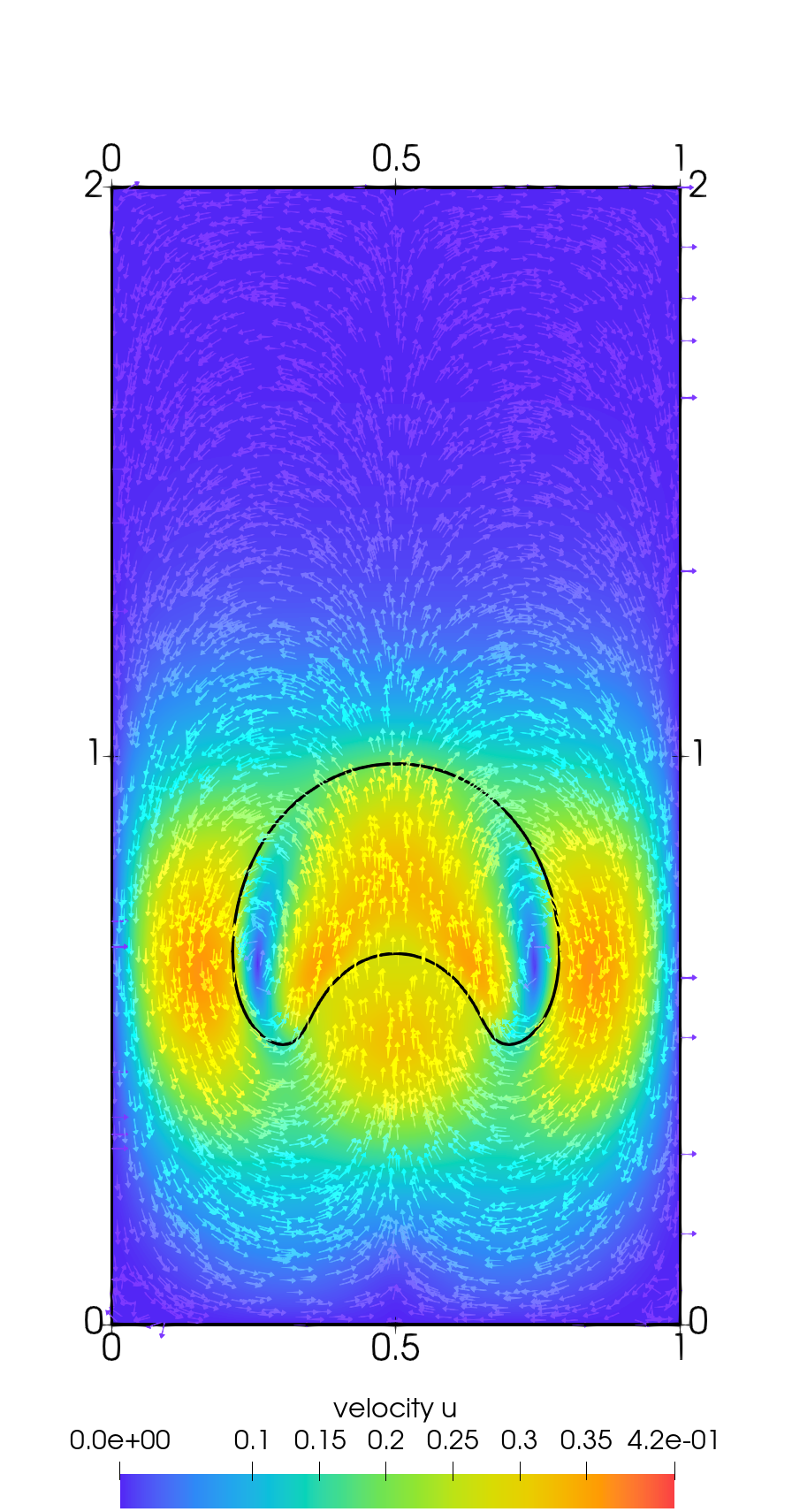}
\hspace{-20pt}
\includegraphics[width=1.9in]{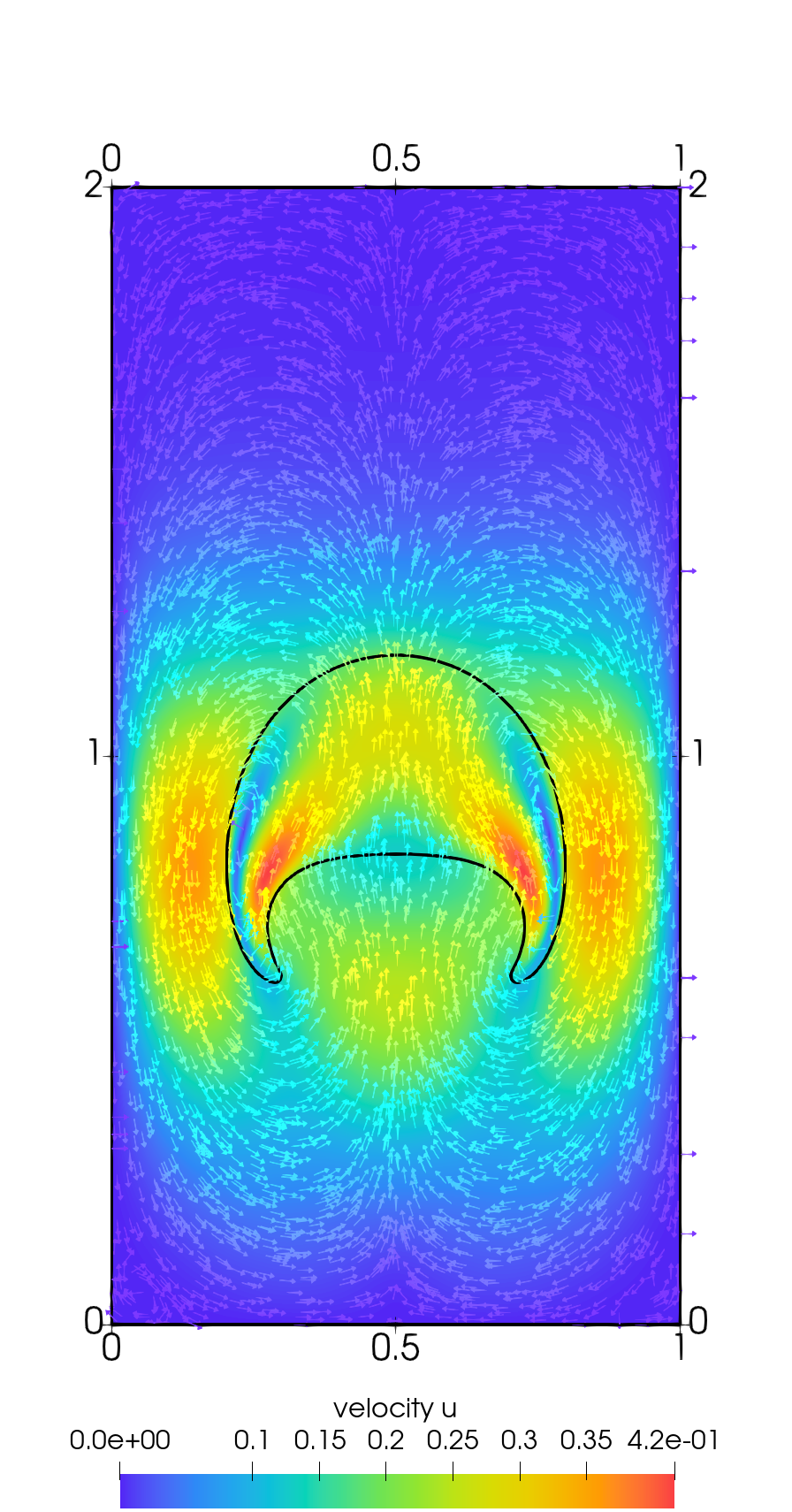}
}
\caption{\it Velocity. Left: $t_1$ = 1. Middle: $t_2$ = 1.5. Right: $t_3$ = 2.5.}
\label{fig_velocity}
\end{figure}

In the computation the mesh may degenerate and in this case one has to regenerate the mesh. 
{
In our computation we regenerate the mesh when the following angle criterion is satisfied:
\begin{align}
\min_{K\in \K_h[\x(t)]}
\min_{\alpha \in  \angle K} \le \frac{\pi}{18}, 
\end{align} 
where $\angle K$ denotes the set of interior angles of a triangle $K$.  The remeshing is performed using the software Gmsh. If the mesh is regenerated then we take the new mesh as the initial mesh to compute the above benchmark quantities.
Figure~\ref{fig_remesh_number} shows the number of remeshing steps over time for two different meshes, using a fixed stepsize $\tau = 1/200$. The meshes before and after  the 3rd and the 6th remeshing steps are shown in Figures~\ref{fig_remesh_compare_M1}-\ref{fig_remesh_compare_M2}.}
In Figure~\ref{fig_benchmark_quantities}, we present the numerical results of circularity, center of mass, rise velocity and energy on the two different meshes  with $\tau = 1/200$. 
From Figures \ref{fig_three_mesh}--\ref{fig_benchmark_quantities} we see that the method with the coarse mesh captures the shape of the interface very well. 
The geometry of the interface at $t = 1$, $1.5$ and $2.5$, as well as the distribution of interface-velocity with the mesh $M_2$, are presented in Figure~\ref{fig_velocity}.

{In summary, with} the numerical comparison of results between coarse and fine meshes in Figures \ref{fig_three_mesh}--\ref{fig_benchmark_quantities}, we have good reasons to trust the numerical simulation results for this benchmark example.

\section*{Acknolwedgement}
All authors contributed equally.
The work of B. Li was partially supported by National Natural Science Foundation of China (grant No. 12201418) and a fellowship award from the Research Grants Council of the Hong Kong Special Administrative Region, China (Project No. PolyU/RFS2324-5S03).
The work of S. Ma and W. Qiu was partially supported by a grant from the Research Grants Council of the Hong Kong Special Administrative Region, China (Project No. CityU 11300621).


\bibliographystyle{abbrv}
\bibliography{two-phase-NS-reference}

\end{document}